\documentclass[12pt]{amsart}

\usepackage{amssymb}
\usepackage{fancyhdr}
\usepackage{amsmath}
\usepackage{mathrsfs}
\usepackage{amsfonts}
\usepackage{enumitem}
\usepackage[all]{xy} 
\usepackage{xcolor}
\usepackage{diagbox}

\newcommand{\disk}{\ensuremath{\mathbb{D}} } 
\newcommand{\sphere}{\bar{\Bbb{C}}} 
\newcommand{\riem}{\Sigma}  
\renewcommand{\Bbb}[1]{\ensuremath{\mathbb{#1}}}

\newcommand{\qs}{\operatorname{QS}}

\newcommand{\kap}{\mathscr{K}}

\theoremstyle{plain}
        \newtheorem{theorem}{Theorem}[section]

        \newtheorem{corollary}[theorem]{Corollary}

\theoremstyle{definition}
        \newtheorem{definition}[theorem]{Definition}

\theoremstyle{remark}
    \newtheorem{remark}[theorem]{Remark}

\numberwithin{equation}{section} 
\numberwithin{figure}{section} 

\usepackage{fullpage}

\dedicatory{}
\begin{document}
\title[Weil-Petersson]{Weil-Petersson Teichm\"uller theory of surfaces of infinite conformal type}

\author{Eric Schippers}
\address{Eric Schippers \\ Department of Mathematics \\
University of Manitoba\\
Winnipeg, Manitoba \\  R3T 2N2 \\ Canada}
\email{eric\_schippers@umanitoba.ca}
\thanks{}

\author{Wolfgang Staubach}
\address{Wolfgang Staubach\\ Department of Mathematics\\
Uppsala University\\
Box 480\\ 751 06 Uppsala\\ Sweden}
\email{wulf@math.uu.se}


\date{\today}

\begin{abstract}
 Over the past two decades the theory of the Weil-Petersson metric has been extended to general Teichm\"uller spaces of infinite type, including for example the universal Teichm\"uller space. 
 In this paper we give a survey of the main results in the Weil-Petersson geometry of infinite-dimensional Teichm\"uller spaces. This includes the rigorous definition of complex Hilbert manifold structures, K\"ahler geometry and global analysis, and generalizations of the period mapping. We also discuss the motivations of the theory in representation theory and physics beginning in the 1980s. Some examples of the appearance of Weil-Petersson Teichm\"uller space in other fields such as fluid mechanics and two-dimensional conformal field theory are also provided.
\end{abstract}

\thanks{Eric Schippers is partially supported by the National Sciences and Engineering Research Council of Canada. }

\maketitle
\begin{section}{Introduction}
\begin{subsection}{History and Motivation}

{We start by recalling the historical origins of the Weil-Petersson scalar product and the Weil-Petersson metric. In analytic number theory and in relation to the theory of automorphic forms, H. Petersson \cite{Petersson} defined a scalar product on the space of cusp forms $S_k$ of weight $k$, as follows:\\
  First let \(\Gamma\) be a congruence subgroup of \(\mathrm{SL}(2,\mathbb{Z})\) and \({F}\) any fundamental domain
for \(\Gamma .\) Let \(\bar{\Gamma}=\Gamma /\{\pm \mathbf{I}\} .\) For the {\it{cusp forms}} \(f, g \in S_{k}(\Gamma)\) (i.e. modular forms of weight $k$ with zero constant coefficient in their Fourier series expansions), one defines the {\it {Petersson scalar product}} as 

\begin{equation}\label{Petersson scalar}
(f, g)=\frac{1}{\left[\mathrm{SL}_{2}(\mathbb{Z}): \bar{\Gamma}\right]} \iint_{{F}} (\mathrm{Im}(z))^{k} f(z) \overline{g(z)} \frac{d A_z}{(\mathrm{Im}(z))^{2}}.
\end{equation}
 Since the hyperbolic area element \(\frac{d A_z}{(\mathrm{Im}(z))^{2}}\) and \((\mathrm{Im}(z))^{k} f(z) \overline{g(z)}\) are both \(\Gamma\)-invariant, the integral is well-defined and independent of the choice
of fundamental domain. This scalar product  converges since both \(f\)
and \(g\) are cusp forms and \(|f(z)|=|g(z)|=O\left(e^{-2 \pi \mathrm{Im}(z)}\right)\) as \(\mathrm{Im}(z)\rightarrow \infty\). Moreover, only one of the forms $f$ or $g$ needs to be a cusp form for the Petersson product to converge.\\

         Using Petersson's scalar product, and in the context of Kodaira-Spencer's deformation theory of compact complex manifolds, A. Weil \cite{Weil1}, \cite{Weil2}, defined an inner product on the space of quadratic differentials on Riemann surfaces as follows. If a point of Teichmüller space is represented by a Riemann surface $\mathscr{R}$, then the cotangent space at that point can be identified with the space of quadratic differentials at $\mathscr{R}$. Since the Riemann surface has a natural hyperbolic metric (say in the case that $\mathscr{R}$ has negative Euler characteristic), then one can define a Hermitian inner product on the space of quadratic differentials by integrating over the Riemann surface. This induces a Hermitian inner product on the tangent space to each point of Teichmüller space, and hence a Riemannian metric.\\
         Explicitly, let
         \(T(\mathscr{R})\) be the Teichmüller space of compact Riemann surface $\mathscr{R}$ of genus \(g\), and let
$\mathscr{R}_{t} \in T(\mathscr{R})$ be a smooth path through \(\mathscr{R}:=\mathscr{R}_{0} \in T(\mathscr{R}) .\) The tangent vector \(\dot{\mathscr{R}}_{0}=\)
\(d \mathscr{R}_{t} /\left.d t\right|_{t=0}\) can be represented uniquely by a harmonic Beltrami differential

\begin{equation}
\mu=\rho^{-1} \bar{\varphi}
\end{equation}
where \(\rho\) is the hyperbolic metric and \(\varphi \in Q(\mathscr{R})\) is a holomorphic quadratic
differential. The {\it Weil-Petersson metric} (or WP metric for short) on \(T(\mathscr{R})\) is given by

\begin{equation}
\left\|\dot{\mathscr{R}}_{0}\right\|_{\mathrm{WP}}^{2}=\|\mu\|_{\mathrm{WP}}^{2}=\int_{\mathscr{R}} |\mu|^{2} \rho=\int_{\mathscr{R}} |\varphi|^{2} \rho^{-1},
\end{equation}
which could be compared with $(f,f)$ using the Petersson scalar product \eqref{Petersson scalar}.}\\

Now let us very briefly sketch Thurston's involvement in Teichm\"uller theory and the theory of Weil-Petersson metric. For a much more complete account of the Thurston's work in this as well as in other contexts, the reader is encouraged to read the essay by K. Ohshika and A. Papadopoulos \cite{Ohshikapapa} or consult the book by A. Fathi, F. Laudenbach, and V. Po\'enaru \cite{Fathi et al}.\\

William Thurston's interest in Teichm\"uller theory stemmed from his interest in the theory of surfaces and deformation theory of Riemannian metrics on surfaces. His motivation was to develop a Teichm\"uller-space theory which was purely based on hyperbolic geometry and did not rely on Teichm\"uller-Ahlfors-Bers' quasiconformal approach and the theory of quadratic differentials\footnote{Indeed since every Teichm\"uller equivalence class contains infinitely many quasiconformal maps, there is much redundant information in the quasiconformal formulation. This motivates the consideration of extremal quasiconformal maps, though as is well known uniqueness fails in general.}, see e.g. \cite{Ohshikapapa}. In accordance to his deformation theory-approach, given a surface $\mathscr{R}$ with the Euler characteristic $\chi(\mathscr{R}) < 0,$ one can always endow $\mathscr{R}$ with Riemannian metrics with curvature $-1.$ Then one says that two Riemannian metrics $g_1$ and $g_2$ are isotopic if $g_1$ stems from $g_2$ via a diffeomorphism isotopic to the identity. Now if $S$ is compact without boundary, the Teichm\"uller space $T(\mathscr{R})$ can be defined as the set of isotopy classes of all Riemannian metrics of curvature $-1$. If $\mathscr{R}$ has boundary, $T(\mathscr{R})$ is the set of isotopy classes of Riemannian metrics of curvature $-1$ on the interior of $\mathscr{R}$ which are modelled on the thin end of a pseudosphere near a boundary component. Thurston referred to such a metric as a hyperbolic metric.\\ 

In connection to his utilization of Teichm\"uller theory in theory of surfaces, one of his remarkable achievements was the proof of the the so-called Nielsen-Thurston classification theorem \cite{Thur2}. This result states that every homeomorphism $f$ from  a compact surface $S$ to itself is isotopic to a homeomorphism $F$ with at least one of the following properties:
\begin{itemize}
    \item $F$ has finite order (i.e. some power of $F$ is the identity),
    \item $F$ preserves some finite union of disjoint simple closed curves on $\mathscr{R}$,
    \item $F$ is pseudo-Anosov (a dynamical term introduced by Thurston, see e.g. \cite {Thur2}).
\end{itemize}

To see how Teichm\"uller theory enters the picture, Thurston considered a closed, connected, orientable surface \(\mathscr{R}\) of genus \(g \geq 2\) . Then he realized that 
T\(\left(\mathscr{R}\right)\), which is by a result of Fricke the space of hyperbolic metrics on \(\mathscr{R}\) up to isotopy (and
is an open ball of dimension \(6 g-6\)), has
a compactification that is homeomorphic to a closed ball. Furthermore the boundary of this ball is
\(\mathrm{PMF}(\mathscr{R})\), i.e. the space of projective classes of measured foliations on \(\mathscr{R}\). Moreover, the mapping class group of $\mathscr{R}$, denoted by
Mod\(\left(\mathscr{R}\right)\), acts on this closed ball. Then, applying the
Brouwer fixed point theorem, Thurston concluded that each element of Mod\(\left(\mathscr{R}\right)\)
fixes some point of T\(\left(\mathscr{R}\right) \cup \mathrm{PMF}(\mathscr{R})\). Finally the analysis of the various cases for the fixed point, yields the classification theorem. An analytic proof for the Nielsen-Thurston classification theorem, based on Teichm\"uller theory of extremal quasiconformal mappings was given by L. Bers in \cite{Bers1}. \\

Thurston also had his own 
approach to the Weil–Petersson metric, which has resulted in a clear understanding of the geometry on the Teichm\"uller space of compact surfaces. In this context, and once again motivated by hyperbolic geometry, he introduced a Riemannian metric on Teichmu\"uller space where
the scalar product of two tangent vectors at a point represented by a hyperbolic surface is the second derivative of the length of a uniformly distributed sequence of closed geodesics on the surface (or the Hessian of length of a random geodesic). Moreover he introduced in \cite{Thur1} the concept of {\it earthquakes} that generalize the Fenchel–Nielsen deformation operation of cutting a hyperbolic surface along a simple closed geodesic and gluing it back with a twist.

In \cite{WolpertThur}, S. Wolpert showed that Thurston’s metric is a multiple of the Weil–Petersson metric. This added yet another component to the understanding of  the picture for identification of
the geometry on Teichm\"uller space coming from the hyperbolic surface geometry, for compact surfaces.

{For a comprehensive overview of many aspects of Teichm\"uller theory (including the Weil-Petersson geometry) in the compact finite dimensional setting, see G. Schumacher's survey \cite{Schumachersurvey}}.
\end{subsection}
\begin{subsection}{Weil-Petersson geometry of surfaces of infinite conformal type} 
This paper is mainly about the Weil-Petersson geometry of the infinite-dimensional quasiconformal Teichm\"uller spaces. Here, we briefly explain its geometric motivations from a bird's-eye perspective. A more full picture of recent advances and their historical context is given in the remainder of the paper. We hope to give an idea of what has compelled many researchers - often looking into Teichm\"uller theory from the outside - to stir the embers of the Ahlfors-Bers theory.\\  

The story begins with an obstacle.
{On surfaces of infinite type such as the disk, the Weil-Petersson pairing does not converge for arbitrary directions in the Teichm\"uller space. So to investigate the Riemannian geometry of the Teichm\"uller spaces of such surfaces, one must remove this obstacle. This requires developing a refinement of Teichm\"uller space whose tangent spaces are Hilbert spaces with respect to the Weil-Petersson product. }\\

Since this is not easy, it is worth first reflecting: why should there be any interest in the infinite-dimensional Teichm\"uller spaces in the first place (let alone their Weil-Petersson refinements)?  After the Ahlfors-Bers period, it is probably safe to say that the focus moved away from quasiconformal Teichm\"uller theory, perhaps mainly because of the overwhelming force of the insights of Thurston. His largely hyperbolic geometric point of view of Riemann surfaces is in some sense Kleinian. That is, the universal cover is a homogeneous space: the quotient of all isometries of the disk by the stabilizer of a point.\\ 

The hyperbolic geometry of Riemann surfaces is rigid. Any local orientation-preserving isometry between two Riemann surfaces is a covering map (see J. H. Hubbard \cite[Proposition 3.3.4]{Hubbard_Book}; in the case of the disk, this reduces to the equality statement in Schwarz lemma). Thus hyperbolic geometry of a Riemann surface is essentially finite-dimensional. 

In contrast, locally conformal transformations need not extend globally. They are much less rigid than hyperbolic isometries - it requires a power series to specify a local biholomorphism fixing a point. Thus the local symmetry group of two-dimensional conformal geometry (complex analysis) is infinite-dimensional.   One therefore sees that the deformation theory of two-dimensional conformal geometry should also be infinite-dimensional. This is indeed the case: for example, the Teichm\"uller space of the disk is a Banach manifold.\footnote{  In fact, by Bers' trick, Teichm\"uller spaces - the deformation spaces - are identified with spaces of conformal symmetries.} \\ 

The local symmetry group also defies easy description as a Lie group. However, although the local symmetry group is not a Lie group according to the modern definition, it does fit S. Lie's original conception of Lie groups as locally defined symmetry groups of PDEs (see the introductions to J. F. Pommaret \cite{Pommaret} and A. Vinogradov \cite{Vinogradov}). In the case at hand, conformal transformations are the orientation-preserving symmetries of Laplace's equation. It is then no surprise that the outside interest in quasiconformal Teichm\"uller theory, especially the universal Teichm\"uller space, has come from areas in which this local point of view cannot be avoided or postponed. These areas include fluid mechanics, infinite-dimensional groups, geometric PDEs, and areas relating to physics, such as two-dimensional conformal field theory, vertex operator algebras, string theory, and geometric quantization.  In these fields the Ahlfors-Bers Teichm\"uller theory makes an appearance.\\  

{We cannot resist making the following comparison. The imposing genius of G. Gould has had such overwhelming effect on listeners, that for a long time some of his interpretations of J. S. Bach, e.g. of the \emph{Goldberg variations} and \emph{the art of the fugue}, have been considered the ultimate renditions of these pieces. Indeed, since the 1950's, Gould's way of rendering Bach has been viewed as a paradigm-shift in the way of playing Bach's keyboard music. However one need neither dismiss Gould in order to see the beauty of entirely different paradigms - nor abandon the other paradigms to appreciate his insights. In the above analogy, in place of Gould one could substitute Klein,  Ahlfors {and Bers}, Grothendieck, or Thurston, and in place of Bach one could substitute Riemann.  In hindsight the Weil-Petersson theory of infinite type surfaces has a simplicity which fits into many paradigms.}  
\\

Analytic problems are a barrier to fitting infinite-dimensional phenomena into the standard geometric paradigms. For example, infinite-dimensional groups might not be continuous (e.g. left multiplication is not continuous in the universal Teichm\"uller space). Similarly, representations necessary for creating vector bundles and connections require solution of problems in functional analysis.   In the case of Teichm\"uller theory, analytic advances have now removed this barrier. Returning to the two examples above, by a theorem of L. Takhtajan and L. P. Teo, the Weil-Petersson universal Teichm\"uller space is a topological group. Also, thanks to a theorem of S. Vodopyanov and S. Nag and D. Sullivan, one obtains a faithful representation of the universal Teichm\"uller space as bounded symplectomorphisms of the Sobolev $H^{1/2}$ space of the unit circle.

We have attempted in this paper to give a survey of some of these advances, and provide a road map of the new Weil-Petersson analytic theory. Since it is still somewhat new, we have emphasized giving a literature review over giving the quickest possible introduction to the results.  Throughout the paper, we emphasize how new geometric and algebraic theorems motivate the analytic theory\footnote{But this also goes the other way - the analysis has informed the geometry and algebra.}.    \\

Most importantly, the Weil-Petersson theory makes the K\"ahler geometry of the moduli spaces of infinite conformal type possible, and thereby also enriches our understanding of the geometry, analysis, representation theory, and physics.  

For K\"ahler geometry one needs a complex manifold, and so much of the development of the Weil-Petersson Teichm\"uller theory has concerned the proof of the existence and compatibility of various complex structures. In particular we outline the results on three complex structures: analogues of two classical complex structures obtained from the Bers embedding and harmonic Beltrami differentials, and a third complex structure. This third complex structure is in some sense a modification of Bers' idea, and is motivated by two-dimensional conformal field theory.  A fourth complex structure arises from generalizations of the classical period map. 

The study of complex structures in Teichm\"uller theory of course has a long history. It is striking how natural the results seem, no matter what angle they are viewed from. Indeed one feels that many of the results, e.g. the generalized period maps and their geometry, could be motivated to mathematicians a century ago. \\

\textbf{Acknowledgements.} {We are grateful to Yi-Zhi Huang for many valuable discussions and for encouraging our work on the moduli spaces, and of course to David Radnell who has been our partner in much of this work.}

\end{subsection}
\begin{subsection}{Outline} 
  We have tried to give a faithful account of the development of the infinite-dimensional Weil-Petersson theory. As far as we know no such account has yet been given.  The subject was initiated in the 80s, but after some foundational results were established in the 2000s, the subject began to lift off. Now results are coming in quickly, and no doubt by the time this appears there will be new developments.

  Whenever possible we of course attempted to use familiar notation, but unfortunately the notation varies wildly in textbooks and papers. To maintain comprehensibility we were forced to be consistent, and it was therefore impossible to adhere to various authors' choices. Needless to say, all theorems are invariant under change of notation.   

  In Section \ref{se:preliminaries} we establish notation and basic definitions. Since a large portion of the paper involves the existence and equivalence of complex structures on the Weil-Petersson Teichm\"uller space, to set the stage we review three models of the usual Teichm\"uller space and give their associated constructions of the complex structure in Sections \ref{se:three_models} and \ref{se:classical_complex_structure}.  Two of these models are classical ones given by the Bers embedding and by harmonic Beltrami differentials. The third is a fiber model of D. Radnell and E. Schippers \cite{RadnellSchippers_fiber} which is motivated by conformal field theory. This model is used in some of the constructions and relates to some of the applications of Weil-Petersson Teichm\"uller theory.  The section concludes with an overview the differences between the classical $L^\infty$ Teichm\"uller theory and the $L^2$ Weil-Petersson Teichm\"uller theory. 

  In Section \ref{se:WP_Universal_Teich} we outline the main results regarding the Weil-Petersson universal Teichm\"uller space of G. Cui and L.A. Takhtajan-L.P. Teo. A description of polarizations and generalizations of the period map into the Siegel disk, is given in Section  \ref{se:period}. Because of their importance in the motivation of the Weil-Petersson theory, we include a historical outline.  Section \ref{birds eye view} briefly outlines some other refinements of Teichm\"uller spaces for the sake of completeness.  Section \ref{se:WP_higher_genus} outlines the Weil-Petersson theory of general surfaces, due to M. Yanagishita and Radnell-Schippers-Staubach.
  
  Section \ref{se:Kahler_and_global} give results on the K\"ahler theory and global analysis of the Weil-Petersson metric. 
  Section \ref{se:in_the_wild} gives a brief account of some applications to physics, fluid mechanics, and stochastic Loewner evolution.

\end{subsection}
\end{section}
\begin{section}{Preliminaries} \label{se:preliminaries}
 
\begin{subsection}{Surfaces, borders, and lifts}
   Let $\mathbb{C}$ and $\sphere$ denote the complex plane and the Riemann sphere respectively. Let $\disk_+ = \{z : |z|<1 \}$ and $\disk_- = \{ z: |z|>1 \} \cup \{\infty \}$ denote the unit disk at $0$ and $\infty$ respectively. 
   We will denote the boundary of the disk by $\mathbb{S}^1$, which we assume to be oriented positively with respect to $\disk_+$.\\
  
  In this paper we consider mostly Riemann surfaces $\riem$ whose universal cover is the disk, that is $\riem=\disk_-/G$ where $G$ is a Fuchsian group (here the quotient map is a covering map so we do not allow our Fuchsian groups to have fixed points). 
  
  The limit set of a Fuchsian group $G$ acting on $\disk_+$ or $\disk_-$ is either the entire boundary $\mathbb{S}^1$ or a nowhere dense subset of $\mathbb{S}^1$; in which case it is called of first or second kind. We will be concerned in this paper with surfaces generated by Fuchsian groups of the second kind.  These are bordered Riemann surfaces (see \cite{Lehto,Ahlfors_Sario}). If $\riem \cong \disk_- /G$ for a Fuchsian group $G$ of the second kind, and we let $B$ be the complement of the limit set of $G$ in $\mathbb{S}^1$, then $\disk_-/G \cup B/G$ is a bordered Riemann surface, with $\disk_-/G \cong \riem$ and $B/G$ the border. Thus every time we mention bordered Riemann surfaces, we are referring to surfaces that are obtained exactly in this way.\\Such a surface also has a double $\riem^d \cong \sphere/G$ so that $\disk_+ /G$ identified with the conjugate surface $\riem^*$.This has an anti-holomorphic involution 
  \[ i: \riem^d \rightarrow \riem^d \]
  which takes $\riem$ to $\riem^*$ and vice versa, which preserves the border. 

  Surfaces of genus $g$ whose border has $n$ connected components, each of which is homeomorphic to $\mathbb{S}^1$, will be referred to as type $(g,n)$ bordered surfaces. We will always assume $n \geq 1$ when using this terminology. For such surfaces, the double is a compact Riemann surface of genus $2g+n-1$, and the border $B/G$ can be naturally viewed as an analytic curve $\partial \riem = \partial \riem^*$ in the double.
\end{subsection}
\begin{subsection}{Differentials}  
 {Given a Riemann surface $\riem$, an $(r,s)$-differential on $\riem$, $(r,s)\in \mathbb{Z}\times \mathbb{Z},$} is a quantity given in local coordinates by $h(z)\, dz^r \odot d\bar{z}^s$, which transforms under a holomorphic change of coordinates $z=g(w)$ via 
 \[  h(z)\, dz^r \odot d\bar{z}^s = \tilde{h}(w)\, dw^r \odot d\bar{w}^s, \ \ \ \ \tilde{h}(w) =  h(g(w)) g'(w)^r \overline{g'(w)}^s.   \]
 Here $\odot$ denotes the symmetric tensor product. We shall denote the space of measurable $(r,s)$-differentials on $\riem$ by $\mathscr{D}_{r,s}(\riem).$
 The regularity will be imposed when the need arises. 
 As is customary we will write the local coordinates form as 
 \[   h(z) dz^r d\bar{z}^s.   \]

 For example, quadratic differentials are $(2,0)$-differentials, and Beltrami differentials are $(-1,1)$ differentials. It is easily seen that products of differentials are well-defined, and that the product of an $(r_1,s_1)$- and an $(r_2,s_2)$- differential is an $(r_1+r_2,s_1+s_2)$-differential.  \\

 The hyperbolic metric on $\riem$ can be viewed as a $(1,1)$-differential, which we will denote by $\rho_\riem$, often written locally as 
 \[ \rho_\riem=  h(z) |dz|^2 \] 
 where $h(z)$ is a strictly positive function. Here it must be remembered that $|dz|^2$ stands for $dz\odot d\bar{z}$.   For example, we have that the hyperbolic metric on $\disk_+$ and $\disk_-$ are 
 \[  \rho_{\disk_+} = \frac{d\bar{z} \odot dz }{(1-|z|)^2}, \ \ \   \rho_{\disk_-} = \frac{d\bar{z} \odot dz }{(|z|^2-1)^2} \]
 respectively.  The metric induces a hyperbolic area measure given locally by
 \[  dA_{\mathrm{hyp}} = h(z) \frac{d\bar{z} \wedge dz}{2i}   \]
 Some sources also write this as $h(z) |dz|^2$.

 The line element $\lambda_\riem = \sqrt{h(z)} |dz|$ is a convenient object which transforms under change of coordinates $z=g(w)$ as a density $\sqrt{h(z)}|dz|=\sqrt{h(g(w))}||g'(w)||dw|$. 
 The density can be used to construct an invariant $L^2$ and $L^\infty$ norm.\\

 For $\alpha \in \mathscr{D}_{r,s}(\Sigma)$ locally given by $g(z)\, dz^r d\bar{z}^s$, it is easily checked that 
 $|\alpha / \lambda_\riem^{r+s}|$ 
 is a well-defined global function which agrees with the local expressions $\sqrt{h(z)}^{-r-s} |g(z)|$. Denote its essential supremum by $\| \alpha \|_\infty$.  
 We then denote 
 \[  L^\infty_{r,s}(\riem) = \{ \alpha \in \mathscr{D}_{r,s}(\Sigma) : \| \alpha \|_\infty < \infty \} \]
 Similarly we obtain the $L^2$ norm
  \[   \| \alpha \|_2 = \iint_{\riem} | \alpha|^2 \lambda^{-2r-2s} dA_{\mathrm{hyp}}     \]
 and the associated Hermitian inner product on $L^2$ differentials, and denote
 \[  L^2_{r,s}(\riem) = \{ \alpha \in \mathscr{D}_{r,s}(\Sigma) :  \| \alpha \|_{2} < \infty \}. \] 
 For example 
\[  L^2_{-1,1}(\disk_-) := \left\{  \mu \in \mathscr{D}_{-1,1}(\disk_-) \,:\, \iint_{\disk_-} \frac{|\mu(z)|^2}{(|z|^2-1)^2} d\mathrm{A}_z < \infty \right\}   \]
where $d\mathrm{A}_z$ denotes Euclidean area measure.

 If $\riem = \mathbb{D}_\pm /G$ where $G$ is a Fuchsian group, then any $(r,s)$-differential $\alpha$ has a unique lift $\hat{\alpha} = a(z) dz^r d\bar{z}^s$ to $\disk_\pm$, satisfying the invariance 
 \begin{equation} \label{eq:defn_differential_invariance}
   a \circ g \cdot (g')^r \cdot (\overline{g'})^s = a \ \ \ \text{for all} \ \ g \in G. 
 \end{equation}
 Denote by $L^\infty_{r,s}(\disk_\pm,G)$ and $L^2_{r,s}(\disk_\pm,G)$ the spaces of differentials which are invariant in the sense of (\ref{eq:defn_differential_invariance}) which are $L^\infty$ and $L^2$ respectively {\it on a fundamental domain $D \subset \disk_\pm$ of $G$}. 
 Note that an element of $L^\infty_{r,s}(\disk_\pm)$ is essentially bounded on $\disk_\pm$ by the invariance. On the other hand, and element of  $L^2_{r,s}(\disk_\pm,G)$ need not be in $L^2_{r,s}(\disk_\pm)$.
 
 We then have a natural isomorphism (directly induced by the covering map
 \begin{align} \label{eq:differential_lift_direct_isomorphisms}
  L^\infty_{r,s}(\riem) & \cong L^\infty_{r,s}(\disk_-,G) \nonumber \\
  L^2_{r,s}(\riem) & \cong L^2_{r,s}(\disk_-,G)
 \end{align}
 and 
 \begin{align} \label{eq:differential_lift_direct_isomorphisms_double}
  L^\infty_{r,s}(\riem^*) & \cong L^\infty_{r,s}(\disk_+,G) \nonumber \\
  L^2_{r,s}(\riem^*) & \cong L^2_{r,s}(\disk_+,G)
 \end{align}
 \begin{remark} \label{re: } We adhere to the convention that $\riem$ is covered by $\disk_-$ and its double $\riem^*$ is covered by $\disk_+$. 
 \end{remark}
 
 The norms and inner products are conformally invariant, as are their straightforward generalizations to $L^p$ spaces. These norms are found throughout complex analysis, Teichm\"uller theory, and number theory with various names and notations.  

 Finally, we will use the notation $A^\infty_r(\riem)$ and $A^2_r(\riem)$ for the subspaces of holomorphic elements of $L^\infty_{r,0}(\riem)$ and $L^2_{r,0}(\riem)$ respectively.   For example 
 \[  A_2^\infty(\disk_+) =  \left\{ Q(z)dz^2  \ \text{ holomorphic on } \disk_+ \,:\,
    \sup_{z \in \disk_+} (1-|z|^2)^2|Q(z)| < \infty \right\}   \]
    is the familiar space of Nehari-bounded quadratic differentials. Also
  \[  A_2^2(\disk_+) :=  \left\{ Q(z)dz^2  \ \text{ holomorphic on } \disk_+ \,:\,
   \iint_{\disk_+} (1-|z|^2)^2|Q(z)|^2 d\mathrm{A}_z < \infty \right\}  \]
  will play an important role in the Weil-Petersson theory.
\end{subsection}
\begin{subsection}{Deformations: quasiconformal maps and quasisymmetries}
 A quasiconformal map $f:\riem_1 \rightarrow \riem_2$ between Riemann surfaces (sometimes abbreviated as $f\in \mathrm{QC}(\riem_1, \riem_2)$), is an orientation-preserving homeomorphism such that the $(-1,1)$-differential  $\mu(f) = \overline{\partial} f/\partial f$ satisfies
 \[   \| \mu(f) \|_\infty < 1.   \]
 Denote the unit ball in $L^\infty_{-1,1}(\riem)$ by $L^\infty_{-1,1}(\riem)_1$. Quasiconformal maps are the generic deformations of Riemann surfaces which do not preserve the complex structures. 
  
 An orientation-preserving homeomorphism $h$ of $\mathbb{S}^1$ is called a \emph{quasisymmetric mapping}, iff there is a constant $k>0$, such that for every $\alpha$ and every $\beta$ not equal to a multiple of $2\pi$, the inequality \[  \frac{1}{k} \leq \left| \frac{h(e^{i(\alpha+\beta)})-h(e^{i\alpha})}{h(e^{i\alpha})-h(e^{i(\alpha-\beta)})} \right|
    \leq k \]
holds. Let $\qs(\mathbb{S}^1)$ denote the set of quasisymmetric maps from $\mathbb{S}^1$ to $\mathbb{S}^1$.  Every quasiconformal map $\Phi:\disk_- \rightarrow \disk_-$ has a continuous extension to $\mathbb{S}^1$, and this extension is a quasisymmetry. Conversely every quasisymmetry is boundary values of a (highly non-unique) quasiconformal map. Let $\text{M\"ob}(\mathbb{S}^1)$ denote the M\"obius transformations which preserve the circle and its orientation.\\ 

Let $\riem$ be a Riemann surface, possessing a border homeomorphic to $\mathbb{S}^1$ which we denote by $C$. Then there is a biholomorphism $\psi:A \rightarrow \mathbb{A}_r$ where $A$ is a doubly-connected region in $\riem$ bounded on one side by $C$ and $\mathbb{A}_r = \{z: r <|z|<1 \}$ is an annulus in the plane. By Carath\'eodory's theorem $\psi$ extends continuously to a homeomorphism from $C$ to one of the boundary curves of $\mathbb{A}_r$, which we assume to be $\mathbb{S}^1$.   
\begin{definition} \label{de:quasisymmetries_surfaces}
We say that a map $\phi:C \rightarrow C_1$ is a quasisymmetry if $\psi_1 \circ \phi \circ \psi^{-1}$ is a quasisymmetry of $\mathbb{S}^1$ for some choice of collar chart $\psi_1$ and $\psi$ of $C_1$ and $C$ respectively. The set of such maps is denoted by $\mathrm{QS}(C,C_1)$. 

Given a bordered Riemann surfaces $\riem$, $\riem_1$ of type $(g,n)$, let $\mathrm{QS}(\partial \riem,\partial \riem_1)$ denote the set of bijections $\phi:\partial \riem \rightarrow \partial \riem_1$ such that the restriction of $\phi$ to each connected component is a quasisymmetry in the sense above. 
\end{definition}

It is not hard to show that if $\psi_1 \circ \phi \circ \psi^{-1}$ is a quasisymmetry for a choice of collar charts $\psi_1$ and $\psi$, then it is a quasisymmetry for all choices. 
\end{subsection} 
\end{section}
\begin{section}{Three models of Teichm\"uller space}\label{se:three_models}  
 In this section, we give three models of Teichm\"uller space: in terms of quasiconformal deformations of the surface itself; in terms of conformal deformations on the mirror image; and in terms of conformal deformations of ``caps''. The first two pictures are classical, and the last picture originates in conformal field theory. The connection of this model to Teichm\"uller theory is due to Radnell-Schippers. 

 Within these models, one may describe the Teichm\"uller space either directly on the surface itself or in terms of objects on the cover which are invariant under the Fuchsian group. We refer to these the {\bf direct picture} and {\bf lifted picture} respectively.

 \begin{subsection}{Definition of Teichm\"uller space} 
   In the direct picture, we fix a Riemann surface $\riem$, covered by the disk via a Fuchsian group of the second kind. Let $\partial \riem$ denote the border (which is possibly empty). 
   \begin{definition}[Teichm\"uller equivalence]
   Given two quasiconformal maps $F_1 :\riem \rightarrow \riem_1$ and $F_2: \riem \rightarrow \riem_2$, we say that $(\riem,F_1,\riem_1)$ and $(\riem,F_2,\riem_2)$ are Teichm\"uller equivalent if there is a biholomorphism $\sigma:\riem_1 \rightarrow \riem_2$ such that $F_2^{-1} \circ \sigma \circ F_1$ is homotopic to the identity rel boundary. Here homotopy rel boundary means that $F_2^{-1} \circ \sigma \circ F_1$ is the identity on the boundary, and this holds throughout the homotopy.
   \end{definition}

   The Teichm\"uller space is then 
   \[ T(\riem) = \{ (\riem,F_1,\riem_1) \}/\sim  \]
   where $\sim$ is the Teichm\"uller equivalence relation.  
   Given a representative $(\riem,F_1,\riem_1)$, the differential $\overline{\partial} F_1/\partial F_1$ is in $L^\infty_{-1,1}(\riem)_1$.  Passing the equivalence relation over to $L^\infty_{-1,1}(\riem)_1$ we obtain 
   \[   T(\riem) \cong L^\infty_{-1,1}(\riem)_1 /\sim. \]

   In the lifted model, assume that $\riem = \disk_- /G$ where $G$ is a Fuchsian group of the second kind. 
   $F_1,F_2 :\riem \rightarrow \riem$ are homotopic rel boundary if and only if, for the suitably normalized lifts $w_1, w_2$ to $\disk_-$, $w_1 = w_2$ on $\mathbb{S}^1$   
 \cite{Lehto}. This leads to the following equivalent definition of Teichm\"uller space.
 \begin{theorem} \label{th:Fuchsian_qcmap_model} The Teichm\"uller space $T(\riem)$ is in one-to-one correspondence with 
    \[ T(G):= \{ w\in \mathrm{QC}(\disk_- ,\disk_-) : w \ \text{fixes $-1$, $i$, $1$}, \  w \circ g \circ w^{-1} \in \emph{M\"ob}(\mathbb{S}^1) \ \forall \ g \in G \}/\sim \]
    where $w_1 \sim w_2$ if $w_1 = w_2$ on $\mathbb{S}^1$.
 \end{theorem}
 By (\ref{eq:differential_lift_direct_isomorphisms}) we have that 
 \[ T(G) \cong L^\infty(\disk_-,G)/\sim. \]
 Note that we can identify $T(\disk_-)$ with $T(\{ 1 \})$. 
 
 There is also an important variant known as the {\it quasi-Fuchsian model} arising from the so-called Bers' trick (see Subsection \ref{Berstrick} ahead), which we will not treat here.  

 Teichm\"uller space has a topology induced by the Teichm\"uller metric. 
 Recall that given the points $p_1$ and $p_2$ in \(T(\riem)\), the Teichmüller metric
is defined by
\begin{equation}\label{Teichmetric}
    d_{T(\riem)}(p_1, p_2) = \inf \mathrm{arctanh} \, \| \mu_{F_2 \circ F_1^{-1}} \|_\infty  .
\end{equation} 
where the infimum is taken over all representatives $[\riem,F_k,\riem_k]$ of $p_k$ for $k=1,2$.

This can be also formulated in the lifted picture. 

For later use, we define here the modular group. 
Let $Q(\riem):=\mathrm{QC}(\riem, \riem)$ and $Q_0(\riem)$ denote the sub-collection of $Q(\riem)$ which are homotopic to the identity rel boundary. This is a normal subgroup with respect to composition. The Teichm\"uller modular group is the quotient 
\[ \mathrm{Mod}(\riem) = Q(\riem)/Q_0(\riem) \]
which anti-acts on the Teichm\"uller space via
\[  [\rho][\riem,F,\riem_1] \mapsto [\riem,F \circ \rho,\riem_1].   \]
It is of course well-known that the Riemann moduli space can be identified with $T(\riem)/\mathrm{Mod}(\riem)$.

The following subgroups of the modular group play an important role ahead.  $\mathrm{ModI}(\riem)$ consists of those equivalence classes of $\mathrm{Mod}(\riem)$ whose representatives are the identity on the boundary.  This group is generated by two types of quasiconformal maps. A basis for the homotopy can be given by $\partial_1 \riem,\ldots,\partial_n\riem,a_1,\ldots,a_g,b_1,\ldots,b_g$ where $a_1,\ldots,a_g,b_1,\ldots,b_g$ are a basis for the homotopy of the compact surface obtained by sewing on caps to obtain a compact surface of genus $g$. Then $\mathrm{ModI}(\riem)$ is generated by integer twists around these curves. 
The group generated by integer twists around $\partial_1 \riem,\ldots,\partial_n \riem$ is denoted by $\mathrm{DB}(\riem)$.  The group generated by integer twists around non-trivial simple closed interior curves $a_1,\ldots,b_g$ is denoted by $\mathrm{DI}(\riem)$. 

Now by the Teichm\"uller equivalence relation, there is a lot of redundancy in the quasiconformal deformation. Which Beltrami differentials are ``essential''?  There is no global answer to this question for general surfaces, but there is a local one.  
Given a Riemann surface $\riem$, define
\begin{equation} \label{eq:harmonic_Beltrami_definition}
  \Omega_{-1,1}(\riem) = \{ \mu \in \mathscr{D}_{-1,1}(\riem) \,: \,  \rho_\riem \mu \in \overline{A^\infty_2(\riem)}.  \}  
\end{equation}
or, in the lifted picture, 
\begin{equation*}  
  \Omega_{-1,1}(\disk_-,G) = \{ \mu \in \mathscr{D}_{-1,1}(\disk_-) \,: \,  \rho_\riem \mu \in \overline{A^\infty_2(\disk_-,G)}.  \}  
\end{equation*}
These differentials are transverse to the Teichm\"uller equivalence relation at the identity. A sufficiently small ball stays within the Teichm\"uller space. This forms the basis for the deformation model fo the complex structure and tangent space, as we will see ahead.  
 \end{subsection}
 \begin{subsection}{Bers embedding model}\label{Berstrick}
   The underlying idea of this model comes from what is called Bers' trick. Given a Riemann surface, rather than reflect the quasiconformal map to the double, one obtains a conformal map of the conjugate surface by extending the Beltrami differential by zero. By considering Schwarzian derivatives of these conformal maps, we obtain a model of Teichm\"uller space by quadratic differentials on the conjugate.
   
   Given any quasiconformal map $F:\riem \rightarrow \riem_1$, let 
   \[  \mu = \frac{\overline{\partial} F}{\partial F}  \]  be its Beltrami differential. Bers' trick is to extend $\mu$ by zero to the double $\riem^d$, that is
   \[  {\mu}^d(z)  = \left\{ \begin{array}{ll} \mu & z \in \riem \\ 0 & z \in \riem^*. \end{array} \right. \]
   Solving the Beltrami equation on $\riem^d$, one obtains a quasiconformal map 
   \begin{equation} \label{eq:Bers_trick}
     {F}^d: \riem^d \rightarrow \riem^d_1 
   \end{equation}
   where $\left. {F}^d \right|_{\riem^*}$ is conformal. 

   To associate a quadratic differential, represent $\riem$ as a quotient $\disk_-/G$ for a Fuchsian group $G$, and let $F_{l}:\disk_+\rightarrow \sphere$ be a lift of $F^d$. 
   The Schwarzian derivative 
   \[\mathcal{S}(F_l) \, dz^2:=  \left( \frac{F_l'''}{F_l'} - \frac{3}{2} \left( \frac{F_l''}{F_l'} \right)^2 \right) dz^2     \]
   is in $A^\infty_2(\disk_+,G)$ and thus gives a well-defined element of $A^\infty_2(\riem^*)$.  The resulting map 
   \[  \beta: T(\riem) \rightarrow A_2^\infty(\riem^*)  \]
   is called the {\it Bers embedding}.
   The Bers embedding depends on the base surface $\riem$. When it is necessary to emphasize this, we will write $\beta_\riem$, or $\beta_G$ in the lifted picture when $\riem \cong \disk_-/G$. 
   
   The importance of the Bers embedding stems from the following fact
   \begin{theorem}
     $\beta$ is well-defined and injective.
   \end{theorem}  
 
 \end{subsection}
 \begin{subsection}{Caps fiber model}  \label{se:fiber_model_general_description}

   In this section, we consider bordered surfaces of genus $g$ with $n$ borders homeomorphic to $\mathbb{S}^1$; recall that we refer to these as type $(g,n)$ bordered surfaces.  
D. Radnell and the first author showed that the Teichm\"uller space of such surfaces fibers over the Teichm\"uller space of compact surfaces with $n$ punctures, and the fibers can be identified with a collection of $n$-tuples of conformal maps into the surface \cite{RadnellSchippers_fiber}. The idea is motivated by conformal field theory (see Section \ref{se:conformal_field_theory} ahead). This fibration features in one model of general Weil-Petersson Teichm\"uller spaces, and also relates to some of the physical motivation for the Weil-Petersson theory.\\

The idea is as follows. Given a type $(g,n)$ Riemann surface $\riem$, instead of sewing on the conjugate surface $\riem^*$ to obtain the double $\riem^d$, one sews a punctured disk to each boundary curve to obtain a compact surface $\riem^P$ with $n$ punctures.  Instead of considering conformal maps of $\riem^*$ as in the Bers embedding, we consider conformal maps of the sewn disks. We will see that variations in Teichm\"uller space $T(\riem)$ are spanned by variations in $T(\riem^P)$ together with variations obtained by changing these conformal maps (equivalently, by changing the boundary values of the quasiconformal deformation of $\riem$).\\ 

The fibres are locally modelled by powers of the following space:   
 \[ \mathcal{O}^{\mathrm{qc}} = \{ f: \disk_+ \rightarrow \mathbb{C} : f \ \mathrm{holomorphic\,\, and} \ \mathrm{quasiconformally\,\, extendible} \}.     \]
 We then consider $n$-tuples of quasiconformally extendible maps into a Riemann surface 
\begin{definition} \label{de:standard_Oqc} Let ${\riem}^P$ be a compact Riemann surface of genus $g$  with $n$ punctures $p_1,\ldots,p_n$. 
 The class $\mathcal{O}^{\mathrm{qc}}({\riem}^P)$ is the set of $n$-tuples of maps 
 $(f_1,\ldots,f_n)$ where 
 \begin{enumerate}
     \item $f_k:\disk_+ \rightarrow {\riem}^P$ are holomorphic and quasiconformally extendible to an open neighbourhood of the closure of $\disk_+$ for $k=1,\ldots,n$;
     \item $f_k(0)=p_k$ for $k=1,\ldots,n$;
     \item the closures of $f_k(\disk_+)$ are pairwise disjoint.
 \end{enumerate}
\end{definition} 
 
To obtain a fibration of $T(\riem)$ over $T(\riem^P)$, we sew on copies of the punctured disk (``caps'').  Fix a bordered surface $\riem$ of type $(g,n)$ and quasisymmetric parametrizations $\tau_k:\mathbb{S}^1 \rightarrow \partial_k \riem$ for $k=1,\ldots,n$. These can be chosen to be analytic if desired. We sew on $n$ copies of the disk $\disk_+$, by identifying $p \in \mathbb{S}^1$ with $\tau(p) \in \partial_k \riem$. The result is a compact surface $\riem^P$ whose complex structure is the unique one compatible with $\mathscr{R}$ and the copies of $\disk_+$. The parametrizations $\tau_k$ then extend to conformal maps $\disk_+ \rightarrow \riem^P$ (the inclusion maps into the sewn surface), and we specify the punctures $p_k=\tau_k(0)$ for $k=1,\ldots,n$. 

The fibre map from $T(\riem)$ to $T(\riem^P)$ is now defined as follows, assuming that $2g-2+n>0$. Given $[\riem,F,\riem_1] \in T(\riem)$, let $\hat{F}:\riem^P \rightarrow \riem^P_1$ be a quasiconformal map such that 
\[  \frac{\overline{\partial} \hat{F}(z)}{\partial \hat{F}(z)} = \left\{  \begin{array}{cc}  \overline{\partial}{F}(z)/\partial F(z)  & z \in \riem \\ 0 & z \in \tau_k(\disk_+)  \end{array} \right.  \]

Define 
\begin{align} \label{eq:fibre_C_definition}
    \mathcal{C}:T(\riem) & \rightarrow T(\riem^P)  \nonumber \\
    [\riem,F,\riem_1] & \mapsto [\riem^P,\hat{F},\riem^P_1].
\end{align}
Note that $\mathcal{C}$ depends on $\tau_1,\ldots,\tau_n$.  Later we will see that $\mathcal{C}$ is holomorphic.

\begin{remark} \label{re:Bers_trick_extended}
 One may think of the definition of $\hat{F}$ as a natural variation on ``Bers' trick''. That is, $\hat{F}$ plays the same role as   (\ref{eq:Bers_trick}). The distinction is that for $\hat{F}$ one sews on caps, while for $F^d$ one sews on the double.
\end{remark}

Fix any $q=[\riem^P,F,\riem_1^P]$. Up to the action of $\mathrm{DB}$, the fibres of $\mathcal{C}$ can be identified with $\mathcal{O}^{\mathrm{qc}}(\riem^P_1)$.
Explicitly, $f=(f_1,\ldots,f_n) \in \mathcal{O}^{\mathrm{qc}}(\riem^P_1)$, let $\riem_1$ be obtained by removing the disks $f_k(\disk_+)$. 
Let
\begin{equation*}  
 F_f':\riem^P \rightarrow \riem_1^P 
\end{equation*}
be a quasiconformal map which is homotopic to $\left. F \right|_\riem$ and such that $F \circ \tau_k = f_k$.  (It can be shown that such a map exists \cite{RadnellSchippers_fiber}.)
Set 
\begin{equation} \label{eq:F_f_definition}
 F_f = \left. F_f' \right|_{\riem} : \riem \rightarrow \riem_1
\end{equation} 
Then $[\riem,F,\riem_1]$ is in the fibre $\mathcal{C}^{-1}(q)$, and any two such elements are related by an element of the modular group $\mathrm{DB}$.
We then define 
\begin{align} \label{eq:fibre_F_definition}
  \mathbf{F}_q:\mathcal{O}^{\mathrm{qc}}(\riem_1^P) & \rightarrow \mathcal{C}^{-1}(q)/\mathrm{DB}  \nonumber \\
  f & \mapsto [\riem,F_f,\riem_1]. 
\end{align}
\begin{theorem}[\cite{RadnellSchippers_fiber}] \label{th:Fq_is_bijection} Let $\riem$ be a bordered surface of type $(g,n)$ with $2g-2+n>0$, and let $\mathcal{C}$, $\mathbf{F}_q$, etc. be as above.  
   For any point $q=[\riem^P,F,\riem_1^P] \in T(\riem^P)$, $\mathbf{F}_q$ is a bijection.
\end{theorem}
\begin{remark} \label{re:fiber_exceptional_cases}
  Since we assume that $n \geq 1$, the restriction that $2g-2+n>0$ only rules out the following two cases: $g=0, n=1$ (the disk) and $g=0,n=2$ (the annulus). These cases are exceptional because when disks are sewn on, the Teichm\"uller spaces of the resulting punctured surfaces, the once and twice punctured sphere respectively, each reduce to a point. The fiber model must therefore be adjusted slightly.
  
  The first case is just the disk, and the construction above reduces to one of the standard models of the universal Teichm\"uller space (cf Remark \ref{re:Bers_trick_extended}). Explicitly, letting $\riem^P = \sphere \backslash \{0 \}$, we have 
  $T(\disk_-) \cong \mathcal{O}^{\mathrm{qc}}(\sphere \backslash \{0 \})/G$ where $G$ are M\"obius transformations fixing $0$ (equivalently, one may apply two further normalizations to $\mathcal{O}^{\mathrm{qc}}(\sphere \backslash \{0 \})$). 

  A similar fiber model occurs in the case of an annulus $\mathbb{A}$; one can identify $T(\mathbb{A})/\mathbb{Z}$ with elements of $\mathcal{O}^{\mathrm{qc}}(\sphere \backslash \{0,\infty\})$ with one further normalization. This is called the Neretin-Segal semigroup, see \cite{RadnellSchippers_annulus}. 
\end{remark}
 
 \end{subsection} 
\begin{subsection}{Right translation and change of base point}

 A quasiconformal map $F_\mu:\riem \rightarrow \riem_\mu$ induces a bijection between Teichm\"uller spaces $T(\riem)$ and $T(\riem_\mu)$. Here we follow a common notational convention which uses the Beltrami differential of $F$ simultaneously as a label. 
 Namely, we have
 \begin{align*}
  \mathrm{R}_{F_\mu^{-1}}: T(\riem) & \rightarrow T(\riem_\mu) \\
  [\riem,F,\riem_1] & \mapsto [\riem_\mu,F \circ F_\mu^{-1},\riem_1].
 \end{align*} 
 
 In the lifted picture, The change of base-point takes the following form. Given a quasiconformal map $w_\mu:\disk_- \rightarrow \disk_-$ fixing $-1$, $i$ and $1$, we get 
 \begin{align*}
  \mathrm{R}_{w_\mu^{-1}}: T(G) & \rightarrow T(G_\mu) \\
  [w] & \mapsto [w \circ w_\mu^{-1}].
 \end{align*}
 Here, 
 \begin{equation} \label{eq:change_base_point_def_group}
  G_\mu = \{  w_\mu \circ g \circ w_\mu^{-1}  : g \in G   \}. 
 \end{equation}
 The Beltrami differential of $w_\mu \circ g \circ w_\mu^{-1}$ is zero; thus it is in $\text{M\"ob}(\mathbb{S}^1)$, and independent of the representative $\mu$. 
 Change of base-point is a homeomorphism; indeed, when complex structures are introduced, it is a biholomorphism. 
\end{subsection}
\end{section}
\begin{section}{Complex manifold and tangent space structure} 
\label{se:classical_complex_structure}
\begin{subsection}{Summary of the $L^\infty$ theory}
\begin{subsubsection}{\underline{Bers embedding}}
 \label{subsubse:Bers_embedding_complex}
 The nature of the Bers embedding was studied for many years, and eventually understood to be a homeomorphism onto an open set. In summary,
 \begin{theorem} The image of $\beta$ is an open set in $A_2^\infty(\riem^*)$. Furthermore $\beta$ is a homeomorphism from $T(\riem)$ $($with the topology induced by the Teichm\"uller distance$)$ onto its image.  
 \end{theorem}

 In other words, every element of the Teichm\"uller space can be represented by a unique quadratic differential on $\riem^*$ (equivalently, by a unique quadratic differential in $A_2^\infty(\disk_+,G)$). Since the image is open and $A_2^\infty(\riem^*)$ is a Banach space, this makes $T(\riem)$ a complex Banach manifold with a global coordinate.

 Let 
 \[ \Phi:L^\infty_{-1,1}(\riem)_1 \rightarrow T(\riem) \] denote the map obtained by solving the Beltrami equation. This map depends on the base surface $\riem$; when it is necessary to emphasize this fact we will write $\Phi_\riem$, or $\Phi_G$ in the lifted picture when $\riem=\disk_-/G$.  

 We will see that $\Phi$ is a submersion, that is, it has local holomorphic sections. 
 \begin{theorem}
     The map $\beta \circ \Phi$ is a holomorphic submersion. 
 \end{theorem}
 This is actually at the root of the deformation model. 
  
\end{subsubsection}
\begin{subsubsection}{\underline{Deformation model}}  \label{subsubse:deformation_complex_structure}
 The complex structure in the deformation model is more involved, because the Beltrami differentials contain redundant information, which the Teichm\"uller equivalence relation removes. We also see this in the description of the tangent space in Section \ref{subse:tangent_space} ahead.
 
At the root of the classical approach is the Ahlfors-Weill reflection. The abstract geometric definition of this reflection is quite simple: 
\begin{align} 
 \Lambda_\riem:A_2^\infty(\riem^*) & \rightarrow L^\infty_{-1,1}(\riem) \nonumber \\
 \alpha & \mapsto - \frac{1}{2}  \rho_\riem^{-1} i^* \alpha.
\end{align} 
In the lifted picture, using the fact that the lift of $i$ is the map $z \mapsto 1/\bar{z}$, this can be written explicitly as follows:
\begin{align} \label{eq:Ahlfors_Weill_reflection_definition}
 \Lambda_{\disk_-}:A_2^\infty(\disk_+) & \rightarrow L^\infty_{-1,1}(\disk_-) \nonumber \\
 Q(z) dz^2 & \mapsto - \frac{1}{2} \frac{1}{\bar{z}^4} \frac{Q(1/\bar{z})}{(|z|^2-1)^2} \frac{d\bar{z}}{dz}.
\end{align} 
It is easily checked that $\Lambda_{\disk_-}$ preserves $G$-invariance. Thus we can define 
\[  \Lambda_G = \left. \Lambda_{\disk_-} \right|_{A_2(\disk_+,G)}: A_2^\infty(\disk_+,G) \rightarrow L_{-1,1}^\infty(\disk_-,G)   \]
so that $\Lambda_G$ is the lift of $\Lambda_\riem$ when $\riem$ is $\disk_-/G$.

Recalling the definition of harmonic Beltrami differentials (\ref{eq:harmonic_Beltrami_definition}) we have 
\begin{align*}
 \Omega_{-1,1}(\riem) & = \Lambda(A_2^\infty(\riem^*)) \\
 \Omega_{-1,1}(\disk_-;G) & = \Lambda(A_2^\infty(\disk_+;G)).  
\end{align*} 

An alternate complex structure in terms of Beltrami differentials is provided by the following.  Let 
\[ \Phi_{\riem}:L^\infty_{-1,1}(\riem) \rightarrow T(\riem)  \]
denote the map obtained by solving the Beltrami equation (equivalently, quotienting by the Teichm\"uller equivalence relation). Denote the lifted map by $\Phi_G$.
The Ahlfors-Weill map is a local inverse of the Bers embedding near zero; namely, there is an open neighbourhood $U_\riem \subseteq A_2^\infty(\riem^*)$ of $0$ (in fact, the ball of radius $2$ centred on $0$) such that $\beta_{\riem} \Phi_\riem  \Lambda_{\riem} =\mathrm{Id}$ on $U_{\riem}$.  Similarly, $\Phi_G  \Lambda_G$ is a local inverse of $\beta_G$ on an open set $U_G$ in $A_2^\infty(\riem)$. Together with the translation maps, we obtain a complex atlas.   
\begin{theorem}  
 Let $F_\mu: \riem \rightarrow \riem_\mu$. 
 The open sets $R_{F_\mu} \Phi_{\riem_\mu} \Lambda_{\riem_\mu} (U_{\riem_\mu}) \subset T(\riem)$ and charts $\beta_{\riem_\mu}  R_{F_\mu^{-1}}$ form a complex atlas of $T(\riem)$. 

 In the lifted picture, the open sets $R_{w_\mu} \Phi_G \Lambda_{G_\mu} (U_{G_\mu})$ and charts $\beta_{G_\mu}  R_{w_\mu^{-1}}$ form a complex atlas of $T(G)$.  
\end{theorem}
In other words, we can use the fact that $\Phi_{\riem_\mu} \Lambda_{\riem_\mu}$ is a biholomorphism to obtain a chart near $[\riem,F_\mu,\riem_\mu]$. 
 
\begin{theorem} Let $\riem$ be a Riemann surface covered by the disk. 
 \begin{enumerate}[font=\upshape]
  \item The two complex structures are compatible. In particular, the Bers embedding is a biholomorphism onto its image, with respect to the complex structure on $T(\riem)$ induced by the charts.
  \item $\beta \Phi$ is a holomorphic submersion. In particular, the complex structure of $L^\infty_{-1,1}(\riem)$ and that of $T(\riem)$ are compatible. 
  \item Change of base point is a biholomorphism.
 \end{enumerate} 
\end{theorem}
\end{subsubsection}
\begin{subsubsection}{\underline{Fiber model}}  \label{se:fiber_model_local_coordinates}

We outline the basic results about the holomorphic fibration here for the case of classical $L^\infty$ Teichm\"uller theory.

The local model of the conformal maps is given by 
\[ \mathcal{O}^{\mathrm{qc}} = \{ f: \disk_+ \rightarrow \mathbb{C} : f \ \mathrm{holomorphic\,\, and\,\, quasiconformally \,\,extendible} \}.     \]
It follows from classical function theory that 
\begin{align*}
    \chi: \mathcal{O}^{\mathrm{qc}} & \mapsto A_1^\infty(\disk_+) \oplus \mathbb{C} \\
    f & \mapsto \left(\frac{f''}{f'} \, dz, f'(0) \right),
\end{align*}
is bounded with respect to the direct sum norm $|f'(0)| + \left\| f''/f' \, dz \right\|_\infty$. Furthermore, the image is open (see \cite{RSnonoverlapping}) and thus $\mathcal{O}^{\mathrm{qc}}$ is a complex Banach manifold.  

We construct a coordinate chart in $\mathcal{O}^{\mathrm{qc}}(\riem^P)$ as follows. Let $\zeta = (\zeta_1,\ldots, \zeta_n)$ be coordinates $\zeta_k:B_k \rightarrow \mathbb{C}$ of $p_k$ for $k=1,\ldots,n$. Let $K=(K_1,\ldots,K_n)$ be compact sets in $B_k$ whose interiors are open sets containing $p_k$.  Define 
\begin{equation} \label{eq:V_open_sets_definition}
  V_{\zeta,K} = \{ f=(f_1,\ldots,f_n) \in \mathcal{O}^{\mathrm{qc}}(\riem^P) : \mathrm{cl} (\, f_k(\disk_+)) \subset K_k, k=1,\ldots,n \}.    
\end{equation}

We then define coordinates
\begin{align}   \label{eq:T_coordinates_definition}
 \mathcal{E}:V_{\zeta,K} & \rightarrow (\mathcal{O}^{\mathrm{qc}})^n  \nonumber\\
 f & \mapsto \zeta \circ f. 
\end{align} 
Using the coordinates $\chi^n \circ \mathcal{E}$ we obtain the following. 
\begin{theorem}[\cite{RSnonoverlapping}]
  $\mathcal{O}^{\mathrm{qc}}({\riem}^P)$ is a Banach manifold  locally modelled on $\bigoplus^nA_1^\infty(\disk_+) \oplus \mathbb{C}$. 
\end{theorem}

The complex structure on $\mathcal{O}^{\mathrm{qc}}(\riem^P)$ is compatible with that induced by the fibration $\mathcal{C}:T(\riem) \rightarrow T(\riem^P)$ a complex structure compatible with that of $T(\riem)$, by the following theorem:  

\begin{theorem}[\cite{RadnellSchippers_fiber}]  \label{th:fiber_proj_holo_classical} Let $\riem$ be a bordered surface of type $(g,n)$ such that $2g-2+n>0$.
\begin{enumerate}[font=\upshape]
 \item The map $\mathcal{C}$ is holomorphic and has local holomorphic sections. In particular the fibre $\mathcal{C}^{-1}(q)$ is a complex submanifold of $T(\riem)$ for any $q \in T(\riem^P)$. 
 \item 
   For any point $q=[\riem^P,F,\riem_1^P] \in T(\riem^P)$, $\mathbf{F}_q$ given by \eqref{eq:fibre_F_definition}, is a biholomorphism. 
\end{enumerate}
\end{theorem}

Using this, it is possible to obtain holomorphic coordinates on $T(\riem)$. Let $T,V_{\zeta,K}$ be coordinates on $\mathcal{O}_{\mathrm{qc}}(\riem^P_1)$ as above, and set $U_{\zeta,K}= T (V_{\zeta,K})$. By an adaptation of Schiffer variation developed by F. Gardiner, there is an open set $U \subseteq \mathbb{C}^d$ containing $0$ where $d$ is the dimension of $T(\riem^P)$, and a biholomorphism 
\begin{align*}
    U & \rightarrow T(\riem^P) \\
    \epsilon & \mapsto [\riem^P,
\nu_\epsilon \circ F,\nu_\epsilon(\riem^P_1)]
\end{align*}
where $\nu_\epsilon$ is conformal on the domain of each chart $\zeta_i$. 

Now, given any point $p \in \mathcal{C}^{-1}(q)$, it can be written $[\riem,F_f,\riem_1] \in  {T}(\riem)$ for some $f \in \mathcal{O}^{\mathrm{qc}}(\riem^P)$ (since $\mathbf{F}_q$ is a bijection by Theorem \ref{th:Fq_is_bijection}, one need only adjust by an element of $\mathrm{DB}$.)   One then defines 
\begin{align} \label{eq:fiber_coordinates_definition}
    \mathcal{G}: U \times V_{\zeta,K} & \rightarrow T(\riem) \nonumber \\
    (\epsilon,f) & \mapsto   [\riem,\nu_\epsilon \circ F_f,\nu_\epsilon(\riem_1)]. 
\end{align}
This is a local biholomorphism.
\begin{theorem}[\cite{RadnellSchippers_fiber}] \label{th:fiber_coordinates_classical}
Let $\riem$ be a bordered surface of type $(g,n)$ such that $2g-2+n>0$. 
 $\mathcal{G}$ is a biholomorphism onto its image. In particular, the collection of maps of the form $\mathcal{G}^{-1}$ forms an atlas of a complex structure on $T(\riem)$ compatible with the standard complex structures. 
\end{theorem}

By composing the second component with $\chi^n \circ \mathcal{E}$, we can get charts directly to the Banach space $\mathbb{C}^d \oplus \bigoplus^n \left[ A_1^\infty(\disk_+) \oplus \mathbb{C} \right]$, if desired. 

\begin{remark} 
 Although the notation is involved, the idea is quite simple.  The Teichm\"uller space of $\riem$ decomposes locally into deformations of the capped surface $\riem^P$ and the transverse deformations of the boundary values of the quasiconformal maps of $\riem$.  
 The variations of the boundary values can be modelled by conformal maps into $\riem^P$, by a modification of Bers' trick.
\end{remark}
\begin{remark}  We saw in Remark \ref{re:fiber_exceptional_cases} that the exceptional cases are $\riem=\disk_-$ and $\riem=\mathbb{A}$ where $\mathbb{A}$ is an annulus. Tracing this through in the case of the disk just recovers the standard model in terms of the pre-Schwarzian version of the Bers embedding. Versions of Theorems \ref{th:fiber_proj_holo_classical} and \ref{th:fiber_coordinates_classical} was given for $T(\mathbb{A})$ in \cite{RadnellSchippers_annulus}. 
\end{remark}
\begin{remark}
 The holomorphicity of $\mathbf{F}_q$ is \cite[Corollary 3.3]{RadnellSchippers_fiber} while the biholomorphicity of $\mathcal{G}$ is \cite[Theorem 4.1]{RadnellSchippers_fiber}, with rather different notation. 
\end{remark}
\begin{remark}
 The idea for these coordinates is due to D. Radnell \cite{Radnell_thesis}. 
\end{remark}
 
 \end{subsubsection} 
\end{subsection}

\begin{subsection}{The tangent space in the three models}\label{subse:tangent_space}

{\underline{Deformations model}}: There are many quasiconformal deformations in a given equivalence class. Thus if we are to use  $L^\infty_{-1,1}(\riem)$ in modelling the tangent space to the Teichm\"uller space, we must quotient out directions tangent to  $L^\infty_{-1,1}(\riem)$.  These redundant directions form what is called the space of ``infinitesimally trivial'' differentials. Namely, define
 \[  \mathcal{N}(\riem) =  \mathrm{Ker} D_0 (\beta_\riem \circ \Phi_\riem) \]
  \begin{equation}\label{derivative} D_0 (\beta_\riem \circ \Phi) : L^\infty_{-1,1} \rightarrow A_2^\infty(\riem^*)   \end{equation}
   denotes the derivative of  $\beta_\riem \circ \Phi_\riem$ at the base point $[0]=[\riem,\mathrm{Id},\riem]$. 
   Thus $\mathcal{N}(\riem)$ is the set of directions tangent to the equivalence relation, by the fact that the Bers embedding is a homeomorphism onto its image.  
   It is a fundamental result that
   \begin{equation} \label{eq:harmonic_Bel_decomposition}
     L^\infty_{-1,1}(\riem) = \mathcal{N}(\riem) \oplus \Omega_{-1,1}(\riem).   
   \end{equation}
   Thus if we let $[0]=[\riem, \mathrm{Id},\riem]$ denote the ``origin'' of Teichm\"uller space, the tangent space at $[0]$ is 
   \[ T_{[0]} T(\riem) \cong \Omega_{-1,1}(\riem).  \]

   The tangent space at an arbitrary point $q=[\riem,F,\riem']$ can be obtained using right translation. Since 
   \[ R_{[F]}[\riem',\mathrm{Id},\riem']=[\riem,F,\riem'] \]
   the derivative satisfies
   \[ D_0 R_{[F]}: T_{[0]} T(\riem') \rightarrow T_q T(\riem).  \]
   Since $\Omega_{-1,1}(\riem')$ models $T_{[0]} T(\riem')$ we obtain
   that 
   \begin{equation}  \label{eq:arbitrary_tangent_space_model}
     T_q T(\riem) \cong R_{[F]} \Omega_{-1,1}(\riem').   
   \end{equation}
   One can view the coordinates established in Section \ref{subsubse:deformation_complex_structure} as local exponentiations of this model of the tangent space. 

   {\underline{Bers embedding model}}:  Since the Bers embedding is a bijection from $T(\riem)$ to an open subset of $A_2^\infty(\riem^*)$, there are no redundant directions in $A_2^\infty(\riem^*)$. Thus the tangent space at $[0]$ can identified with $A_2^\infty(\riem^*)$, that is
   \[ T_{[0]} T(\riem) \cong A_2^\infty(\riem^*).   \]
   Since $A_2^\infty(\riem^*)$ is a linear space, we also have 
   \[   T_q T(\riem) \cong T_{\beta(q)} A_2^\infty(\riem^*) \cong A_2^\infty(\riem^*). \]

 {\underline{Fiber model}}: We remark here that it is possible to get a model of the tangent space using the fibration. This will not be used elsewhere in the paper.
  Let $p = [\riem,F,\riem_1]$ and $q = \mathcal{C}(p) = [\riem^P,\hat{F},\riem_1^P]$. By Theorem \ref{th:fiber_coordinates_classical}, using the linearity of $A_1^\infty(\disk_+)$, 
  \[ T_{[\riem,F,\riem_1]} T(\riem) \cong T_q T(\riem^P) \times \bigoplus_{k=1}^n \left(A_1^\infty(\disk_+) \oplus \mathbb{C} \right)   \]
  via the derivative map of $\mathcal{G}$.  
  Since $T_q T(\riem^P)$ is finite dimensional, it poses no analytic difficulties. Under $\mathcal{G}$ it is identified with $\mathbb{C}^d$ using Schiffer variation, but any model of the tangent space could be used.   
\end{subsection}
\begin{subsection}{Weil-Petersson Teichm\"uller space in a nutshell}

In each of the models given in Section \ref{se:three_models}, an analytic choice leads to the complex manifold structure on Teichm\"uller space.  The Weil-Petersson Teichm\"uller space is, roughly, obtained by replacing the $L^\infty$ objects by $L^2$ objects. For example, in the Weil-Petersson theory the space of $L^2$ harmonic Beltrami differentials 
\[ H_{-1,1}(\riem) = \left\{ \mu \in \mathscr{D}_{-1,1}(\riem) \,: \, 
 \rho_\riem \mu \in \overline{A_2^2(\riem)}  \right\} \]
plays the role of $\Omega_{-1,1}(\riem)$.

 \begin{center}
\begin{tabular}{||c || c | c ||} 
 \hline 
 & {\bf Classical} & {\bf Weil-Petersson}  \\ [0.5ex] 
 \hline\hline
  {\bf Beltrami differentials} & $L^\infty_{-1,1}(\riem)_1 $ & $L^\infty_{-1,1}(\riem)_1 \cap L^2_{-1,1}(\riem)$    \\ [1ex] 
 \hline
   {\bf Deformation model coord.}  & $\Omega_{-1,1}(\riem)$ & $H_{-1,1}(\riem)$    \\  [1ex] 
 \hline 
 {\bf Bers embedding model coord.}  & $A_2^\infty(\riem^*)$ & $A_2^2(\riem^*)$  \\  [1ex] 
 \hline  
 {\bf Fiber model coord.} & $(A_1^\infty(\disk_+) \oplus \mathbb{C})^n \times T(\riem^P)$ & $(A_1^2(\disk_+) \oplus \mathbb{C})^n \times T(\riem^P)$   \\ [1ex] 
 \hline
\end{tabular}
\end{center}
 The Weil-Petersson inner product of the tangent space at $[0]=[\riem,\mathrm{Id},\riem]$ is as follows. In the Bers embedding model, we will see ahead that the tangent space can be identified with $A_2^2(\riem^*)$.  In this model, the pairing is just the inner product induced by that in $A_2^2(\riem^*)$: 
 \[   \left< \phi_1,\phi_2 \right>_{\mathrm{WP},0} : = \left< \phi_1,\phi_2 \right>_{A^2_2(\riem^*)} = \iint_{\riem^*} \rho_{\riem^*}^{-2} \phi_1 \overline{\phi_2} \, dA_{\mathrm{hyp}}. \] 
 In the lifted picture we have, for $h_1(z)dz^2,h_2(z)dz^2 \in A_2^2(\disk_+,G)$ 
 \[   \left< h_1(z)dz^2,h_2(z)dz^2 \right>_{\mathrm{WP},0} : =  \left< h_1(z)dz^2,h_2(z)dz^2 \right>_{A^2_2(\disk_+,G)} =  \iint_{N} (1-|z|^2)^2 h_1(z) \overline{h_2(z)} \frac{d\bar{z} \wedge dz}{2i} \] 
 where $N$ is a fundamental domain of $G$.  

 In the deformation model, the tangent space is modelled by harmonic Beltrami differentials. In that case, for $\mu_1,\mu_2 \in H_{-1,1}(\riem)$ the Weil-Petersson pairing is 
 \[   \left< \mu_1,\mu_2 \right>_{\mathrm{WP},0} = 
        \iint_{\riem} \mu_1 \overline{\mu_2} \, dA_{\mathrm{hyp}}.    \]
 It is easily checked that up to a constant, under the Ahlfors-Weill reflection $\Lambda$, this agrees with the pairing on quadratic differentials.  
 
 Before this can be made sense of, we first need to define the Weil-Petersson Teichm\"uller space for some class of surfaces $\riem$, and show that it has a Hilbert manifold complex structure. 
 The existence and equivalence of the various complex structures on Teichm\"uller space was an important research question, with origins in Riemann's work.  In the finite-dimensional and classical Ahlfors-Bers theory this encompasses many important results. The difficulties are no less in the Weil-Petersson theory.

 \begin{remark}
 The construction of the Weil-Petersson Teichm\"uller space and its equivalent complex structures requires repeating all the analysis in the classical construction, but in the $L^2$ case. Although there are quite a few surprises, there is also a long list of ``expected'' theorems, some of which have technical proofs. In the initial stages, the value of the endeavour was perhaps not as obvious as it is in hindsight.
 \end{remark}
\end{subsection}
\begin{subsection}{The tangent space to Teichm\"uller space: classical and Weil-Petersson theory}

\end{subsection}
\end{section} 
\begin{section}{The Weil-Petersson Universal Teichm\"uller space}
\label{se:WP_Universal_Teich}
\begin{subsection}{History and overview} \label{se:WPUniversal_history}
  The early motivations for the Weil-Petersson Teichm\"uller space came from investigations of the group $\mathrm{Diff}(\mathbb{S}^1)$ of diffeomorphisms of the circle. This group plays and important role in loop groups, representation theory, string theory, and conformal field theory \cite{Khesin_Wendt, Pressley_Segal,BowickRajeev, BowickRajeev_short, Pekonen, Schottenloher, Witten}. Its coadjoint orbits have been investigated in connection with representation theory by A. Kirillov and D. Yuri'ev \cite{KY2,Kirillov_short}. There are two coadjoint orbits, $\mathrm{Diff}(\mathbb{S}^1)/\text{M\"ob}(\mathbb{S}^1)$ and $\mathrm{Diff}(\mathbb{S}^1)/\mathbb{S}^1$ (where $\mathbb{S}^1$ acts by the subgroup of rotations in $\text{M\"ob}(\mathbb{S}^1)$). Kirillov-Yuri'ev identified the symplectic forms on these orbits arising from the Kirillov-Kostant-Souriau theory \cite{Kirillov_short,KY2}.  The tangent space at the identity of $\mathrm{Diff}(\mathbb{S}^1)$ consists of the smooth vector fields on the circle $\mathrm{Vect}(\mathbb{S}^1)$, generated by $L_n =-i e^{in \theta} \partial/\partial \theta$, $n \in \mathbb{Z} \backslash \{-1,0,1\}$. There is a natural complex structure $J$ given by 
  \[  J L_m =  -i \, \mathrm{sgn}(m) L_m.    \]
  
  The unique homogeneous symplectic form on $\mathrm{Diff}(\mathbb{S}^1)/\text{M\"ob}(\mathbb{S}^1)$ is given at the identity by \cite{KY2,Kirillov_short,Bowick_Lahiri} 
  \[  \omega(L_m,L_n) = a (m^3-m) \delta_{m+n,0}   \]
  for some real constant $a$. 
  There is thus a unique K\"ahler metric compatible with the complex structure $J$.\\ 

  Motivated by the work of M. Bowick \cite{Bowick} and M. Bowick and S. Rajeev \cite{BowickRajeev},
  S. Nag and A. Verjovsky \cite{NagVerjovsky} observed the connection of $\mathrm{Diff}(\mathbb{S}^1)/\text{M\"ob}(\mathbb{S}^1)$ to the universal Teichm\"uller space. Using the deformation model of the tangent space, they showed that the K\"ahler metric agrees with the formal expression of the Weil-Petersson pairing of Beltrami differentials in $T(\disk_-)$, and that this pairing therefore converges at points in $\mathrm{Diff}(\mathbb{S}^1)/\text{M\"ob}(\mathbb{S}^1) \subset T(\disk_-)$ on directions tangent to $\mathrm{Diff}(\mathbb{S}^1)/\text{M\"ob}(\mathbb{S}^1)$.  Furthermore, they showed that the inclusion of $\mathrm{Diff}(\mathbb{S}^1)/\text{M\"ob}(\mathbb{S}^1)$ in $T(\disk_-)$ is G\^{a}teaux holomorphic. 

  Setting geometry and algebra aside for now, the question which immediately leaps to mind is: what is the completion of $\mathrm{Diff}(\mathbb{S}^1)/\text{M\"ob}(\mathbb{S}^1)$ with respect to the Weil-Petersson pairing of Nag-Verjovsky?\\ 
  
   Since $T(\disk_-)$ is a Banach manifold, we can certainly not expect it to have a K\"ahler metric, and indeed the Weil-Petersson pairing does not converge in arbitrary directions tangent to $T(\disk_-)$ even at the identity. Thus one is led to find a subset or refinement of the universal Teichm\"uller space, such that the Weil-Petersson pairing is convergent on every tangent space. 
    A number of researchers have developed this ``Weil-Petersson universal Teichm\"uller space''.  In this section, we describe the manifold and tangent space properties, which were developed by Cui \cite{Cui} and Takhtajan-Teo \cite{Takhtajan_Teo_Memoirs}.\\

  The seminal paper on the Weil-Petersson Teichm\"uller space was that of Cui \cite{Cui}.  Cui defined Weil-Petersson class quasisymmetries (calling them ``integrably asymptotically affine homeomorphisms''), constructed a Weil-Petersson universal Teichm\"uller space with a complex structure obtained from a holomorphic Bers embedding, and showed that it is the completion of $\text{Diff}(\mathbb{S}^1)/\text{M\"ob}(\mathbb{S}^1)$. This completely established the Bers embedding model of the WP universal Teichm\"uller space.  It was also shown that the WP universal Teichm\"uller space is complete with respect to the induced geodesic distance.\\ 
  
  Some of these results were later proven independently by Takhtajan-Teo \cite{TT_arxiv_1}, later published as the first of two chapters of  \cite{Takhtajan_Teo_Memoirs}. They also constructed the complex structure from the deformation model, and thoroughly developed the tangent space structure. Moreover they showed that a natural infinite-dimensional analogue of the period mapping into the Siegel disk is holomorphic both in the classical and Weil-Petersson settings. This infinite-dimensional analogue was given by G. Segal \cite{SegalUnitary}, and taken up by many authors including Kirillov-Yuri'ev \cite{KY2}, Nag \cite{Nag_bulletin}, Nag and Sullivan \cite{NagSullivan} and D. Hong and S. Rajeev \cite{HongRajeev}. This is described in Section \ref{se:polarizations_Siegel} ahead.
  
Takhtajan-Teo also established the first results on the K\"ahler geometry and global analysis of the WP universal Teichm\"uller space. This topic is deferred to Section \ref{se:Kahler_and_global}. 
\end{subsection}
\begin{subsection}{Weil-Petersson universal Teichm\"uller space}
\label{se:WP_universal_Cui_Takhtajan} 
  
\begin{subsubsection}{\underline{Bers embedding model}}
   
 A quasisymmetry $\phi \in \qs(\mathbb{S}^1)$ is said to be a Weil-Petersson quasisymmetry if it has a quasiconformal extension to $\disk_-$ with Beltrami differential in $L^2_{-1,1}(\disk_-)$.  Denote the class of such quasisymmetries by $\qs_{\mathrm{WP}}(\mathbb{S}^1)$.
 We then consider the following subset of Teichm\"uller space.
\begin{definition}  \label{de:WP_univ}
   The{\it{ $\mathrm{WP}$-class universal Teichm\"uller space}} is defined to be
    \[  T_{\mathrm{WP}}(\disk_-) := \qs_{\mathrm{WP}}(\mathbb{S}^1)/\text{M\"ob}(\mathbb{S}^1).  \] 
\end{definition}
 Cui obtained a characterization of Weil-Petersson quasisymmetries in terms of the Schwarzian and the pre-Schwarzian derivative. 
 Recall that given a $\phi$ with quasiconformal extension $F$, there is a corresponding map $\hat{F}$ obtained from Bers' trick, such that 
 \begin{equation} \label{eq:Bers_trick_temp}
   f := \left. \hat{F} \right|_{\disk_+} 
 \end{equation} 
 is conformal and independent of the choice of extension (up to a  normalization).  

 Define the pre-Schwarzian and Schwarzian derivatives 
 \[  \mathcal{A}(f)=\frac{f''}{f'}, \ \ \ \ \ \  \ \ \ \mathcal{S}(f) =\frac{f'''}{f'} - \frac{3}{2} \left( \frac{f''}{f'} \right)^2.  \]
 We have
 \begin{theorem}[\cite{Cui,Takhtajan_Teo_Memoirs}] \label{th:WP_quasisymmetry_characterize}
  The following are equivalent:
  \begin{enumerate}[font=\upshape]
  \item  $\phi:\mathbb{S}^1 \rightarrow \mathbb{S}^1$ is a $\mathrm{WP}$-class quasisymmetry;
  \item $f$ satisfies $\mathcal{S}(f)(z)\,dz^2 \in A_2^2(\disk_+)$;
  \item if $f$ is normalized so that $\infty$ is in the interior of the complement of $f(\disk_+)$ then $(f''(z)/f'(z))\, dz \in A_1^2(\disk_+)$. 
  \end{enumerate}
 \end{theorem}

  {It was shown by Cui \cite{Cui} and later by Takhtajan-Teo  \cite{Takhtajan_Teo_Memoirs} that}
  \begin{theorem}[\cite{Cui,Takhtajan_Teo_Memoirs}]   \label{th:disk_inclusion_bounded}
$A_2^2(\disk_+) \subset A_2^\infty(\disk_+)$
and the inclusion is bounded.
\end{theorem} 
 
Furthermore, the following holds. 
\begin{theorem}[\cite{Cui,Takhtajan_Teo_Memoirs}]  \label{th:Bers_univ_WP}
 $\beta(T_{\mathrm{WP}}(\disk_-))= \beta(T(\disk_-)) \cap A_2^2(\disk_+)$.  In particular,
 $\beta(T_{\mathrm{WP}}(\disk_-))$ is open.
\end{theorem}
By the fact that the image of the Bers embedding is open, we have that $T_{\mathrm{WP}}(\disk_-)$ has a complex Hilbert manifold structure inherited from $A_2^2(\disk)$.  

 Cui showed that the Bers embedding is holomorphic in the following sense. Modifying Cui's notation somewhat, let $M_* = L_{-1,1}^2(\disk_-)\cap L_{-1,1}^\infty(\disk_-)$ with the direct sum norm 
 \[   \| \mu \|_*= \| \mu \|_\infty + \| \mu \|_2. \]
 This is a Banach space.  Furthermore, if 
 \[ {M_*}_1 := \{ \mu \in M_* : \| \mu \|_\infty <1 \},  \]
then ${M_*}_1$ is obviously open, by the continuity of $\mu \mapsto \| \mu \|_{\infty}$ with respect to the direct sum norm. 
This yields  
 \begin{theorem}[\cite{Cui}] \label{th:Cui_Bers_starting_in_M} The restriction of the Bers isomorphism $\beta \circ \Phi$ to ${M_*}_1$ is a holomorphic map into $A_2^2(\disk_+)$.
 \end{theorem}

 The fact that the Weil-Petersson metric on $T_{\mathrm{WP}}(\disk_-)$ is convergent at every point, and the induced distance is complete, was first shown by Cui. Using the complex structure obtained from $A_2(\disk_+)$ via Theorem \ref{th:Bers_univ_WP}, we obtain that the tangent space at any point can be modelled by $A_2(\disk_+)$.  Given two elements $\phi,\psi \in A_2(\disk_+) \cong T_0 A_2(\disk_+)$, the Weil-Petersson pairing at the identity is naturally defined by the $L^2$ pairing in $A_2(\disk_+)$
 \[   \left< \phi,\psi \right>_{\mathrm{WP},0} = \iint_{\disk_+} \frac{1}{(1-|z|^2)^2} \overline{\psi(z)} \phi(z) dA_z.   \]
 At an arbitrary point the pairing is defined by right translation/change of base point. For $h \in T_{\mathrm{WP}}(\disk_-)$ let $\check{R}_h$ denote the right translation induced on $\beta(T_{\mathrm{WP}}(\disk_-))$ induced by right translation $R_{h}$ on $T_{\mathrm{WP}}(\disk_-)$. 
 We have 
 \begin{theorem}[\cite{Cui}] \label{th:Cui_right_translation}
   $\check{R}_h$ is a biholomorphism. Furthermore, there is a constant $C$ $($independent of $h$$)$ so that the derivative $D_{0}  \check{R}_{h}$ at the identity satisfies
   \[   \frac{1}{C} \leq \| D_{0} \check{R}_{h} \|_2 \leq C.  \] 
 \end{theorem} 
 
 Given any pair of tangent vectors at $h \in T_{\mathrm{WP}}(\disk_-)$ represented by $\psi, \phi \in T_{\beta(h)} A_2(\disk_+) \cong A_2(\disk_+)$, we define the Weil-Petersson pairing at $h$ by 
 \begin{equation} \label{eq:disk_WP_pairing_qd}
   \left< \phi,\psi \right>_{\mathrm{WP},h} = \left< (D_{0}  \check{R}_{h})^{-1} \phi, (D_{0} \check{R}_{h})^{-1} \psi \right>_{\mathrm{WP},0}.   
 \end{equation} 
 It is immediately seen that this is a convergent Riemannian metric on every tangent space, which is smoothly varying. In other words, we have found the largest subset of the universal Teichm\"uller space, which is a Hilbert manifold and such that the Weil-Petersson pairing converges on each tangent space. This answers the question posed in Section \ref{se:WPUniversal_history}.  
 
 Not only are the tangent spaces complete, the manifold itself is complete with respect to the distance induced by the Weil-Petersson metric.
 \begin{theorem}[\cite{Cui}]  \label{th:Cui_completeness} \label{th:WP_universal_Teich_complete}
   $T_{\mathrm{WP}}(\disk_-)$ is complete with respect to the distance induced by the Weil-Petersson pairing.  
 \end{theorem}

 It was shown by Cui \cite{Cui} that $\qs_{\mathrm{WP}}(\mathbb{S}^1)$ is, like $\qs(\mathbb{S}^1)$, a group. Takhtajan-Teo showed that it is a topological group. 
\begin{theorem}[\cite{Takhtajan_Teo_Memoirs}]   \label{th:WP_is_top_group} $\qs_{\mathrm{WP}}(\mathbb{S}^1)/\emph{M\"ob}(\mathbb{S}^1)$
 is a topological group. 
\end{theorem}
As is well-known, $\qs(\mathbb{S}^1)/\text{M\"ob}(\mathbb{S}^1)$ is not a topological group. 

A key aspect of the approach of Cui is that the Douady-Earle extension of an element of $\qs_{\mathrm{WP}}(\mathbb{S}^1)$ has Beltrami differential in $L^2_{-1,1}(\disk_-)$. The Douady-Earle extension is an explicit quasiconformal extension $E(\phi)$ to the $\disk_-$ of a quasisymmetry $\phi \in \mathrm{QS}(\mathbb{S}^1)$. It is conformally natural, that is, for any $M_1,M_2 \in \text{M\"ob}(\mathbb{S}^1)$, 
\[   E(M_1 \circ \phi \circ M_2) = M_1 \circ E(\phi) \circ M_2. \]
\begin{theorem}[\cite{Cui}] \label{th:Cui_Douady_Earle} 
   If $\phi \in \qs_{\mathrm{WP}}(\mathbb{S}^1)$ then $E(\phi)$ has Beltrami differential $\mu(E(\phi)) \in L^2_{-1,1}(\disk_-)$.  Furthermore, given $\mu \in M_{\ast 1}$ let $\sigma(\mu) \in M_{* 1}$ be the Beltrami differential of $E(\left. w_\mu \right|_{\mathbb{S}^1})$. Then $\sigma:M_{* 1} \rightarrow M_{* 1}$ is continuous. 
\end{theorem}

Finally, we observe that it follows directly from Theorem \ref{th:disk_inclusion_bounded} that 
\begin{theorem}[\cite{Cui,Takhtajan_Teo_Memoirs}] 
 \label{th:inclusion_universal_holomorphic} The inclusion from $T_{\mathrm{WP}}(\disk_-)$ into $T(\disk_-)$ is holomorphic. 
\end{theorem}
This theorem refers to the inclusion into the Banach manifold $T(\disk_-)$. As was mentioned earlier, Nag-Verjovsky \cite{NagVerjovsky} showed that the inclusion from $\text{Diff}(\mathbb{S}^1)/\text{M\"ob}(\mathbb{S}^1)$ into $T(\disk_-)$ is G\^ateaux holomorphic.   The inclusion from $T_{\mathrm{WP}}(\disk_-)$ into $T(\disk_-)$ endowed with Takhtajan-Teo's Hilbert manifold structure described in the next section, is also holomorphic. This follows more or less automatically from their construction.\\ 

Some of the results in this section were generalized to the $L^p$ case for $p \geq 1$ by H. Guo \cite{GuoHui}. S. Tang showed that the Douady-Earle
extension is in $L^p$ and that the Bers embedding is holomorphic with respect to the 
intersection norm on $L^p \cap L^\infty$. {For an overview of the $L^p$ theory as well as other refinements of the Teichm\"uller space see Section \ref{birds eye view}.} 
\begin{remark} In Guo \cite{GuoHui}, a paper of Cui is cited as ``Teichm\"uller spaces and $\text{Diff}(\mathbb{S}^1)/\text{M\"ob}(\mathbb{S}^1)$, to appear.'' This appears to refer to the paper Cui \cite{Cui}, and indeed the authors were once informed of this by an anonymous referee.  
\end{remark}
\begin{remark}
Although neither of \cite{Cui,GuoHui} uses the term ``Weil-Petersson'' in the text or title, these are as far as we know the first papers rigorously defining the Weil-Petersson Teichm\"uller space and developing its complex structure.
At the time of the publication of \cite{Cui} and \cite{GuoHui}, the journal Science in China was unfortunately still not widely available in the west.\\ 
\end{remark}
\end{subsubsection}
\begin{subsubsection}{{\underline{Deformation model}}} \label{se:universal_WP_deformation}

 The deformation model, including the model of the tangent space by harmonic Beltrami differentials, was carried out by Takhtajan-Teo \cite{TT_arxiv_1}, later published as the first of two chapters of  \cite{Takhtajan_Teo_Memoirs}. Their approach has a major difference: a Hilbert manifold structure is given to the entire Teichm\"uller space $T(\disk_-)$, at the cost of disconnecting the space. The connected component of $[0] \in T(\disk_-)$ is $T_{\mathrm WP}(\disk_-)$. 

 Recall that in the classical case the tangent space can be modelled by the harmonic Beltrami differentials $\Omega_{-1,1}(\disk_-)$, defined in terms of the Ahlfors-Weill reflection $\Lambda_{\disk_-}$ (see (\ref{eq:Ahlfors_Weill_reflection_definition})).   Takhajan and Teo define 
 \begin{equation}\label{TTs H-class}
    H_{-1,1}(\disk_-) =  \Lambda_{\disk_-} A_2^2(\disk_+), 
 \end{equation} 
 which makes sense thanks to Theorem \ref{th:disk_inclusion_bounded}.
 
  As in the $L^\infty$ case one can construct an atlas using the Ahlfors-Weill reflection. By Theorem \ref{th:disk_inclusion_bounded} there is a sufficiently small ball $B \subset A_2^2(\disk_+)$ centred at $0$ which is contained in the ball of radius $2$ in $A_2^\infty(\disk_+)$ so that all elements of $\Lambda_{\disk_-}(B)$ are Beltrami differentials of elements of $T(\disk_-)$. For $\mu \in L_{-1,1}^\infty(\disk_-)$ setting $V_\mu = R_{w_\mu} \Phi_{\disk_-} \Lambda_{\disk_-}(B)$, we obtain the following.
 \begin{theorem}[\cite{Takhtajan_Teo_Memoirs}]  Let $w_\mu:\disk_- \rightarrow \disk-$ be a normalized solution to the Beltrami equation. The sets $R_{w_\mu} \Phi_{\disk_-} \Lambda_{\disk_-}(B)$ and charts $\beta_{\disk_-} R_{w_\mu^{-1}}$ form an atlas for a complex Hilbert manifold structure on $T(\disk_-)$.   
 \end{theorem} 
 
 In the topology associated to this complex structure, $T(\disk_-)$ contains infinitely many connected components. These are integral manifolds of a certain distribution, which thus coincide with the tangent spaces. These tangent spaces are precisely the directions in $T(\disk_-)$ in which the Weil-Petersson pairing converges. The connected component of the identity is $T_{WP}(\disk_-)$.

 \begin{theorem}[\cite{Takhtajan_Teo_Memoirs}]  Right translation $($change of base point$)$ is a biholomorphism with respect to this complex structure. 
 \end{theorem}

   \begin{remark} Note that this statement is stronger than the statement given in Cui's Theorem \ref{th:Cui_right_translation} that right translation/change of base point is a biholomorphism, since it applies to all quasiconformal translations, not just Weil-Petersson class ones. If one restricts to Weil-Petersson class deformations, one obtains this part of Cui's Theorem as a special case.  
    
     On the other hand, the lower bound on the differential given in Cui's theorem, which leads to his short proof of completeness, is not contained in Theorem \ref{th:right_translation_TT}.  
   \end{remark}
 
  Recalling the decomposition of equation (\ref{eq:harmonic_Bel_decomposition}), we have
   \begin{equation}
     L^\infty_{-1,1}(\disk_-) \cap L^2_{-1,1}(\disk_-) = \left( \mathcal{N}(\disk_-) \cap L^2_{-1,1}(\disk_-) \right) \oplus H_{-1,1}(\disk_-).   
   \end{equation}
   Here we have used the fact that $H_{-1,1}(\disk_-) \subset \Omega_{-1,1}(\disk_-)$, which follows directly from Theorem \ref{th:disk_inclusion_bounded}.
   This shows that we can identify the tangent space at $[0]$ with the $L^2$ harmonic Beltrami differentials:
   \[  T_{[0]} T_{\mathrm{WP}}(\disk_-) \cong H_{-1,1}(\disk_-) \] 
   Takhtajan-Teo showed
   \begin{theorem}[\cite{Takhtajan_Teo_Memoirs}] \label{th:right_translation_TT} For any normalized quasiconformal map $w_\mu:\disk_- \rightarrow \disk_-$,
     $D_0 R_{w_\mu}$ is a bounded isomorphism. 
   \end{theorem}
   
   We then define a distribution $\mathfrak{D}_T$ over $T(\disk_-)$ as follows: for each $[w_\mu] \in T(\disk_-)$, we assign the Hilbert space $D_0 R_{w_\mu}  H_{-1,1}(\disk_-)$.  

   Takhtajan-Teo's surprising idea was to disconnect the space $A_2^\infty(\disk_+)$ and view it as a union of translates of $A_2^2(\disk_+)$. That is, consider the foliation of $A_2^\infty(\disk_+)$ 
   \[      A_2^\infty(\disk_+) = \cup_{\phi \in A^\infty_2(\disk_+)} \, \left[ \phi + A_2^2(\disk_+) \right].         \]
   Each leaf is Hilbert manifold in the obvious way, and one can view the vector space $A^\infty(\disk_+)$ as an uncountable collection of disconnected Hilbert manifolds (ignoring the original norm). The vector spaces at $\phi \in A_2^\infty(\disk_+)$  tangent to $A_2^2(\disk_+)$ form a distribution $\mathfrak{D}_A$, and the leafs of the foliation are integral manifolds of this distribution. 

   Summarizing some of their results:
   \begin{theorem}[\cite{Takhtajan_Teo_Memoirs}] \label{th:TT_Bers_embedding} 
    Consider $T(\disk_-)$ and $A_2^\infty(\disk_+)$ as Hilbert manifolds as above. Then 
     \begin{enumerate}[font=\upshape]
      \item The Bers embedding $\beta:T(\disk_+) \rightarrow A_2^2(\disk_-)$ is a biholomorphism onto its image. 
      \item The integral manifolds of the distribution $\mathcal{D}_T$ are the  inverse images of the leaves $\phi + A_2^2(\disk_+)$ of the foliation of $A_2^\infty(\disk_+)$ under the Bers embedding. 
      \item The connected component of $[0] \in T(\disk_-)$ is $\beta^{-1}(A_2^2(\disk_+))$. 
     \end{enumerate}
   \end{theorem}
   In particular, we see that $T_{\mathrm WP}(\disk_-)$ is the connected component of the identity of $T(\disk_-)$ with respect to Takhtajan-Teo's complex Hilbert manifold structure.  
   
   We also see that the tangent space can be modelled by 
   \[   T_{[w_\mu]} T(\disk_-) \cong  D_0 R_{\mu} H_{-1,1}(\disk_-)   \]
   and as we saw before, by
   \[  T_{[w_\mu]} T(\disk_-) \cong  D_0 \check{R}_{\mu} A_2^2(\disk_+).  \]
   Thus we have two different expressions for the Weil-Petersson pairing; 
   equation (\ref{eq:disk_WP_pairing_qd}) and 
   \begin{equation}  \label{eq:disk_WP_pairing_bd}
     \left< \phi,\psi \right>_{\mathrm{WP},[w_\mu]} = \left< (D_{0} R_{w_\mu})^{-1} \phi, (D_{0} R_{w_\mu})^{-1} \psi \right>_{\mathrm{WP},[0]}.
   \end{equation} 
   Tracing through the definitions shows that these are obviously equal up to a constant.

   \begin{remark}
    It should be observed that the contents of Theorem \ref{th:TT_Bers_embedding} part (1) and Theorem \ref{th:Cui_Bers_starting_in_M} are quite different, even when the Bers embedding is restricted to $T_{\mathrm{WP}}(\disk_-)$. 
   \end{remark}
\end{subsubsection}  
\end{subsection} 

\begin{subsection}{Some further characterizations of Weil-Petersson quasisymmetries}\label{Further char} 

  In this section we mention some other recent characterizations of $\mathrm{QS}_{\mathrm{WP}}(\mathbb{S}^1)$. Here we are brief, because there is already an excellent survey on the topic in C. Bishop's  paper \cite{Bishop2}. In the same paper, Bishop gave many new and striking characterizations of the WP-class, as well as generalizations to higher dimensions. The reader is encouraged to read the original papers by Bishop \cite{Bishop1, Bishop2} in order to fully appreciate the scope and depth of the results that were obtained therein, and to also see the connection of those results to many branches of mathematics as well as physics.\\
  
Takhtajan-Teo \cite{Takhtajan_Teo_Memoirs} asked a question regarding the characterization of WP-class quasisymmetries.  In this connection, Y. Shen \cite{Shen_characterization} proved the following characterization of the WP-class quasisymmetries.
  \begin{theorem}\label{shen charac}  Let \(f\) be a conformal mapping from \(\mathbb{D}_+\) onto
a quasidisk and \(h=g^{-1} \circ f\) be a corresponding quasisymmetric
conformal welding on $\mathbb{S}^{1}$. Then the following statements are equivalent
\begin{enumerate}[font=\upshape]
    \item  \(h\) is \emph{WP}-class
    \item \(h\) is absolutely continuous $($with respect to the arc-length measure$)$ and \(\log h^{\prime}\) belongs to the homogeneous Sobolev space
 \(\dot{H}^{\frac{1}{2}}\left(\mathbb{S}^{1}\right)\)
\end{enumerate}
\end{theorem} 
Later, Shen in collaboration with Y. Hu and L. Wu \cite{Shenwuhu} gave an intrinsic characterization of the elements of the Weil-Petersson class, without using quasiconformal extensions. These results answered the question of Takhtajan and Teo mentioned above.

Recently, Shen and X. Liu \cite{shenliu} proved the following characterization of the so-called $p$--integrable quasicircles for $p>1$ which in the case of WP-class quasicircles reads

\begin{theorem}\label{shen-liu}
    Let \(f\) be a conformal mapping from \(\mathbb{D}_+\) onto
a quasidisk and \(h=g^{-1} \circ f\) be a corresponding quasisymmetric
conformal welding 
 for \(\Gamma := f\left(\mathbb{S}^{1}\right)\). Then the following statements are equivalent

\begin{enumerate}[font=\upshape]
    \item \(\Gamma\) is a \emph{WP}-class quasicircle
\item  \(\log g^{\prime} \in \mathcal{D}\left(\mathbb{D}_{-}\right)\);
\item  \(\Gamma\) is rectifiable with length \(l\) and there exists some real-valued function \(b \in \dot{H}^{\frac{1}{2}}\left(\mathbb{S}^{1}\right)\) such
that an arclength parameterization \(z: \mathbb{S}^{1} \rightarrow \Gamma\) satisfies \(z^{\prime}(\zeta)=\frac{l}{2 \pi} e^{i b(\zeta)}\);
\item \(\Gamma\) is rectifiable and the unit tangent direction \(\tau\) to \(\Gamma\) satisfies \(\tau(z)=e^{i u(z)}\) for some
real-valued function \(u \in\dot{H}^{\frac{1}{2}}(\Gamma)\).
\end{enumerate}

\end{theorem}

In \cite{Bishop1} and \cite{Bishop2}, C. Bishop gave a more geometric approach to the proof of Theorem \ref{shen charac}. He also showed the following characterization result:
\begin{theorem}
 With \(f\) and $\Gamma$ as in \emph{Theorem \ref{shen-liu}} one has
\begin{enumerate}[font=\upshape]
    \item \(\Gamma\) is a \emph{WP}-class quasicircle iff it has finite M\"obius energy, i.e.,
\begin{equation}
    \emph{Möb}(\Gamma):=\int_{\Gamma} \int_{\Gamma}\left(\frac{1}{|x-y|^{2}}-\frac{1}{\ell(x, y)^{2}}\right)\, d x \, d y<\infty, 
\end{equation}
where \(\ell(x, y)=\) arclength distance between \(x, y\) along curve $\Gamma$.

\item \(\Gamma\) is a \emph{WP}-class quasicircle iff it is chord-arc and the arclength parameterization is in the Sobolev space \(H^{\frac{3}{2}}(\mathbb{S}^1)\).

\end{enumerate}

\end{theorem}
Here, a rectifiable curve \(\Gamma\) is called {\it chord-arc} if for all \(x, y \in \Gamma\) one has  \(length(\gamma)=O(|x-y|)\)
where \(\gamma \subset \Gamma\) is the shortest sub-arc with endpoints \(x, y\).\\

{\begin{remark}
    The fact that the WP-class quasicircles are chord-arc curves was stated by the authors in the article with D. Radnell \cite{RSS_WP_multiply}, though the proof had a gap. Bishop provided a rigorous proof of that fact in \cite{Bishop1}. 
\end{remark}}

Finally, we also note that the tangent space to \(T_{\mathrm{WP}}(\mathbb{D}_-)\) at the identity consists precisely of the \(H^{\frac{3}{2}}\) vector fields on the unit circle with some normalization conditions (see \cite{NagVerjovsky}, \cite{Takhtajan_Teo_Memoirs}).

\end{subsection} 

\end{section}
\begin{section}{The period mapping}  \label{se:period}  
\begin{subsection}{Polarizations and the Siegel disk} \label{se:polarizations_Siegel} 
  A number of researchers have considered an embedding of the group $\text{Diff}(\mathbb{S}^1)/\text{M\"ob}(\mathbb{S}^1)$ in an infinite-dimensional version of the Siegel disk. Later this was extended to 
  the Teichm\"uller spaces $T_{\mathrm WP}(\disk_-)$ and $T(\disk_-)$.
  
  A paper of Nag and Sullivan established the analytic point of view appropriate for Teichm\"uller theory. To simplify the presentation, we adopt their point of view. Along the way we will fill in some of the history of the ideas.
  
 Let $H^{1/2}_{\mathbb{R}}(\mathbb{S}^1)$ and $H^{1/2}(\mathbb{S}^1)$ denote Sobolev $1/2$ space of real- and complex-valued functions on the circle respectively. Denote spaces with constants modded out with a dot, e.g. $\dot{H}^{1/2}_{\mathbb{R}}$. 
 It was independently shown by K. Vodopy'anov \cite{Vodopyanov} and S. Nag and D. Sullivan \cite{NagSullivan} that the composition by a homeomorphism $\phi$
 \[  h \mapsto  h \circ \phi  \]
 (mod constants) is a bounded operator on $\dot{H}^{1/2}(\mathbb{S}^1)$ if and only if $\phi$ is a quasisymmetry 
 (Vodopy'anov formulates this result in an equivalent form on the real line).
 One can show that $\dot{H}^{1/2}(\mathbb{S}^1)$
 is a symplectic space with respect to the completion of the pairing
  \begin{equation} \label{eq:symplectic_form_Honehalf}
     \omega(f,g) = \frac{1}{\pi} \int_{\mathbb{S}^1} f \cdot dg  
  \end{equation}
 for smooth functions. 
 Thus we obtain a representation of quasisymmetries by symplectomorphisms on $\dot{H}^{1/2}(\mathbb{S}^1)$. In summary
 \begin{theorem}[\cite{Vodopyanov,NagSullivan}] \label{th:symp_are_quasi}
  A homeomorphism $\phi:\mathbb{S}^1 \rightarrow \mathbb{S}^1$ induces a bounded composition operator on $\dot{H}^{1/2}(\mathbb{S}^1)$ if and only if $\phi \in \mathrm{QS}(\mathbb{S}^1)$. This bounded operator is a symplectomorphism.
 \end{theorem}

 The fact that $\mathrm{Diff}(\mathbb{S}^1)$ acts symplectically with respect to this pairing was already known. Segal \cite{SegalUnitary} investigated it in conjunction with a metaplectic representation of the restricted symplectic group. The action was there restricted to smooth real-valued functions on the circle mod constants (the $C^\infty$ elements of $\dot{H}^{1/2}_{\mathbb{R}}$).
  Segal defined an infinite-dimensional analogue of the Siegel disk, which parametrizes the positive polarizations of a complexified symplectic space.
  Nag and Sullivan \cite{NagSullivan} showed that the action has a natural extension to $\mathrm{QS}(\mathbb{S}^1)$, and initiated the analytic theory which is ultimately necessary to rigorously investigate the K\"ahler geometry of Teichmuller space, the period map, {and symplectic actions}. As we will see, the correct analytic leads to a clearer geometric and algebraic picture (in much the same way that the algebra of Fourier series in $L^2(\mathbb{S}^1)$ is clearer than the algebra of Fourier series in, say, $C^\infty(\mathbb{S}^1)$).  For example, it leads to the precise connection to the moduli spaces $T(\disk_-)$ and $T_{\mathrm{WP}}(\disk_-)$, as we will see shortly.\\
  
  We now describe the infinite Siegel disk more precisely.
  We say that $(V^\mathbb{C},\omega,J)$ is a complex symplectic Hilbert space if the following hold:
  \begin{enumerate} 
   \item $V^\mathbb{C}$ is the complexification of a real vector space $V$ with complex structure;
   \item $J$ and $\omega$ are the complex linear extensions of the complex structure and a symplectic form on $V$;
   \item $V^{\mathbb{C}}$ is a separable Hilbert space with respect to the pairing $\left<v,w \right> := \omega(v, J\overline{w})$;
   \item $\omega$ and $J$ are continuous. 
  \end{enumerate} 
   \begin{definition}  
  Let $(V^{\mathbb{C}},J,\omega)$ be a complex symplectic Hilbert space. We say that a complex subspace $W$ of $V^{\mathbb{C}}$ is a positive polarizing subspace if 
  \begin{enumerate}
      \item $W$ is isotropic;
      \item $V^{\mathbb{C}} = W \oplus \overline{W}$; and
      \item $\omega(v,i\overline{w})$ is a positive-definite sesquilinear form on $W$.  
  \end{enumerate}
  We also refer to the decomposition (2) as a polarization.
  \end{definition}
  It is easily verified that if $W_0$ and $\overline{W_0}$ are the $-i$- and $i$-eigenspaces of $J$, then $W_0 \oplus \overline{W_0}$ is a polarization, which we call the standard polarization.  

  The set of polarizations is parametrized by the infinite Siegel disk, which we now define. 
  \begin{definition} Let $(V^{\mathbb{C}},J,\omega)$ be a complex symplectic Hilbert space and let $V^{\mathbb{C}} = W_0 \oplus\overline{W_0}$ denote the standard polarization. The infinite Siegel disk is the set of bounded linear maps $Z:\overline{W_0} \rightarrow W_0$ satisfying
  \begin{enumerate}
   \item $\omega(v,Zw)=\omega(w,Zv)$ for all $v,w \in \overline{W_0}$; 
   \item $I-\overline{Z} Z$ is positive definite.
  \end{enumerate}
  The restricted Siegel disk is the set of elements of the infinite Siegel disk which in addition satisfy 
  \begin{enumerate} \setcounter{enumi}{2}
   \item $Z$ is Hilbert-Schmidt.
  \end{enumerate}
  \end{definition}
  \begin{remark}
  Condition (1) is equivalent to $Z^T=Z$ in an orthonormal basis, or alternatively ${Z}^*=\overline{Z}$.
  \end{remark}
  \begin{remark} Although we have not seen the term ``restricted Siegel disk'', it is an obvious extension of terminology in use for the Shale group. In the literature the term ``Siegel disk'' sometimes agrees with the terminology here and sometimes is used for the restricted Siegel disk.    
  \end{remark}
  It can be shown that the set of positive polarizations is in one-to-one correspondence with the infinite Siegel disk, by setting $W$ to be the graph of $Z$.  
  It is also easy to see that the group of real symplectomorphisms of $V$ acts on the space of polarizations by complex linear extension, and thus act on the Siegel disk.\\

  The Siegel disk was considered first by C. L. Siegel in finite dimensions \cite{Siegel_symplectic_book}. Segal \cite{SegalUnitary} considered the restricted Siegel disk in association with the representations of the symplectic group induced by $\mathrm{Diff}(\mathbb{S}^1)$ as described above.
  The condition that $Z$ be Hilbert-Schmidt, included by Segal, is derived from Shale \cite{Shale}. Shale showed it to be necessary and sufficient that an automorphism of the CCR algebra be implementable as a unitary operator on Fock space - that is, the necessary condition to produce the projective unitary representation associated to the metaplectic group. The corresponding subgroup $\mathrm{Sym}_{\mathrm{res}}(V)$ of the symplectic group of a vector space $V$ is called the restricted symplectic group, and appears throughout the representation theory of infinite-dimensional groups, see e.g. J. Ottesen \cite{Ottesen} or M. Schottenloher \cite{Schottenloher} in addition to the above references.\\ 

  A standard example of a polarization is given by the holomorphic and anti-holomorphic decomposition of one-forms on a Riemann surface. Indeed this is the basis of the classical period map, and was a motivation for Siegel's construction. In the literature this is usually given in terms of the equivalent Siegel upper half-plane model.  
  
  Another example is provided by the positive and negative Fourier modes of functions on the circle \cite{SegalUnitary}, which we provide here in the setting of Nag and Sullivan. Set $V^{\mathbb{C}} = \dot{H}^{1/2}(\mathbb{S}^1) \subset L^2(\mathbb{S}^1)$, where the symplectic form is the completion of (\eqref{eq:symplectic_form_Honehalf}) and the complex structure is the Hilbert transform 
  \[  J e^{i n \theta} = -i \, \mathrm{sgn}(n) \, e^{i n \theta}  \]
  for $n \in \mathbb{Z} \backslash \{0\}$. The standard polarization is given by the functions with only positive and only negative Fourier modes. These are boundary values of holomorphic and anti-holomorphic functions with square integrable derivatives on the disk (the Dirichlet spaces $W_0=\dot{\mathcal{D}}(\disk_+)$, $\overline{W_0}=\overline{\dot{\mathcal{D}}(\disk_+)}$), that is   
  \[    \dot{H}^{1/2}(\mathbb{S}^1) = \dot{\mathcal{D}}(\disk_+) 
  \oplus \overline{\dot{\mathcal{D}}(\disk_+)}. \]
  
  It is easily shown that $\text{M\"ob}(\mathbb{S}^1)$ is the precise subset of $\text{QS}(\mathbb{S}^1)$ preserving the standard polarization. 
  Thus by Theorem \ref{th:symp_are_quasi} of Vodopy'anov, Nag-Sullivan, we obtain an embedding of the universal Teichm\"uller space into the infinite Siegel disk. 
  This was considered earlier for the smooth subspace $\mathrm{Diff}(\mathbb{S}^1)$ by many authors, including Segal \cite{SegalUnitary}  as mentioned above, Kirillov-Yuri'ev \cite{KY2}, and Nag \cite{Nag_bulletin}. It was called the KYNS (Kirillov-Yuri'ev-Nag-Sullivan) embedding by Takhtajan-Teo \cite{Takhtajan_Teo_Memoirs}\footnote{though Segal-Kirillov-Yuri'ev-Nag-Sullivan-Takhtajan-Teo embedding might be more appropriate.}. \\
  
  The finite-dimensional Siegel disk is a homogeneous complex manifold with a natural symplectically invariant Hermitian metric, which at the identity is given by $Tr (Z_1 \bar{Z}_2)$ \cite{Siegel_symplectic_book}.  (Note that by the symmetry of the $Z$ matrices, this is the Frobenius inner product for matrices).  This naturally generalizes to the infinite-dimensional restricted Siegel disk, because the condition for its convergence is precisely that $Z_1$ and $Z_2$ are Hilbert-Schmidt.\\ 

  It's not hard to show that diffeomorphisms induce Hilbert-Schmidt composition operators \cite{SegalUnitary}. 
  It was shown by Kirillov-Yur'iev \cite{KY2}, and later Nag and Sullivan \cite{Nag_bulletin}, that the pull-back of the generalized Siegel K\"ahler form under the KYNS period embedding agrees with the Weil-Petersson metric for $\text{Diff}(\mathbb{S}^1)/\text{M\"ob}(\mathbb{S}^1)$. The full result for $T_{\mathrm{WP}}(\disk_-)$ is due to Takhtajan-Teo.\\
  
   Now the operator $Z$ in the KYNS embedding was shown to be the Grunsky matrix by Kirillov-Yuri'ev \cite{KY2}.
   The Grunsky matrix takes many forms in the literature. One of these, due to Bergman and Schiffer \cite{BergmanSchiffer}, is as an integral operator on the anti-holomorphic Bergman space.  Given an element of $T(\disk_-)$ let $f:\mathbb{D}_+ \rightarrow \sphere$ be the associated normalized conformal map given by Bers trick (\ref{eq:Bers_trick_temp}). We associate a Grunsky operator to $f$ (and hence to an element of $T(\disk_-)$) as follows. For square integrable antiholomorphic $\overline{h(z)}$ we define 
   \begin{equation}  \label{eq:Grunsky_integral}  
   {\bf{Gr}}_f \overline{h}(z) 
     = \frac{1}{\pi} \iint_{\disk^+}   \left( \frac{f'(w) f'(z)}{(f(w)-f(z))^2} - \frac{1}{(w-z)^2} \right) \overline{h(w)}\,\frac{d\bar{w} \wedge dw}{2i}.   
  \end{equation} 
  Since differentiation is an isometry from the Dirichlet space mod constants to the Bergman space, one can equivalently formulate this on the homogeneous Dirichlet space. In this case we get 
  \begin{equation} \label{eq:Grunsky_Schiffer_connection} 
   Z = \Pi([\mu]) = d^{-1} \mathbf{Gr}_f d: \overline{\dot{\mathcal{D}}(\disk_+)} \rightarrow \dot{\mathcal{D}}(\disk_+).    
  \end{equation}
   From the integral expression (\ref{eq:Grunsky_integral}), it is easily seen that this operator is M\"obius invariant, that is $\mathbf{Gr}_{T \circ f}=\mathbf{Gr}_f$ for $T$ M\"obius, so that the normalization of $f$ is irrelevant as expected.  
   {\begin{remark}  The Grunsky matrix can be defined for holomorphic functions $f$ which are one-to-one in a neighbourhood of $0$ \cite{Durenbook,Pommerenkebook}.  In \cite{KY2}, there is a potentially very misleading typographical error. It is stated that $I- \overline{Z}_f Z_f>0$ is equivalent to the condition that the locally conformal map $f$ is univalent on the disk. In fact, it is classically known that, rather, the fact that $I-\overline{Z} Z \geq 0$ when applied to vectors in $\mathbb{C}^n$ for any $n$ is equivalent to univalence on $\disk$. (In the formulation here, $I-\overline{Z} Z \geq 0$ as a linear transformation on polynomials.) The condition that $I-\overline{Z} Z > 0$  is equivalent to univalence and quasiconformal extendibility by a result of Pommerenke and K\"uhnau.  See for example Pommerenke \cite{Pommerenkebook}. 
   \end{remark}
   }
   
The following theorem, independently obtained by Y. Shen and Takhtajan-Teo, strongly motivates Weil-Petersson class Teichm\"uller theory. 
   \begin{theorem}[\cite{ShenGrunsky,Takhtajan_Teo_Memoirs}] \label{th:symp_res_is_WP}
    The Grunsky operator induced by $\phi \in \mathrm{QS}(\mathbb{S}^1)$ is Hilbert-Schmidt if and only if $\phi$ is a Weil-Petersson class quasisymmetry. In particular, a homeomorphism induces a composition operator on $\dot{H}^{1/2}(\mathbb{S}^1)$ which is in the restricted symplectic group if and only if it is a Weil-Petersson class quasisymmetry.
   \end{theorem} 
   Thus $T_{\mathrm{WP}}(\mathbb{D}_-)$ is exactly the group of elements of $T(\disk_-)$ which map into the restricted Siegel disk under the period map $\Pi$. Denote the restriction of $\Pi$ to $T_{\mathrm{WP}}(\mathbb{D}_-)$ by $\Pi_{WP}$.  

   The hierarchy of conditions is striking, in light of the geometric interpretation. $I- \overline{Z}_f Z_f \geq 0$ implies univalence of $f$, but $I- \overline{Z}_f Z_f>0$ implies in addition that $f$ is  quasiconformally extendible (and thus $f$ is associated to an element of the universal Teichm\"uller space). Adding the further condition that $Z$ be Hilbert-Schmidt then restricts to the Weil-Petersson class Teichm\"uller space.  
   
   \begin{remark}
   Segal-Wilson and Kirillov-Yuri'ev consider the action of $\mathrm{Diff}(\mathbb{S}^1)$ on $L^2(\mathbb{S}^1)$, so that the extensions are in fact in the direct sum of the Hardy space and anti-holomorphic Hardy space of the disk. This is valid for smooth homeomorphisms, but quasisymmetries do not appear to be bounded on $L^2(\mathbb{S}^1)$. For quaisymmetries one has to consider the space $\dot{H}^{1/2}(\mathbb{S}^1)$ space as Nag and Sullivan do.  
   \end{remark}

One may ask whether it is possible to recover the Weil-Petersson metric on finite-dimensional spaces from that on $T_{\mathrm{WP}}(\disk_-)$. 
Nag and Verjovsky \cite{NagVerjovsky} showed that the Teichm\"uller space of any Fuchsian group is transverse to the fibre. In particular, the pull-back of the Weil-Petersson metric on $T_{\mathrm{WP}}(\disk_-)$ to the finite-dimensional Teichm\"uller spaces cannot be the Weil-Petersson metric of the finite-dimensional space. However, they show how an averaging procedure due to S.J. Patterson recovers the Weil-Petersson metric on the finite-dimensional spaces from that on $\text{Diff}(\mathbb{S}^1)/\text{M\"ob}(\mathbb{S}^1)$. Takhtajan-Teo use this technique to obtain Wolpert's curvature formulas from their curvature formulas on $T_{\mathrm{WP}}(\disk_-)$ \cite{Takhtajan_Teo_Memoirs} (see Section \ref{se:Kahler_and_global} ahead).
 \end{subsection}
 \begin{subsection}{Interpretation of the period map} 
 
  Among a number of interpretations, Nag and Sullivan \cite{NagSullivan} observed that this period map is analogous to the period map for compact surfaces. 
 
  {Historically, B. Riemann \cite{Riemann} used the period matrices and his bilinear relations to provide a necessary and sufficient condition for a set of $2g$ linearly independent (over $\mathbb{R}$) vectors in $\mathbb{C}^{g}$ to be periods of $g$ holomorphic differentials on a $g$-dimensional Abelian variety. He also solved the {\it Jacobi inversion problem} for an arbitrary algebraic curve (using his Theta functions). Riemann's work was developed further by Torelli \cite{Torelli}, where biholomorphic equivalence of compact Riemann surfaces of genus $g\geq 1$ was connected to the so-called {\it polarized Jacobian} of the Riemann surfaces in question.}  
{The work of Torelli was developed further by A. Andreotti and A. Weil, and turned into one of the cornerstones of complex algebraic geometry.}\\

We begin with holomorphicity of the period map.
Ahlfors proved that the period map of a compact surface $R$ depends holomorphically on the surface, and that the complex structure of $T(R)$ is the unique one such that this holds \cite{Ahlfors_period}. Takhtajan-Teo \cite{Takhtajan_Teo_Memoirs} showed that both period maps $\Pi$ and $\Pi_{\mathrm{WP}}$ are holomorphic. 
  The codomain of $\Pi$ is 
  the Banach space $\mathfrak{B}(\overline{\mathcal{D}(\disk_-)} \rightarrow \mathcal{D}(\disk_+)$ of bounded operators from $\overline{\mathcal{D}(\disk_-)}$ to $\mathcal{D}(\disk_+)$. We have the following.
  \begin{theorem}[\cite{Takhtajan_Teo_Memoirs}] 
   \label{th:TT_period_Banach_holomorphic}
   The period      embedding $\Pi:T(\disk_-) \rightarrow \mathfrak{B}(\overline{\dot{\mathcal{D}}(\disk_-)},\dot{\mathcal{D}}(\disk_-))$
   is injective and holomorphic.
  \end{theorem}
  \begin{remark}
      Their statement is in terms of $\ell^2$. The formulation here is an easy consequence of equation (\ref{eq:Grunsky_Schiffer_connection}), after identifying $\ell^2$ with the Bergman space.  We altered the formulation for consistency with Nag and Sullivan \cite{NagSullivan} and to offer a geometric interpretation ahead.
  \end{remark}
   The codomain of $\Pi_{\mathrm{WP}}$, the restricted Siegel disk, is a subset of the {Sato}-Segal-Wilson Grassmannian $\mathfrak{S}$. This possesses a complex Hilbert manifold structure induced by the Hilbert-Schmidt norm.  We have 
  \begin{theorem}[\cite{Takhtajan_Teo_Memoirs}]  \label{th:TT_period_WP_holomorphic}
   The \emph{KYNS} period embedding
   \[ \Pi_{\mathrm{WP}}: T_{\mathrm{WP}}(Z) \rightarrow \mathfrak{S}  \]
   is holomorphic. 
  \end{theorem}

  The authors showed that the boundary values of an element of the homogeneous Dirichlet space of a quasidisk $\kap$ are in turn boundary values of an element of the homogeneous Dirichlet space of its complement. This ``overfare'' is bounded precisely for quasidisks. This allows the following result interpreting the polarizations of Nag and Sullivan. 
  \begin{theorem}[\cite{RadnellSchippersStaubach_Grunsky_genuszero}]  \label{th:graph_of_Grunsky}
    Let $f$ be associated to an element $[\mu] \in T(\disk_-)$ and $Z= \Pi([\mu])$. Let $\kap = f(\disk_+)$ and $\kap^* = \sphere \backslash \mathrm{cl} (f(\disk_+))$.   Then the graph $W$ of $Z$ is 
    \[  \{ h:\disk_+ \rightarrow \mathbb{C} \text{ harmonic} \,: \, \left. h \right|_{\mathbb{S}^1} = H \circ f \text{ for some } H \in \dot{\mathcal{D}}(\kap^*) \}.  \]
  \end{theorem}
  It is understood in the above that $H$ extends to $\partial \kap^*=\partial \kap$, so that this composition makes sense (there are analytic subtleties, which we suppress for the sake of brevity).
  In other words, the graph $W$ of the Grunsky operator $Z$ is the pull-back of the boundary values of Dirichlet-bounded holomorphic functions on $\kap^*$. {Thus the new polarization  can be identified with the pull-back of the polarization on $\kap^*$. Identifying $\kap^*$ with an element of Teichm\"uller space, this is the interpretation of the period map given in Nag and Sullivan \cite{NagSullivan} (after conjugating by $d$ to rephrase in terms of forms).   }\\
  
  The interpretation of the polarizations as the pull-back of boundary values of holomorphic functions also arises in two-dimensional conformal field theory \cite[Appendix D]{Huang}. 
  Radnell-Schippers-Staubach discuss the relation of this result to conformal field theory in \cite{RadnellSchippersStaubach_CFT_review}; however there the theorem was stated only for Grunsky matrices associated to elements of the Weil-Petersson Teichm\"uller space. Later the authors discovered that, surprisingly, this extends to the entire Teichm\"uller space \cite{RadnellSchippersStaubach_Grunsky_genuszero}. In fact it was extended there to bordered surfaces of type $(0,n)$.\footnote{The earlier result for WP-class quasicircles was given in \cite{RSS_WP_multiply}, which was left unpublished.} The theorem for $(0,n)$ is as follows in the notation used here. Let $f \in \mathcal{O}^{\mathrm{qc}}(\sphere \backslash \{p_1,\ldots,p_n\})$ be a collection of quasiconformal maps with non-overlapping images, each taking $0$ to $p_n$. Denote the induced Grunsky operator by 
  \[  Z_f :\overline{\dot{\mathcal{D}}(\disk_+)}^n \rightarrow   \dot{\mathcal{D}}(\disk_+)^n   \]
  (this can be defined using the integral expression (\ref{eq:Grunsky_integral}) for example). Setting 
  \[  \riem = \sphere \backslash \cup_{k=1}^n \mathrm{cl} \, f_k(\disk_+)   \]
  then have that, denoting $\mathbb{D}_+^n = \mathbb{D}_+  \times \cdots \times \mathbb{D}_+$, the graph of the generalize Grunsky operator is  
  \[  \{ h:\disk_+^n \rightarrow \mathbb{C} \text{ harmonic} \,: \, \left. h \right|_{\left(\mathbb{S}^1\right)^n} = H \circ f \text{ for some } H \in \dot{\mathcal{D}}(\riem) \}.  \]
  and Theorem \ref{th:graph_of_Grunsky} is a special case.   
  
  In light of the above it was natural to ask whether it is possible to use the fiber/CFT model to define period matrices in such a way that the finite-dimensional period map and infinite-dimensional period map are naturally unified.  The next step, extending Grunsky operators to arbitrary genus surfaces and number of boundary curves, was accomplished by M. Shirazi \cite{Shirazi_Grunsky}, and Theorem \ref{th:graph_of_Grunsky} was shown to hold in that setting. 
  There, adjustments must be made for topological obstructions arising from the classical periods.  The relation of the generalized period mappings to the classical period matrices for compact surfaces beyond an analogy was an open question, but recent work \cite{Schippers_Staubach_scattering} {and work in progress} unifies them. 

\begin{remark} 
 In higher genus one needs to work with one-forms and not just functions. To this end the authors showed that the boundary values of $L^2$ harmonic one-forms can be identified with the Sobolev $H^{-1/2}$ space in a conformally invariant way, and developed a theory of overfare of one-forms \cite{Schippers_Staubach_scattering,Schipp-Staub-Kenig}.
\end{remark}

 It is thus of interest to ask whether the extended period mappings are holomorphic maps of Teichm\"uller space. Radnell-Schippers-Staubach \cite{RSS_period_genuszero} showed that Theorem \ref{th:TT_period_Banach_holomorphic} holds for surfaces of type $(0,n)$, using the atlas arising from the fiber model. (Note that there the exceptional cases $g=0, n=1$ and $g=0, n=2$ are taken care of as in Remark \ref{re:fiber_exceptional_cases}.) The methods in that paper should extend fairly easily to type $(g,n)$. A generalization of Theorem \ref{th:TT_period_WP_holomorphic} to the Weil-Petersson Teichm\"uller space  of type $(g,n)$ bordered surfaces (described in Section \ref{se:WP_higher_genus} ahead) is highly desirable. 
\end{subsection}     
\end{section}

\begin{section}{A brief overview of other refinements of Teichm\"uller space}\label{birds eye view}

Here, for the sake of completeness and comparison, we also briefly mention some other refinements of the Teichm\"uller space other than the Weil-Petersson.\\

{\underline{Gardiner-Sullivan's asymptotically conformal Teichm\"uller space.}}
 The asymptotic Teichmüller space \(\mathrm{A T}(\mathscr{R})\)
was introduced by F. Gardiner and D. Sullivan \cite{Gardiner-Sullivan} for the upper half-plane, and by
C. Earle, F. Gardiner and N. Lakic \cite{EarleGardinerLakic1} for arbitrary hyperbolic Riemann surfaces. A univalent function \(f\) in the unit disc $\mathbb{D}$ admiting a quasiconformal extension to the whole plane, with \(f(\infty)=\infty\) is called {\it{asymptotically conformal}} \cite{Pommerenke}
  if $$
\lim_{|z| \rightarrow 1^{-}}(1-|z|)\left|\frac{f^{\prime \prime}(z)}{f^{\prime}(z)}\right|= 0. 
$$
  Now given the two quasiconformal mappings \(f\) and \(g\) on a hyperbolic Riemann surface $\mathscr{R}$, consider the following equivalence relation:\\
  $f\sim g$, if  there is an
asymptotically conformal mapping \(h\) of \(f(\mathscr{R})\) onto \(g(\mathscr{R})\) such that \(g^{-1} \circ h \circ f\)
is homotopic to the identity relative to the ideal boundary of \(\mathscr{R}\).
        The {\it {asymptotic
Teichmüller space}} \(\mathrm{A T}(\mathscr{R})\) is the set of equivalence classes of quasiconformal mappings of \(\mathscr{R}\) under this relation. 
       
Note that this definition is a variation of the definition of
Teichmüller space \({T}(\mathscr{R})\) where the mapping \(h\)
is required to be conformal.
Since conformal mappings are asymptotically conformal, there is a well-defined
projection $\mathcal{P}: T(\mathscr{R}) \rightarrow \mathrm{AT}(\mathscr{R})$. 
If $\mathscr{R}$ is a closed Riemann surface of finite genus with a finite number of points removed, then
$\mathrm{AT}(\mathscr{R})$ consists of one point.
Earle, Gardiner and Lakic \cite{EarleGardinerLakic1, EarleGardinerLakic2} proved that $(\mathrm{A T}(\mathscr{R})$ has a complex manifold
structure
so that the quotient map \(\mathcal{P}: T(\mathscr{R}) \rightarrow \mathrm{A T}(\mathscr{R})\) is holomorphic. For more information on this topic see Gardiner-Lakic's monograph \cite{GardinerLakic}. \\

      {\underline{Astala-Zinsmeister's BMO Teichm\"uller space.}}
         This version of the Teichm\"uller space was introduced by K. Astala and M. Zinsmeister in \cite{AstZins}. To define it requires some preparation.
         First recall that a positive measure \(\lambda\) defined in a simply connected domain \(\kap\) is called a {\it Carleson measure} if
\begin{equation}
    \label{Carleson norm}
    \|\lambda\|_{c}=\sup \left\{\frac{\lambda(\kap \cap \mathbb{D}(z, r))}{r}: z \in \partial \kap, 0<r<\operatorname{diameter}(\kap)\right\}<\infty,
\end{equation}

where \(\mathbb{D}(z, r)\) is the disk with center \(z\) and radius \(r\). By \(C M(\kap)\) one denotes the set
of all Carleson measures on \(\kap\).

Now let  \(\mathcal{L}_\mathrm{BMO}(\mathbb{D}_+)\) be the Banach space of essentially bounded measurable functions \(\mu\) on \(\mathbb{D}_+\) such that the
measure
$$
\lambda_{\mu}=\frac{|\mu|^{2}(z)}{1-|z|^{2}}\, d\mathrm{A}_z
$$
is in \(C M(\mathbb{D}_+)\). The norm on \(\mathcal{L}_\mathrm{BMO}(\mathbb{D}_+)\) is defined by
$$\Vert \mu\Vert_{\mathcal{L}_\mathrm{BMO}(\mathbb{D}_+)}
:= \Vert \mu\Vert_{\infty}+\Vert \lambda_{\mu} \Vert_{c}^{1 / 2}
$$
where \(\left\|\lambda_{\mu}\right\|_{c}\) is the Carleson norm of \(\lambda_{\mu}\) defined in \eqref{Carleson norm}. Set
$$\mathcal{M}_{\mathrm{BMO}}(\mathbb{D}_+)=\left\{\mu \in \mathcal{L}_{\mathrm{BMO}}(\mathbb{D}_+):\|\mu\|_{\infty}<1\right\}.$$
Cui-Zinsmeister \cite{CuiZinsmeister} {extended this to surfaces associated with a Fuchsian group $G$ as follows}. Let \(G\) be a Fuchsian group, i.e. a properly discontinuous fixed point free group of Möbius transformations
which keeps \(\mathbb{D}_+\) invariant. For such a group set
\begin{equation}
    \label{defn: MG}
    M(G):=\left\{\mu \in L^{\infty}(\mathbb{D}_+):\|\mu\|_{\infty}<1 \text { and } \forall g \in G, \mu=\mu \circ g \frac{\bar{g}^{\prime}}{g^{\prime}}\right\}.
\end{equation}
Now for the corresponding Riemann surface \(\mathscr{R}=\mathbb{D}_+ / G\), define $$\mathcal{M}_{\mathrm{BMO}}(G):=M(G) \cap \mathcal{M}_{\mathrm{BMO}}(\mathbb{D}_+)$$ with
the same equivalence relation as in the classical Teichm\"uller space. Given $\mathcal{M}_{\mathrm{BMO}}(G)$, the BMO-{\it{Teichm\"uller space}} denoted here by $\mathrm{BMOT}(\mathscr{R}),$ is defined as the
quotient space associated to the equivalence relation.

Note that if \(G\) is co-compact then, as for the
asymptotic Teichmüller spaces, $\mathrm{BMOT}(\mathscr{R})$ is trivial. But this is the extent of the analogy and the
equivalence relation in the case of asymptotic Teichmüller space is rather different. For further developments in this context the reader is referred to \cite{MatsWei2, {WeiZins}}. 
\\
         
{\underline{Guo-Tang's $L^p$-Teichm\"uller space.}}
For $p\in [1, \infty)$, the $L^p$-Teichm\"uller spaces were first introduced by Guo in \cite{GuoHui}, and developed by Tang \cite{Tang}, are defined as follows. First let  \(\mathcal{L}^p(\mathbb{D}_+)\) be the Banach space of essentially bounded measurable functions \(\mu\) on \(\mathbb{D}_+\) such that the
$$
\iint_{\mathbb{D}_+}\frac{|\mu|^{p}(z)}{(1-|z|^{2})^2}\, d\mathrm{A}_z <\infty.
$$       
 Set $\mathcal{M}_{p}(\mathbb{D}_+)=\left\{\mu \in \mathcal{L}^p(\mathbb{D_+}):\|\mu\|_{\infty}<1\right\},$ and for a Fuchsian group $G$ define
 \begin{equation}
  \mathcal{M}_p(G):=M(G) \cap \mathcal{M}_{p}(\mathbb{D}_+),   
 \end{equation}
where $M(G)$ is as in
     \eqref{defn: MG}, and with the same equivalence relation as in the classical Teichm\"uller space. For the Riemann surface Riemann surface \(\mathscr{R}=\mathbb{D} _+/ G\), the $L^p$-{\it{Teichm\"uller space}} denoted here by $\mathrm{L^p T}(\mathscr{R}),$ is defined as the
quotient space associated to the equivalence relation.  In \cite{Yanagishita}, Yanagishita  introduced, for any $p\geq 2$, a complex structure on the $G$-invariant $p$-integrable Teichm\"uller space associated to an arbitrary Fuchsian group $G$ satisfying the Lehner's condition.

For further developments in this context the reader is referred to \cite{shenliu, TangShen, MatsYana1, MatsWei}. 
The case of $p=2$ is of course the WP-{class Teichm\"uller space.}
  
\end{section}

\begin{section}{Weil-Petersson Teichm\"uller spaces of general surfaces}  \label{se:WP_higher_genus} 
\begin{subsection}{Overview}
 The Weil-Petersson Teichm\"uller theory was extended to more general surfaces by Yanagishita and Radnell-Schippers-Staubach independently.  The approaches are somewhat different, so that the results largely complement each other. In the end, the complex structures in all three models are shown to exist and be equivalent. We give a brief outline here. The precise results are given in the next section. \\ 
 
 The first step was taken by Radnell-Schippers-Staubach in \cite{RSS_Filbert_arxiv}, later published in two parts \cite{RSS_Filbert1,RSS_Filbert2}.  There, a definition of Weil-Petersson class Teichm\"uller spaces of type $(g,n)$ was given, by restricting to equivalence classes of quasiconformal deformations whose boundary values are Weil-Petersson class quasisymmetries.  The topology and a Hilbert manifold structure was given using the fiber model.\\ 

 Yanagishita \cite{Yanagishita} and Radnell-Schippers-Staubach \cite{RSS_WP_arxiv} (later published as \cite{RSS_WP1,RSS_WP2})  simultaneously and independently developed two further aspects. Yanagishita gave a complex structure on the Weil-Petersson Teichm\"uller space for surfaces satisfying what he calls Lehner's condition, which puts a lower bound on the length of the simple closed geodesics.  This applies to bordered surfaces of type $(g,n)$. His approach is similar to the approach of Cui \cite{Cui}, using the Bers embedding model of Section \ref{subsubse:Bers_embedding_complex}. In fact, as in Guo \cite{GuoHui} and Tang \cite{Tang} in the case of the universal Teichm\"uller space, the results are shown to hold for the $L^p$ theory for $p\geq 1$ (which in general results in a Banach manifold structure). The case $p=2$ is the WP space and one obtains a Hilbert manifold structure. 

On the other hand, Radnell-Schippers-Staubach \cite{RSS_WP_arxiv}, later published in two parts as \cite{RSS_WP1,RSS_WP2}, gave the Weil-Petersson Teichm\"uller space a complex Hilbert manifold structure based on harmonic Beltrami differentials; that is, based on the deformation model of Section \ref{subsubse:deformation_complex_structure}, for type $(g,n)$ bordered surfaces with $2g-2+n>0$. Two definitions were advanced for this Weil-Petersson Teichm\"uller space: that the boundary values of the quasiconformal map are WP-class quasisymmetries, and that there is a representative with $L^2$ Beltrami differential. The latter definition agrees with Yanagishita's definition restricted to type $(g,n)$ surfaces.  Using sewing techniques and boundary behaviour of the hyperbolic metric, it was shown that the two definitions are equivalent. The deformation complex structure was also shown to be compatible with the complex Hilbert manifold structure obtained from the fiber model. Note that even in the classical $L^\infty$ case, the compatibility of the complex structures derived from the fiber/CFT point of view and any of the classical structures was a fairly recent result \cite{RadnellSchippers_fiber}. 
As a consequence, for type $(g,n)$ bordered surfaces with $2g-2+n>0$,  Radnell-Schippers-Staubach \cite{RSS_WP_arxiv} gave the tangent space structure and show that the Weil-Petersson pairing converges on each tangent space.\\  

Yanagishita \cite{Yanagishita_Kahler} later independently gave the construction of the complex Hilbert manifold structure based on harmonic Beltrami differentials, for surfaces more generally satisfying Lehner's condition. He then showed that this is compatible with the complex structure derived from the Bers embedding, and that the Weil-Petersson metric converges on each tangent space. He furthermore showed that the metric is K\"ahler, generalizing the theorem of Takhtajan-Teo from the case of the universal Teichm\"uller space.  The results relating to K\"ahlericity and curvature will be discussed in Section \ref{se:Kahler_and_global} ahead.\\ 

 In summary, for surfaces satisfying Lehner's condition, the Bers embedding model and deformation model give equivalent complex structures. For bordered surfaces of type $(g,n)$ $2g-2+n>0$, the deformation model and fiber model give equivalent complex structures.
 Thus in the special case of type $(g,n)$ bordered surfaces with $2g-2+n>0$, all three complex structures are equivalent. 
 \begin{remark}
 One can then ask whether there is a meaningful sense in which the fiber model extends to general surfaces satisfying Lehner's condition (or perhaps some larger subset than type $(g,n)$ bordered surfaces). 
 \end{remark}
\end{subsection}
\begin{subsection}{Manifold and tangent space structure}
\begin{subsubsection}{{\underline{Fiber model}}}

The Weil-Petersson Teichm\"uller spaces for bordered surfaces of type $(g,n)$ was first defined and given a complex Hilbert manifold structure by Radnell- Schippers-Staubach \cite{RSS_Filbert_arxiv}, later appearing in \cite{RSS_Filbert1,RSS_Filbert2}. This was obtained by refining the fiber model for $T(\riem)$ to Weil-Petersson class mappings. Recall that one obtains a genus $g$ surface $\riem^P$ with $n$ punctures by sewing on caps. The idea is to restrict the construction of Section \ref{se:fiber_model_local_coordinates}) to more regular Weil-Petersson class maps. 
The local model of the conformal maps is given by 
\[ \mathcal{O}^{\mathrm{qc}}_{\mathrm{WP}} = \{ f \in \mathcal{O}^{\mathrm{qc}} : \| f''/f' \|_{A_1^2(\disk_+)} < \infty  \}.     \]
We thus define 
\begin{align*}
    \chi_{\mathrm{WP}} = \left. \chi \right|_{\mathcal{O}^{\mathrm{qc}}_{\mathrm{WP}}}:  \mathcal{O}^{\mathrm{qc}}_{\mathrm{WP}} & \mapsto A_1^2(\disk_+) \oplus \mathbb{C} \\
    f & \mapsto \left( f''/f',f'(0) \right).
\end{align*}
Since inclusion in $A_1^2(\disk_+) \rightarrow A_1^\infty(\disk_+)$ is bounded by Theorem \ref{th:disk_inclusion_bounded} of Cui and Takhtajan-Teo, it follows directly that the inclusion $\mathcal{O}_{\mathrm{WP}}^{\mathrm{qc}}$ is bounded.   
We have
\begin{theorem}[\cite{RSS_Filbert1}]
 Let $\riem^P$ be a surface of genus $g$ with $n$ punctures.
  $\chi_{\mathrm{WP}}(\mathcal{O}^{\mathrm{qc}}_{\mathrm{WP}})$ is open in $A_1^2(\disk_+) \oplus \mathbb{C}$.  
\end{theorem}

We say that $f \in \mathcal{O}^{\mathrm{qc}}_{\mathrm{WP}}(\riem^P)$ if $f =(f_1,\ldots,f_n) \in \mathcal{O}^{\mathrm{qc}}(\riem^P)$ (see Definition \ref{de:standard_Oqc}) and furthermore, there are coordinates
$\zeta_k:B_k \rightarrow \mathbb{C}$ of $p_k$ for $k=1,\ldots,n$ such that $\zeta_k \circ f_k \in \mathcal{O}_{\mathrm{WP}}^{\mathrm{qc}}$. 

We construct a coordinate chart in $\mathcal{O}^{\mathrm{qc}}_{\mathrm{WP}}(\riem^P)$ by restricting the open sets (\ref{eq:V_open_sets_definition}) and charts (\ref{eq:T_coordinates_definition}) of Section \ref{se:fiber_model_local_coordinates} to the Weil-Petersson case. With notation as in that section, set
\[  V_{\zeta,K} = \{ f=(f_1,\ldots,f_n) \in \mathcal{O}^{\mathrm{qc}}_{\mathrm{WP}}(\riem^P) : \mathrm{cl} (\, f_k(\disk_+)) \subset K_k, k=1,\ldots,n \}.    \]
These sets form a base for a Hausdorff, secound countable topology on $\mathcal{O}_{\mathrm{WP}}^{\mathrm{qc}}$ \cite{RSS_Filbert_arxiv}. 
We then define 
\begin{align*}  
 \mathcal{E}:V_{\zeta,K} & \rightarrow (\mathcal{O}^{\mathrm{qc}}_{\mathrm{WP}})^n \\
 f & \mapsto \zeta \circ f.
\end{align*} 
\begin{theorem}[\cite{RSS_Filbert1}]
 Let $\riem^P$ be a surface of genus $g$ with $n$ punctures.
  $\mathcal{O}^{\mathrm{qc}}({\riem}^P)$ is a Banach manifold  locally modelled on $\bigoplus^n A_1^\infty(\disk_+) \oplus \mathbb{C}$, with respect to the atlas $\chi_{\mathrm{WP}}^n \circ \mathcal{E}$.  
\end{theorem}
We also have
\begin{theorem}[\cite{RSS_Filbert1}] 
 Let $\riem^P$ be a surface of genus $g$ with $n$ punctures.
  The inclusion $\mathcal{O}_{\mathrm{WP}}^{\mathrm{qc}}(\riem^P) \rightarrow \mathcal{O}^{\mathrm{qc}}(\riem^P)$ is holomorphic. 
\end{theorem}
 Recalling Definition \ref{de:quasisymmetries_surfaces} of quasisymmetries of borders, we define Weil-Petersson class quasisymmetries as follows. 
 \begin{definition} Let $\riem$ and $\riem_1$ be bordered surfaces, with components $C$ and $C_1$ of the borders homeomorphic to $\mathbb{S}^1$. 
 We say that a quasisymmetry $\phi:C \rightarrow C_1$ is Weil-Petersson class if, for any collar charts $\psi$ and $\psi_1$ of $C$ and $C_1$ respectively, $\psi_1 \circ \phi \circ \psi^{-1} \in \mathrm{QS}_{\mathrm{WP}(\mathbb{S}^1)}$.  

 If $\riem$ and $\riem_1$ are bordered surfaces of type $(g,n)$, 
 We say that a map $\phi \in \mathrm{QS}(C,C_1)$ is Weil-Petersson if its restriction to each component of $\partial \riem$ is Weil-Petersson class in the above sense. Denote the set of such Weil-Petersson class maps by $\mathrm{QS}_{\mathrm{WP}}(\partial \riem,\partial \riem_1)$. 
\end{definition}

The Weil-Petersson class Teichm\"uller space of a type $(g,n)$ bordered surface was defined by Radnell-Schippers-Staubach as follows. 
\begin{definition}  \label{de:WP_Teich_space_us_first_version} Let $\riem$ be a type $(g,n)$ bordered surface.  The Weil-Petersson class Teichm\"uller space of $\riem$ is 
 \[  T_{\mathrm{WP}}(\riem) = \left\{ [\riem,F,\riem_1] \in T(\riem) : \left. F \right|_{\partial \riem} \in \mathrm{QS}_{\mathrm{WP}}(\partial \riem,\partial \riem_1)    \right\}.  \]
\end{definition}  

Recall the definitions of $\mathcal{C}$ and $\mathbf{F}_q$ given in (\ref{eq:fibre_C_definition}), (\ref{eq:fibre_F_definition}), and surrounding text. 
Let 
\[  \mathcal{C}_{\mathrm{WP}} = \left. \mathcal{C} \right|_{T_{\mathrm{WP}}(\riem)}.  \]  
For a point $q$ let $\mathbf{F}_{\mathrm{WP},q}$ be the restriction of $\mathbf{F}_q$ to $\mathcal{O}_{\mathrm{WP}}^{\mathrm{qc}}$.
 
This is well-defined because $\mathrm{ModI}(\riem)$, and in particular $\mathrm{DB}(\riem)$, obviously preserve $T_{\mathrm{WP}}(\riem)$.
Finally let $U_{\zeta,K}=T(V_{\zeta,K})$ be the open sets induced on $\left(\mathcal{O}^{qc}_{\mathrm{WP}}\right)^n$ induced by the open sets $V_{\zeta,K}$ in $\left(\mathcal{O}^{qc}_{\mathrm{WP}}(\riem_1)\right)$, and let 
\[ \mathcal{G}_{\mathrm{WP}}:U \times V_{\zeta,K} \rightarrow T_{\mathrm{WP}}(\riem) \]
be obtained by restricting the second component of $\mathcal{G}$ to the open subset $V_{\zeta,K}$ of $\mathcal{O}^{\mathrm{qc}}_{WP}(\riem_1)$. $\mathcal{G}$ of course depends on the choices of $\zeta$, and $K$, and the point $p \in \mathcal{C}^{-1}(q)$, so one obtains a collection of such maps.  Recall that the entire construction depends on a choice of quasisymmetries $\tau_1,\ldots,\tau_n$ parametrizing the boundaries. 

\begin{theorem}[\cite{RSS_Filbert2}]  \label{th:us_WP_atlas_fiber} Let $\riem$ be a type $(g,n)$ bordered surface with $2g-2+n>0$, and fix an $n$-tuple of quasisymmetries $\tau=(\tau_1,\ldots,\tau_n)$ where $\tau_k \in \mathrm{QS}_{\mathrm{WP}}(\mathbb{S}^1,\partial_k \riem)$. 
 The collection of maps $\mathcal{G}^{-1}$ induced by $\tau$ as above form an atlas, so that $T_{\mathrm{WP}}(\riem)$ is a complex Hilbert manifold with respect to this atlas.  This complex structure is independent of $\tau$.  
\end{theorem}
\begin{remark} It was shown that the domains of maps in this  atlas form the base for a Haudorff, second countable, and separable topology. 
\end{remark}
\begin{remark}
    By composing with $\chi^n_{\mathrm{WP}} \circ \mathcal{E}$, the local model of $T_{\mathrm{WP}}(\riem)$ is $\mathbb{C}^d \oplus \left( A_1^2(\disk_+) \oplus \mathbb{C} \right)^n$. 
\end{remark}
\begin{remark} With respect to this complex structure, the fibers are complex manifolds and $\mathcal{C}_{\mathrm{WP}}$ is a holomorphic map with local holomorphic submersions. This follows more or less directly from the definitions. 
\end{remark}
\begin{remark} There are two exceptional cases ruled out by $2g -2+n>0$, because of the degeneration of $T_{\mathrm{WP}}(\riem^P)$ to a point, as we saw in Remark \ref{re:fiber_exceptional_cases}. The exceptional case $g=0$ and $n=1$ is just $T_{\mathrm{WP}}(\disk_-)$, so that the fiber model identifies $T_{\mathrm{WP}}(\disk_-)$ with suitably normalized elements of $\mathrm{O}^{\mathrm{qc}}_{\mathrm{WP}}(\sphere \backslash \{0\})$. Thus this case is just the pre-Schwarzian model of  $T_{\mathrm{WP}}(\disk_-)$ already given by Cui and Takhtajan-Teo. 

The case $g=0$, $n=1$ can be dealt with in the same way as it was in Radnell-Schippers \cite{RadnellSchippers_annulus} by repeating those arguments for $T_{\mathrm{WP}}(\mathbb{A})$, though this has not been done.
\end{remark}

We also have 
\begin{theorem}[\cite{RSS_Filbert2}] \label{th:change_of_base_holo_fiber} 
 Let $\riem$ be a type $(g,n)$ bordered surface such that $2g-2+n>0$, and let $F:\riem \rightarrow \riem_0$ be a quasiconformal map such that $[\riem,F,\riem_0] \in T_{\mathrm{WP}}(\riem)$. The change of base point map \[ R_F:T_{\mathrm{WP}}(\riem_0) \rightarrow T_{\mathrm{WP}}(\riem) \]
 is a biholomorphism with respect to the complex structure obtained in \emph{Theorem \ref{th:us_WP_atlas_fiber}}.  
\end{theorem}
This of course will play a role in proving the compatibility of the different complex structures.

The following theorem is also essential.
\begin{theorem}[\cite{RSS_Filbert2}]
  Let $\riem$ be a type $(g,n)$ bordered surface for $2g-2+n>0$. The modular group $\mathrm{ModI}$ preserves $T_{\mathrm{WP}}(\riem)$. Furthermore, it acts properly discontinuously and fixed-point-freely by biholomorphisms.  
\end{theorem}

In \cite{RSS_WP_arxiv}, later appearing in \cite{RSS_WP1}, Radnell-Schippers-Staubach developed a local characterization of differentials in terms of their behaviour on collars. Let $\riem$ be a type $(g,n)$ Riemann surface. By definition of bordered surface, there is an atlas of charts on $\riem \cup \partial \riem$ mapping into the closed disk (equivalently, the upper half plane), such that the transition maps are homeomorphisms with respect to the relative topology on the closed disk and holomorphic on the interior (see e.g. Ahlfors Sario \cite{Ahlfors_Sario}).  We call a chart an interior chart if it maps into the open disk (equivalently, its domain contains only points in $\riem$); we call it a border chart if its domain contains a point in $\partial \riem$. 

Let $\alpha \in \mathscr{D}_{l,m}(\riem)$. Given a chart $\psi:U \rightarrow V$ in this atlas, we say $\psi^* \alpha = h(z) dz^l d\bar{z}^k \in L^2_{l,m}(\disk_+;V)$ if and only if 
\[   \int_V \frac{|h(z)|^2}{(1-|z|^2)^{2-2l-2m}} \frac{d\bar{z} \wedge dz}{2i} <\infty.   \]
In other words we demand that the pull-back under the chart be in $L^2$ with respect to the hyperbolic metric on $\disk_+$.  

\begin{theorem}[\cite{RSS_WP1}] 
 Let $\riem$ be a bordered surface of type $(g,n)$. The following are equivalent.
 \begin{enumerate}[font=\upshape]
  \item $\alpha \in L^2_{l,m}(\riem)$;
  \item $\psi^*\alpha \in L^2_{l,m}(\disk_+;V)$ for all charts $\psi:U \rightarrow V$ in the border atlas;
  \item $\psi^*\alpha \in L^2_{l,m}(\disk_+;V)$ for a specific collection of interior charts which cover $\riem$ and for specific collar charts $\psi_1,\ldots,\psi_n$ of $\partial_1 \riem,\ldots,\partial_n \riem$. 
 \end{enumerate}
\end{theorem} 
Essentially, the idea is that the hyperbolic metric on a collar neighbourhood of a boundary behaves analytically the same as the hyperbolic metric on the disk near $\mathbb{S}^1$. This is a useful tool for applying the sewing point of view.

This local characterization can be used to show the following. 
\begin{theorem}[\cite{RSS_WP1}] \label{th:us_WP_homotopy_characterize}
 Let $\riem$ and $\riem_1$ be bordered surfaces of type $(g,n)$ and let $f:\riem \rightarrow \riem_1$ be quasiconformal.  Then $\left. F \right|_{\partial \riem} \in \mathrm{QS}_{\mathrm{WP}}(\riem)$ if and only if $F$ is homotopic rel boundary to a quasiconformal map $F_0:\riem \rightarrow \riem_1$ such that 
 \[ \mu(F_0) \in L^2_{-1,1}(\riem) \cap L^\infty_{-1,1}(\riem)_1. \]
\end{theorem}
Ahead, we will see that a result of Yanagishita shows that such a representative is obtained from the Douady-Earle extension.

\begin{remark} \label{re:Cuis_characterization_extended_us} It was also shown in \cite{RSS_WP_arxiv,RSS_WP2} that 
 $\phi_k = \left. F \right|_{\partial_k \riem} \in \mathrm{QS}(\mathbb{S}^1,\riem)$ for $k=1,\ldots,n$ if and only if the corresponding rigging $f = (f_1,\ldots,f_n)$, $f:\disk_+\rightarrow \riem^P$ obtained from the fiber map $\mathcal{C}$ is in $\mathcal{O}_{\mathrm{WP}}^{\mathrm{qc}}$. As a consequence, Theorem \ref{th:us_WP_homotopy_characterize} is a generalization of Theorem \ref{th:WP_quasisymmetry_characterize} of Cui and Takhtajan-Teo to type $(g,n)$ surfaces.  
\end{remark}

In particular, 
\begin{corollary} \cite{RSS_WP_arxiv,RSS_WP1} \label{co:WP_Teich_recharacterize_us}
 Let $\riem$ be a bordered Riemann surface of type $(g,n)$. Then
 \begin{equation} \label{eq:alternate_definition_WP_Teich}
  T_{\mathrm{WP}}(\riem) = \{ [\riem,F,\riem_1] \in T(\riem) \,:\, \mu(F) \in L^2_{-1,1}(\riem) \cap L^\infty_{-1,1}(\riem) \}. 
 \end{equation}
\end{corollary}
Thus (\ref{eq:alternate_definition_WP_Teich}) gives an alternate definition of the Weil-Petersson Teichm\"uller space for these surfaces.

\end{subsubsection}
\begin{subsubsection}{{\underline{Bers embedding model}}}
 Yanagishita \cite{Yanagishita} constructed a complex structure on Weil-Petersson Teichm\"uller space using the Bers embedding model. 
 We outline Yanagishita's results on the Bers embedding in the $L^2$ case. In fact, these results hold in $L^p$ for $p \geq 1$ (or in a few cases $p \geq 2$). 
 
 Let $G$ be an arbitrary Fuchsian group. 
 Let $T^p(G)$ denote the subset of elements $[\mu]$ of $T(G)$ with representatives $\mu \in L^p_{-1,1}(N)$ where $N$ is a fundamental domain of $G$.  
 \begin{remark}
 Equivalently, $\mu \in L^p_{-1,1}(\riem)$ on the resulting quotient surface $\riem = \disk_- /G$. Thus for $p=2$ this agrees with the Definition of equation (\ref{eq:alternate_definition_WP_Teich}) of Radnell-Schippers, (though that definition was only stated for type $(g,n)$ bordered surfaces).  From here on, we will restrict to $p=2$, which is the case with a complex Hilbert manifold structure and Riemannian metric. 
 \end{remark}
 
 Yanagishita showed that the map 
 \[  \sigma: T^2(G) \rightarrow L^\infty_{-1,1}(\disk_-;G)  \]
 taking $[\mu]$ to the Douady-Earle extension of the boundary values of $w_\mu$ is in $L^2_{-1,1}(\disk_-;G)$.  

 Denote $\mathrm{Ael}^2(\disk_-;G) = L^\infty_{-1,1}(\disk_-) \cap L^2_{-1,1}(\disk_-;G)$, where this is given the direct sum norm 
 \[ \| \cdot \|_{2,\infty} = \| \cdot \|_\infty + \| \cdot \|_2
 \]
 as in Cui (note that Cui used the notation $M_{1,*}$).  We also define the perturbed group $G_\psi$ for $\psi \in T(\disk_-;G)$ as in (\ref{eq:change_base_point_def_group}).
 
 Recall that Cui showed that the Weil-Petersson quasisymmetries form a group. This is not true for $G$-invariant quasisymmetries, but the following theorem generalizes the result of Cui in some sense. 
 \begin{theorem} \cite{Yanagishita}  \label{th:Yana_Douady_Earle_group} For any Fuchsian group $G$,   the following are equivalent.
 \begin{enumerate}[font=\upshape]
 \item
  $\psi \in T^2(G)$.  
 \item $\psi^{-1} \in T^2(G_\psi)$.
 \item $\sigma(\psi) \in L^\infty_{-1,1}(\disk_-;G) \cap L^2_{-1,1}(\disk_-;G)$.  
 \end{enumerate}
 \end{theorem} 
 Written in the direct picture, setting $\riem$ and $\riem_\psi$ to be $\disk_-/G$ and $\disk_-/G_\psi$ respectively, this shows that change of basepoint $R_\psi$ takes $T^2(\riem)$ to $T^2(\riem_\psi)$. 
 \begin{remark}
 Note that this is a different theorem than Theorem \ref{th:change_of_base_holo_fiber} even in the type $(g,n)$ case, since the complex structures are not yet known to be the same. Ahead we will see that they are in fact the same. 
 \end{remark}

  Yanagishita also defined the following analogue of the Teichm\"uller distance. For $\psi_1,\psi_2 \in T^2(G)$, define
  \[ \ell_{2,\infty}(G)(\psi_1,\psi_2) = \inf \left\{ \left\| \frac{\mu_1 - \mu_2}{1- \mu_1 \overline{\mu_2}} \right\|_{2,\infty} \,: \, \mu_k \in \mathrm{Ael}^2(G), [\mu_k]=\psi_k, k=1,2 \right\}.
  \]
 We also have the following theorem, partly generalizing Theorem \ref{th:WP_quasisymmetry_characterize} of Cui and Takhtajan-Teo. 
 \begin{theorem}[\cite{Yanagishita}]
  For any Fuchsian group, 
  \[ \beta(T^2(G)) \subset \beta(T(G)) \cap A^2_2(\disk_+;G). \]
 \end{theorem} 
 \begin{remark} This also generalizes one inclusion of Theorem \ref{th:WP_quasisymmetry_characterize}, but interestingly, it does so in a different way than Remark \ref{re:Cuis_characterization_extended_us}.
 \end{remark}

 Thus one can ask about topological properties of the Bers embedding. Indeed we have the following. 
 \begin{theorem}[\cite{Yanagishita}]
  Let $G$ be an arbitrary Fuchsian group.
  \begin{enumerate}
   \item $\sigma$ is continuous from $\left(T^2(G),\ell_{2,\infty}\right)$ to $\left(\mathrm{Ael}^2(G), \| \cdot \|_{2,\infty}\right)$. 
   \item The Bers embedding $\beta$ is continuous from $(T^2(G),\ell_{2,\infty})$ to $A^2_2(\disk_+;G)$. 
  \end{enumerate}
 \end{theorem}

 Much more can be said, if one furthermore assumes that $G$ satisfies what Yanagishita calls ``Lehner's condition''. 
 \begin{definition} A Fuchsian group $G$ is said to satisfy Lehner's condition if the infimum of the lengths of all simple closed geodesics in $\disk_-/G$ is finite. 
 \end{definition}
 He then shows the following, based on results of R. Rajeswara, D. Niebur and M. Sheingorn, and J. Lehner.
 \begin{theorem}[\cite{Yanagishita}] \label{th:Lehner_and_inclusion}
  A Fuchsian group $G$ satisfies Lehner's condition if and only if the inclusion from $A^2_2(G)$ into $A^\infty_2(G)$ is bounded. 
 \end{theorem}
 Observe that a type $(g,n)$ bordered surface can be seen directly from the definition to satisfy Lehner's condition.

 Yanagishita showed that under this condition, the Bers embedding is a homeomorphism onto its image.
 \begin{theorem}[\cite{Yanagishita}]  \label{th:Yanagishita_Bers_embedding}
  Let $G$ be a Fuchsian group satisfying Lehner's condition. Then the Bers embedding $\beta$ is a homeomorphism from $(T^2(G),\ell_{2,\infty})$ onto its image in $A^2_2(\disk_+;G)$. Thus $T^2(G)$ possesses a complex Hilbert manifold structure. 
 \end{theorem} 
 This completes the construction of the complex structure from the Bers embedding, analogous to the classical $L^\infty$ Banach manifold structure of Section \ref{subsubse:Bers_embedding_complex}. 

 Furthermore, with respect to this complex structure, change of base point is a biholomorphism. 
 \begin{theorem}[\cite{Yanagishita}] 
  For any Fuchsian group $G$ and $\psi \in T^2(G)$, 
  the change of base point map $R_\psi$ is a homeomorphism from $(T^2(G),\ell_{2,\infty}(G))$ onto $(T^2(G_\psi),\ell_{2,\infty}(G_\psi))$.  If $G$ satisfies Lehner's condition, then $R_\psi$ is a biholomorphism with respect to the complex Hilbert manifold structure induced from $A^2_2(\disk_+;G)$ by Theorem \ref{th:Yanagishita_Bers_embedding}. 
 \end{theorem}
 On the Riemann surface, the change of base point takes $T_{\mathrm{WP}}(\disk_-/G)$ to $T_{\mathrm{WP}}(\disk_-/G_\psi)$. 
 \begin{remark}
  It should be noted that the statement regarding complex structure in Theorem \ref{th:Yanagishita_Bers_embedding} does not follow from Theorem \ref{th:us_WP_atlas_fiber} even in the case of type $(g,n)$ surfaces, nor does Theorem \ref{th:us_WP_atlas_fiber} follow from Theorem \ref{th:Yanagishita_Bers_embedding}. This is because the complex structures had not yet been shown to be the same.  Ultimately, it was shown that the three complex structures agree, as we will see below. 
 \end{remark}
\end{subsubsection}
\begin{subsubsection}{{\underline{Deformation model}}}

 The deformation model of the complex structure and tangent space was given in \cite{RSS_WP_arxiv}, and later published in \cite{RSS_WP2}, for type $(g,n)$ bordered surfaces such that $2g-2+n>0$. The theory is based the Ahlfors-Weill reflection, as in the classical $L^\infty$ model of Section \ref{subsubse:deformation_complex_structure}.  
 
 Below, we use the notation 
 \[  \mathrm{TBD}(\riem) = L^\infty_{-1,1}(\riem) \cap L^2_{-1,1}(\riem).   \]
 Here TBD stands for ``tangent Beltrami differential'', that is tangent vectors to curves in both $L^\infty_{-1,1}(\riem)_1$ and $L^2_{-1,1}(\riem)$. The following existence theorem, which was of great use in describing the tangent space, justifies considering this space. 
 \begin{theorem}[\cite{RSS_WP2}]
  Let $\riem$ be a type $(g,n)$ bordered surface such that $2g-2+n>0$. Let $t \mapsto [\riem,F_t,\riem_t]$ be a holomorphic curve in $T_{\mathrm{WP}}(\riem)$ with respect to the complex structure obtained in Theorem \ref{th:us_WP_atlas_fiber}, such that $[\riem,F_0,\riem_0]=[\riem,\mathrm{Id},\riem]$. There is a $\delta>0$ and a one-parameter family of representatives $(\riem,F_t\riem_t)$ for $|t|<\delta$ such that
  \begin{enumerate}[font=\upshape]
   \item $\|\mu(F_t)\|_2$ is uniformly bounded in $t$;
   \item $t\mapsto \mu(F_t)$ is holomorphic in $L^2_{-1,1}(\riem)$; 
   \item $t \mapsto \mu(f_t)$ is holomorphic in $L^\infty_{-1,1}(\riem)$. 
  \end{enumerate}
 \end{theorem} 
  Thus one can find curves which are simultaneously holomorphic in both the $L^\infty$ and $L^2$ spaces. 

  It was also shown that
  \begin{theorem}[\cite{RSS_WP2}] Let $\riem$ be a type $(g,n)$ bordered surface. We have $H_{-1,1}(\riem) \subset \Omega_{1,1}(\riem)$ and the inclusion is bounded.
  \end{theorem}
  This is equivalent to the boundedness of the inclusion $A_2^2(\riem^*) \rightarrow A_2^\infty(\riem^*)$. This also follows from Theorem \ref{th:Lehner_and_inclusion}, using the fact that type $(g,n)$ bordered surfaces satisfy Lehner's condition. Alternately, one can say that these two theorems imply that type $(g,n)$ surfaces satisfy Lehner's condition.
  
  As a corollary, we obtain
  \begin{corollary}[\cite{RSS_WP2}] \label{co:harmonic_inclusion_bounded}
   If $\riem$ is a type $(g,n)$ bordered surface, then 
   \[  \Omega_{-1,1}(\riem) \cap \mathrm{TBD}(\riem) = H_{-1,1}(\riem).   \]
  \end{corollary}
  
  It is natural then to define 
  \[  \mathcal{N}_r(\riem) = \mathcal{N}(\riem) \cap L^2_{-1,1}(\riem) \]
  where recall that $\mathcal{N}(\riem)$ are the infinitesimally trivial differentials, that is, those differentials which are tangent to the Teichm\"uller equivalence relation.
  Using Corollary \ref{co:harmonic_inclusion_bounded}, and properties of the Bergman projection on differentials, Radnell-Schippers-Staubach obtained the following characterization of $\mathrm{TBD}(\riem)$ modulo infinitesimally trivial differentials.  
  \begin{theorem}[\cite{RSS_WP2}] \label{th:modding_out_inf_trivial_us}
   Let $\riem$ be a type $(g,n)$ bordered surface. Then we have the direct sum decomposition 
   \[  \mathrm{TBD}(\riem) = \mathcal{N}_r(\riem) \oplus H_{-1,1}(\riem).   \]
   The projection from $\mathrm{TBD}(\riem)$ onto $H_{-1,1}(\riem)$ is bounded. Finally, 
   $\mathcal{N}_r(\riem)$ is the kernel of the derivative of the Bers embedding $D_0(\beta \circ \Phi)$ at $[\riem,\mathrm{Id},\riem]$.  
  \end{theorem}

  This in turn led to a description of the tangent vectors in $T_{\mathrm{WP}}(\riem)$. 
  \begin{theorem}[\cite{RSS_WP2}]
   Let $v$ be a vector tangent to $T_{\mathrm{WP}}(\riem)$ at $[\riem,\mathrm{Id},\riem]$. There is a holomorphic curve $t \mapsto [\riem,F_t,\riem_t]$ for $|t|<\delta$ such that $[\riem,F_0,\riem_0]=[\riem,\mathrm{Id},\riem]$ and 
   \begin{enumerate}[font=\upshape]
    \item $\mu(F_t) \in H_{-1,1}(\riem)$ for $|t|<\delta$;
    \item $\mu_t$ is holomorphic in $t$ with respect to $H_{-1,1}(\riem)$;
    \item $v$ is tangent to $[\riem,F_t,\riem_t]$ at $t=0$. 
   \end{enumerate}
  \end{theorem}

  This allows us to construct a new atlas. Recall that $\Phi:L^\infty_{-1,1}(\riem)_1 \rightarrow T(\riem)$ takes a Beltrami differential to its Teichm\"uller space representative by solving the Beltrami equation. First, we locally invert $\Phi$. 
  \begin{theorem}[\cite{RSS_WP2}] \label{th:us_WP_atlas_deformation}
   Let $\riem$ be a type $(g,n)$ bordered Riemann surface such that $2g-2+n>0$. There exists a neighbourhood $U$ of $0 \in H_{-1,1}(\riem)$ such that $\left. \Phi \right|_{U}$ is a biholomorphism onto $T_{WP}(\riem)$ with respect to the complex structure obtained from \emph{Theorem \ref{th:us_WP_atlas_fiber}}. 
  \end{theorem}
  Using holomorphicity of change of base point with respect to the complex structure (Theorem \ref{th:change_of_base_holo_fiber}), we have 
  \begin{corollary}[\cite{RSS_WP_arxiv,RSS_WP2}]
   Let $\riem$ be a bordered surface of type $(g,n)$ such that $2g-2+n>0$. Let $F:\riem \rightarrow \riem_0$ be a quasiconformal map. There is an open neighbourhood $B$ of $0$ in $H_{-1,1}(\riem_0)$ such that 
   \[   \Psi_{(\riem,f,\riem_0)} := R_{F_\mu} \Phi_{\riem_0}   \]
   is a biholomorphism onto an open neighbourhood of $[\riem,F,\riem_0]$. The collection of inverses of these maps form an atlas for a complex Hilbert manifold structure on $T_{\mathrm{WP}}(\riem)$ compatible with the complex structure obtained in Theorem \ref{th:us_WP_atlas_fiber}.\footnote{There is a typographical error in the statement of this result in \cite[Theorem 4.8]{RSS_WP2}: $H_{-1,1}(\riem)$ should be replaced by $H_{-1,1}(\riem_0)$ as it is here.} 
  \end{corollary}
  \begin{remark} 
   The inverses of these maps are $\Lambda^{-1}_{\riem_0} \beta_{\riem_0} R_{F_{0}^{-1}}$, so up to the obviously isometric reflection into $A^2_2(\riem_0^*)$ these are the $L^2$ versions of the classical charts in Section \ref{subsubse:deformation_complex_structure}.  Note that because the Ahlfors-Weill reflection $\Lambda$ is obviously an isometry, we also automatically obtain that $\beta_{\riem_0} R_{F_0^{-1}}$ is an atlas of charts into $A_2^2(\riem_0^*)$. 
  \end{remark}

  The results above show that we can identify the tangent space at an arbitrary point $[\riem,F,\riem_0] \in T_{\mathrm{WP}}(\riem_0)$ with   
  \[   T_{[\riem,F,\riem_0]} T_{\mathrm{WP}}(\riem) \cong  D_0 R_{F} H_{-1,1}(\riem_0).     \]
  It is an immediate consequence that the Weil-Petersson metric converges on every tangent space. In fact we can define the Weil-Petersson pairing as follows. Let $\mu,\nu \in H_{-1,1}(\riem)$ be elements in the tangent space of $T_{\mathrm{WP}}(\riem)$ at $[\riem,\mathrm{Id},\riem]$. Then let 
  \[ \left< \mu, \nu \right>_{[\riem,\mathrm{Id},\riem]}   \]
  denote the inner product in $H_{-1,1}(\riem)$.  Then for any arbitrary point $[\riem,F,\riem_0] \in T_{\mathrm{WP}}(\riem)$ and tangent vectors $\mathbf{v},\mathbf{w} \in  T_{[\riem,F,\riem_0]} T_{\mathrm{WP}}(\riem)$
  we define the Weil-Petersson Riemannian metric by 
  \begin{equation} \label{eq:WP_pairing_general_surfaces}  \left< D (R_{F})^{-1}\mathbf{v}, D (R_{F})^{-1} \mathbf{w} \right>_{[\riem,F,\riem_0]}.  
  \end{equation}
  This converges at every tangent space and is immediately seen to be smooth with respect to the by biholomorphicity of the change of base point map. 
  It is the unique metric which is invariant under change of base point and agrees with the pairing in $H_{-1,1}(\riem_0)$ at the identity.

  Later, Yanagishita extended the deformation atlas result to general surfaces $\riem=\disk_-/G$ for Fuchsian groups $G$ satisfying Lehner's condition \cite{Yanagishita_Kahler}.
  \begin{theorem}[\cite{Yanagishita_Kahler}] \label{th:Yanagishita_WP_atlas_deformation} Let $G$ be a Fuchsian group  satisfying Lehner's condition. The charts $\beta_{\riem_0} R_{F_0^{-1}}$ forms an atlas for a complex Hilbert manifold structure compatible with the complex structure obtained from Bers embedding via \emph{Theorem \ref{th:Yanagishita_Bers_embedding}}.
  Change of base point is holomorphic. 
  \end{theorem}
  He also extended the description of the tangent space and the Weil-Petersson pairing to surfaces $\riem = \disk_-/G$ for Fuchsian groups $G$ satisfying Lehner's condition. 
  \begin{theorem}[\cite{Yanagishita_Kahler}]
      \emph{Corollary \ref{co:harmonic_inclusion_bounded}} holds more generally for surfaces satisfying Lehner's condition.
  \end{theorem}
  \begin{theorem}[\cite{Yanagishita_Kahler}]
     \emph{Theorem \ref{th:modding_out_inf_trivial_us}} holds more generally for surfaces satisfying Lehner's condition. 
  \end{theorem}
  He thus obtains the convergent Weil-Petersson pairing (\ref{eq:WP_pairing_general_surfaces}) for surfaces satisfying Lehner's condition.

  The statements in Theorem \ref{th:us_WP_atlas_deformation}  and in Theorem \ref{th:Yanagishita_WP_atlas_deformation} - that the atlas defines a complex structure - are the same in the case of type $(g,n)$ bordered surfaces; indeed Theorem \ref{th:us_WP_atlas_deformation} is a consequence of Theorem \ref{th:Yanagishita_WP_atlas_deformation}. However, the statements about the compatibility of this structure are completely different: Theorem \ref{th:us_WP_atlas_deformation} says that this atlas is compatible with the complex Hilbert manifold structure obtained from the fiber structure, whereas Theorem \ref{th:Yanagishita_WP_atlas_deformation} says that this atlas is compatible with the complex Hilbert manifold structure obtained from the Bers embedding. Putting these results together, we obtain that all three agree. That is, by Theorems \ref{th:us_WP_atlas_fiber},  \ref{th:us_WP_atlas_deformation}, and \ref{th:Yanagishita_WP_atlas_deformation} we have the following.
  \begin{corollary}
   Let $\riem$ be a type $(g,n)$ bordered surface such that $2g-2+n>0$. The complex Hilbert manifold structures obtained from the fiber model, the Bers embedding, and the deformation model, are all equivalent. 
  \end{corollary}
  This was observed by the authors in our survey \cite{SchippersStaubach_comparison}. 
  \begin{remark}
  To the $L^2$ analogues of the two classical complex structures, we have added the new complex structure derived from the fiber/CFT point of view. In the classical $L^\infty$ case this also required some extra work \cite{RadnellSchippers_fiber}. 
  \end{remark}

 Finally, we conclude with another interesting result of Yanagishita. In stark contrast to Cui's completeness theorem (Theorem \ref{th:Cui_completeness}), the Weil-Petersson distance is almost never complete.  
 \begin{theorem}[\cite{Yanagishita_completeness}]
  Let $\riem$ be a Riemann surface satisfying Lehner's condition. If $\riem$ is conformally equivalent to the punctured disk or the disk, then $T_{\mathrm{WP}}(\riem)$ is complete. Otherwise, it is not complete. 
 \end{theorem}
 \begin{remark}
     In fact the theorem is shown to hold for $p\geq 2$; the completeness of the $p$-Weil-Petersson Teichm\"uller space of the disk for $p\geq 2$ is due to K. Matsuzaki \cite{Matsuzaki_rigidity}. 
 \end{remark}
  \end{subsubsection}
 \end{subsection}
\end{section}
\begin{section}{K\"ahler geometry and global analysis of Weil-Petersson Teichm\"uller space} \label{se:Kahler_and_global}
 In this section we give a brief overview of
 some geometric results associated with the Weil-Petersson metric. The emphasis is on Kählericity and Kähler potentials, Chern classes and index theorems. Moreover, we shall also recall some of the results in the finite dimensional case (the case for compact Riemann surfaces), both for the sake of comparison with the currently available infinite dimensional results, and also for their intrinsic interest in the Weil-Petersson context. The presentation given here is by no means exhaustive and we confine ourselves to just a selection of results that are relevant to the WP-class Teichm\"uller theory. 
 
\begin{subsection}{K\"ahlericity, curvatures and K\"ahler potentials}

In the paper \cite{Weil2} Weil asked whether the WP metric is K\"ahler and if so, what is its curvature.\\
This problem was solved by L. Ahlfors \cite{Ahlfors_Kahler} , who showed that for compact Riemann
 surfaces the Weil-Petersson metric is K\"ahler. 
 More specifically he proved the following theorem:
 \begin{theorem}
Let $\mathscr{R}$ be a compact genus $g>1$ Riemann surface. On the Teichmüller space \(T(\mathscr{R})\) Bers' coordinates are geodesic
for the Weil-Petersson metric at the reference point. In particular this metric
is Kähler.     
 \end{theorem}
 
 In the
 commentary to his collected works \cite{Ahlfors collected}, Ahlfors mentioned that Weil also had a proof but had not published it.  Later Ahlfors computed the curvatures of holomorphic sections and the
 Ricci curvature \cite{Ahlfors_curvature}, which resulted in

 \begin{theorem}
The Ricci curvature, the holomorphic sectional curvature and the
scalar curvature of the Weil-Petersson metric are negative.
 \end{theorem}

 In \cite{Wolpert}, S. Wolpert used the classical Maa{\ss} calculus of invariant differential operators to obtain a formula for the Riemannian curvature tensor, showing that the sectional curvatures of the Weil-Petersson metric are negative. The negative curvature result was obtained independently by H. L. Royden \cite{Royden1}, \cite{Royden2} and A. Tromba \cite{Tromba}.
 More explicitly let \(\mathscr{R}\) be a compact Riemann
surface of genus \(g>1\). Denote by \(\left\{\mu_{\alpha} ; \alpha=1, \ldots, 3 g-3\right\}\) a basis of the vector space of harmonic Beltrami differentials on \(\mathscr{R}\) corresponding to a set of tangent
vectors \(\left\{\partial / \partial s_{\alpha} ; \alpha=1, \ldots, 3 g-3\right\}\). Introducing the basis
\begin{equation}
    g_{\alpha \bar{\beta}}=\int_{\mathscr{R}} \mu_{\alpha} \overline{\mu_{\beta}}\, dA_{\mathrm{hyp}},
\end{equation} 
one can write the Weil-Petersson metric as
$$
\omega_{\mathrm{WP}}=i g_{\alpha, \bar{\beta}}\, d z^{\alpha} \wedge d \bar{z}^{{\beta}},
$$
using Einstein's summation convention. The following result was established in \cite{Wolpert} and \cite{Tromba}.

\begin{theorem}\label{thm_wolpert-tromba}

Let $\mathscr{R}$ be a compact Riemann surface of genus $g>1$ then one has

\begin{enumerate}[font=\upshape]

\item The holomorphic sectional and Ricci curvatures of $\omega_{\mathrm{WP}}$ on $T(\mathscr{R})$ for $g>1$
are bounded from above by $\frac{-1}{2 \pi(g-1)}$.

\item The sectional curvature of $\omega_{\mathrm{WP}}$ is negative.

\item The curvature tensor of the Petersson-Weil metric is equal to
\begin{equation}
\begin{aligned} R_{\alpha \bar{\beta} \gamma \bar{\delta}}= & -2 \int_{\mathscr{R}}(D-2)^{-1}\left(\mu_{\alpha} \overline{\mu_{\beta}}\right)\left(\mu_{\gamma} \overline{\mu_{\delta}}\right) \, dA_z \\ & -2 \int_{\mathscr{R}}(D-2)^{-1}\left(\mu_{\alpha} \overline{\mu_{\delta}}\right)\left(\mu_{\gamma} \overline{\mu_{\beta}}\right) \, dA_{\mathrm{hyp}}\end{aligned}
\end{equation}
where \(D\) is the real Laplacian on \(L^{2}\)-functions on \(\mathscr{R}\) with spectrum in $(-\infty, 0]$.
\end{enumerate}
\end{theorem}
In \cite{Schumacher}, Schumacher went beyond the the negativity of the sectional curvature in Theorem \ref{thm_wolpert-tromba}, and proved the following result.
\begin{theorem}
The Weil-Petersson metric on \(T(\mathscr{R})\) for \(g>1\) has
strongly negative curvature in the sense of \emph{Y. T. Siu}, i.e. its Riemann tensor satisfies
$$
R_{\alpha \bar{\beta} \gamma \bar{\delta}}\left(A^{\alpha} \overline{B^{\beta}}-C^{\alpha} \overline{D^{\beta}}\right)\left(\overline{A^{\delta}} B^{\gamma}-\overline{C^{\delta}} D^{\gamma}\right) \geq 0
$$
for all complex vectors \(A^{\alpha}, B^{\beta}, C^{\gamma}, D^{\delta}\), and equality holds only for \(A^{\alpha} \overline{B^{\beta}}=\)
\(C^{\alpha} \overline{D^{\beta}}\) for all \(\alpha\) and \(\beta\).
\end{theorem}

{Recall from Section \ref{se:polarizations_Siegel} that in finite 
dimensions the Siegel disk possesses a natural K\"ahler metric. The K\"ahler potential for this metric is   $-\mathrm{Tr} \log{(1-\overline{Z} Z)}$. Recall that Segal defined an infinite Siegel disk on which this K\"ahler metric is defined \cite{SegalUnitary}.

Note that the condition that $Z$ be Hilbert-Schmidt (part of Segal's definition of the infinite Siegel disk) is precisely the condition that the K\"ahler potential converges.  
 Kirillov and Yuriev \cite{KY2}
  and Nag \cite{Nag_bulletin} showed that the pull-back
 of a natural K\"ahler metric on the infinite Siegel disk is the Weil-Petersson metric. 
 Thus 
 \begin{equation} \label{eq:period_kahler_potential_WP}
 -\operatorname{Tr} \log \left(1-Z^{\ast} Z\right)
 \end{equation}
 is also a K\"ahler potential
 for the Weil-Petersson metric, as was observed explicitly in \cite{KY2,HongRajeev}.  Hong-Rajeev
 noted that $\text{Diff}(\mathbb{S}^1)/\text{M\"ob}(\mathbb{S}^1)$ was not complete with
 respect to the K\"ahler metric, which indicated that it was not the correct analytic setting
 for the Weil-Petersson metric. As we saw in Secton \ref{se:period}, this analytic setting is provided by $T_{\mathrm{WP}}(\disk_-)$.} 
 \\

{A generalization of the results of Ahlfors, Wolpert and Tromba mentioned above to the infinite dimensional setting was Takhtajan-Teo's K\"ahler-Einstein theorem for the universal Teichm\"uller space. To describe it consider the universal Teichm\"uller space $T(\disk_-)$, the universal Teichm\"uller curve $\mathscr{T}(\disk_-)$, and a projection $\pi:\mathscr{T}(\disk_-)\to T(\disk_-).$ The universal Teichm\"uller curve $\mathscr{T}(\disk_-)$ is a holomorphic fiber space over $T(\disk_-)$ considered by Bers.  The fiber over each point $[\mu]\in T(\disk_-)$ is a quasidisk $w_{\mu} (\mathbb{D}_{-})\subset \sphere$ with the complex structure inherited from the Riemann sphere and $\mathscr{T}(\disk_-)=\{([\mu], z);\,[\mu]\in T(\disk_-),\, z\in w_{\mu} (\mathbb{D}_{-}) \}.$ The Weil-Petersson universal Teichm\"uller space ${T}_{\mathrm{WP}}(\disk_-)$ is the connected component of the origin (i.e. the identity) of \(T(\disk_-)\), and the Weil-Petersson universal Teichm\"uller curve $\mathscr{T}_{\mathrm{WP}}(\disk_-)$ is obtained by restricting to $T_{\mathrm{WP}}(\disk_-)$. 
 More explicitly, one can summarize the results of Takhtajan-Teo as follows.
Let \(G=\frac{1}{2}\left(\Delta+\frac{1}{2}\right)^{-1}\), where $\Delta$ is the Laplace-Beltrami operator of the hyperbolic metric on \(\mathbb{D}_{-}\), and let $G(z, w)$ denote the integral kernel of $G$.
This kernel is explicitly given by
\(G(z, w)=\frac{2 u+1}{2 \pi} \log \frac{u+1}{u}-\frac{1}{\pi},\) where \(u=u(z, w)=\frac{|z-w|^{2}}{\left(1-|z|^{2}\right)\left(1-|w|^{2}\right)},\) see e.g. \cite{Hejhal}.\\

Then Takhtajan-Teo \cite{Takhtajan_Teo_Memoirs} proved the following ground-breaking result in Weil-Petersson Teichm\"uller theory of surfaces of infinite conformal type.
\begin{theorem}
Set
\begin{equation}\label{resolvent}
    G(f)(z)=\iint_{\mathbb{D}_{-}} G(z, w) f(w)\, dA_{\mathrm{hyp}}.
\end{equation}

then the following claims hold:
\begin{enumerate}[font=\upshape]
    \item The Weil-Petersson metric is a Kähler metric on the Hilbert manifold $T(\disk_-)$. Moreover the Bers coordinates are geodesic coordinates at the origin of \(T(\disk_-)\).
\item Let \(\mu_{\alpha}, \mu_{\beta}, \mu_{\gamma}, \mu_{\delta} \in H_{-1,1}(\disk_-)\) $($defined in \eqref{TTs H-class} which is isomorphic to the tangent space at the origin of the universal Teichm\"uller space$)$, be orthonormal tangent vectors.
Then the Riemann tensor at the origin of $T(\disk_-)$ is given by
$$
R_{\alpha \bar{\beta} \gamma \bar{\delta}}=-\frac{\partial^{2} g_{\alpha \bar{\beta}}}{\partial t_{\gamma} \partial \bar{t}_{\delta}}=-\left\langle G\left(\mu_{\alpha} \bar{\mu}_{\delta}\right), \mu_{\beta} \bar{\mu}_{\gamma}\right\rangle-\left\langle\mu_{\alpha} \bar{\mu}_{\beta}, G\left(\bar{\mu}_{\gamma} \mu_{\delta}\right)\right\rangle
$$
where $g_{\alpha \bar{\beta}}$ is the Weil-Petersson metric, the pairing $\langle \cdot, \cdot\rangle$ is the pairing in $L^2_{-1,1}(\disk_-)$.

\item $\mathscr{T}_{\mathrm{WP}}(\disk_-)$ and \(T_{\mathrm{WP}}(\disk_-)\) are topological groups in the Hilbert manifold topology of \(T(\disk_-)\), and the Hilbert manifold \(T_{\mathrm{WP}}(\disk_-)\) is also a Kähler-Einstein manifold with the negative definite Ricci
tensor
$$
\mathrm{Ric}_{\mathrm{W P}}=-\frac{13}{12 \pi} \omega_{\mathrm{W P}},
$$
where $\omega_{\mathrm{WP}}$ is the symplectic form of the Weil-Petersson Kähler metric on $T(\disk_-)$.
\end{enumerate}
 
\end{theorem}}

Takhtajan and Teo \cite{Takhtajan_Teo_Memoirs} also introduced the so-called \emph{universal Liouville action} defined by
\begin{equation}\label{Liouville}
\mathrm{S}_{1}([\mu])=\iint_{\mathbb{D}_+}\left|\mathcal{A}\left({f}^{\mu}\right)\right|^{2} d \mathrm{A}_z+\iint_{\mathbb{D}_{-}}\left|\mathcal{A}\left({g}_{\mu}\right)\right|^{2} d\mathrm{A}_z -4 \pi \log \left|{g}_{\mu}^{\prime}(\infty)\right|,
\end{equation}
where \(w_{\mu}={g}_{\mu}^{-1} \circ {f}^{\mu}\) is the conformal welding corresponding to \([\mu]\) in the WP-class Teichm\"uller space, and $\mathcal{A} f:=f''/f'$ is the pre-Schwarzian derivative of $f$.
They showed that this action is a Kähler potential for the WP-metric
on the Weil-Petersson Teichmüller space, in fact
\begin{equation}
    \frac{i}{2}\partial \overline{\partial} \mathrm{S}_{1}=  \omega_{\mathrm{W P}}.
\end{equation}
Furthermore, if $f_\mu$ is the conformal map associated to $[\mu] \in T_{\mathrm{WP}}(\disk_-)$, and $Z$ is the Grunsky map associated to $f$, then 
\begin{equation} \label{eq:Kahler_potential_Grunsky}
 S_2 = \log \det (I - \overline{Z} Z)  
\end{equation}
is a potential for the Weil-Petersson metric. This result of Takhtajan and Teo generalizes the formula (\ref{eq:period_kahler_potential_WP}) for the K\"ahler potential obtained from pulling back that on the finite Siegel disk, and completes the picture in Nag \cite{Nag_bulletin}, Nag-Verjovsky \cite{NagVerjovsky}, and Kirillov-Yuriev \cite{KY2}.\\

{Finally in relation to the potential above, we also mention Takhtajan-Teo's sewing formula for the determinants in terms of the Fredholm determinant of the welding curve. We will say more about this in Subsection \ref{Loewner subsection} (Theorem \ref{TTdetsewing}), where we discuss zeta-regularized determinants.}\\

{
\begin{remark}
 It is interesting that the quantities $S_1$ \cite{Schifferplane,Schiffer_multiply,Schiffer_rend_Fredholm,Schiffer_Pol} and $S_2$ \cite{Schiffer_Hawley_connections} were singled out by Schiffer and Schiffer-Hawley as having special geometric significance (though with greater regularity assumed). A variational formula for $S_1$ was given by Schiffer and Hawley in \cite{Schiffer_Hawley_connections}, see \cite[Remark 3.11]{Takhtajan_Teo_Memoirs}. Schiffer and Hawley also show that $S_1=S_2$. Even with the stronger regularity, the argument is not entirely complete, because it is only shown to be true under a specific set of variations.  

 The paper investigates the construction of invariant domain functionals out of kernel functions and quantities arising in potential theory, using sophisticated and inspired analogies with connections and curvature in Riemannian settings. One idea of the paper is to construct canonical mappings by minimizing functionals involving the curvature of the boundary.\footnote{This principle, with deep historical roots in complex analysis, is that there is a correspondence between extremal problems, invariant quantities, and canonical mappings.} 
 
 Like much of the work by Schiffer and his co-authors, the paper \cite{Schiffer_Hawley_connections} exploits ideas regarding symmetry and invariance in geometry and physics, but one gets the impression that the technology of the time was not yet adequate. Indeed Schiffer was knowledgeable in those subjects, having studied invariant theory under I. Schur \cite{Schiffer_Schur}, and having had a life-long interest in physics, even co-authoring a textbook on general relativity \cite{Adler_et_al_relativity}. In any case, one cannot deny the persistent relevance of his techniques and ideas for complex geometry and physics.  {An attentive reader can even find there a Polyakov-Alvarez type formula for the Fredholm determinant of a curve.} 
\end{remark}}

The invention of the concepts of \emph{overfare} and \emph{scattering} on Riemann surfaces by the authors in \cite{Schippers_Staubach_scattering}, led us to derive a formula for the K\"ahler potential of the WP-metric.
One of the basic objects in our work is the so-called {\it Schiffer operators} whose definition we now recall.

\begin{definition}
For a compact Riemann surface $\mathscr{R}$ with Green's function $\mathscr{G}(w,w_0;z,q),$ the \emph{Schiffer kernel} is defined by
 \[  L_{\mathscr{R}}(z,w) =  \frac{1}{\pi i} \partial_z \partial_w \mathscr{G}(w,w_0;z,q).    \]
Let $\mathscr{R}$ be a compact surface divided by a complex $\Gamma$ of simple closed curves into surfaces $\riem_1$ and $\riem_2$. For $k=1,2$, and $\riem_1$ and $\riem_2$ as above, we define the {\it Schiffer comparison operators} by
 \begin{align*}
  {\bf{T}}_{k}: \overline{\mathcal{A}(\riem_k)} & \rightarrow \mathcal{A}(\riem_1 \cup \riem_2)  \\
  \overline{\alpha} & \mapsto \iint_{\riem_k} L_{\mathscr{R}}(\cdot,w) \wedge \overline{\alpha(w)},
 \end{align*}
 and the {\it restriction operator}
\begin{align*}
  \mathbf{R}^0_{k}: \mathcal{A}(\riem_1 \cup \riem_2) & \rightarrow \mathcal{A}(\riem_k) \\
    \alpha & \mapsto \left. \alpha \right|_{\riem_k}.
 \end{align*}
Moreover, we define for $j,k \in \{1,2\}$
 \begin{equation}\label{defn: T sigmajsigmak}
     {\bf{T}}_{j,k} = \mathbf{R}^0_{k} \mathbf{T}_{j}: \overline{\mathcal{A}(\riem_j)} \rightarrow \mathcal{A}(\riem_k).
     \end{equation}
 \end{definition}
We also have an analogue of the Cauchy operator, namely
 \begin{definition}
     Let $\Gamma$ be a quasicircle in the Riemann sphere, whose complement has components $\riem_1$ and $\riem_2$. 
Let $\dot{\mathcal{D}}(\riem)$ denote the homogeneous Dirichlet space (that is, modulo constants). Let $\Gamma_\epsilon$ approach $\Gamma$ from within $\riem_1$.  For fixed $q \notin \Gamma$, the Cauchy operator is defined by 
\begin{align*}
 \dot{\mathbf{J}}:\dot{\mathcal{D}}(\riem_1) & \rightarrow \dot{\mathcal{D}}(\riem_2)  \\
 h & \mapsto \frac{1}{2\pi i} \lim_{\epsilon \searrow 0} \int_{\Gamma_\epsilon}  {h(\zeta)} \left( \frac{1}{\zeta-z} - \frac{1}{\zeta-q} \right) d\zeta \ \ \ z \in \riem_2  
\end{align*}
where we discard constants.
\end{definition}

 As a consequence of identities for the Schiffer operators and the properties of the Cauchy operators in \cite{Schippers_Staubach_scattering}, the authors derived the folloiwng formula for this K\"ahler potential.\\ 

\begin{theorem} \cite{Schipp-Staub-Kenig}\label{th:new_Kahler_potential}
    The potential for the Weil-Petersson metric on the Weil-Petersson Teichm\"uller space of the disk is 
    \begin{equation}
        \log \det (\mathbf{T}_{1,2}^* \mathbf{T}_{1,2})= \mathrm{log} \, \mathrm{det} \, \left( \dot{\mathbf{J}}^* \dot{\mathbf{J}} \right).
    \end{equation}
\end{theorem}

This can be compared with Theorem \ref{TTdetsewing} below. Our theorem involves the complex Cauchy operator in place of the Neumann jump operator, and does not require the assumption that the curve is $C^3$. \\

Another outcome of the scattering theory developed in \cite{Schippers_Staubach_scattering} is the relationship between the Fredholm indices of the Schiffer comparison operators and the topological invariants of the Riemann surfaces. 

\begin{theorem}[\cite{Schippers_Staubach_scattering}]
 If $\Sigma_1$, $\Sigma_2$ { are connected, and} of genus $g_1$ and $g_2$, 	then $$\mathrm{index}(\mathbf{T}_{1,2}) = g_1 - g_2 .$$\\
 If $\Sigma_2$ { is connected and} of genus $g$, and $\Sigma_1$ consists of $n$ disjoint simply connected regions, then  $$\mathrm{index}(\mathbf{T}_{1,2}) = 1-n+g.$$
\end{theorem}

Returning to the K\"ahler property of the WP-Metric, in the paper \cite{Yanagishita_Kahler} of M. Yanagishita obtained the following generalization of Takhtajan and Teo's result. 
\begin{theorem}[\cite{Yanagishita_Kahler}] Let $G$ be a Fuchsian group satisfying Lehner's condition. Then the Weil-Petersson metric on $T_{\mathrm{WP}}(\disk_-/G)$ is K\"ahler. 
\end{theorem}
Furthermore
\begin{theorem}[\cite{Yanagishita_Kahler}] Let $G$ be a Fuchsian group satisfying Lehner's condition. Then the holomorphic sectional and Ricci curvatures of the Weil-Petersson metric on $T_{\mathrm{WP}}(\disk_-/G)$ are negative.
\end{theorem}  
 
\end{subsection}

  \begin{subsection}{Chern classes, Quillen metric and zeta functions}
  \label{se:Chern_Quillen_Zeta}
     Next we discuss the connection of the WP-Teichm\"uller theory to Chern-Weil theory. In this context, Wolpert \cite{Wolpert} 
computed the first Chern form \(\mathrm{c}_{1}(v)\) of the vertical tangent bundle of the universal Teichmüller curve, with the Hermitian metric induced by the
hyperbolic metric on fibers. Then since integration of powers of \(\mathrm{c}_{1}(v)\) over fibers produces characteristic classes \(\kappa_{n}\) (see \eqref{Mumford} ahead) on Teichmüller space, Wolpert showed that
that \(\kappa_{1}\) is a multiple of the Weil-Petersson Kähler form \(\omega_{\mathrm{WP}}\).\\

    Another important and interesting topic is the Quillen metric on the so-called determiannt line bundle and its connection to WP-Teichm\"uller theory. To descibe this, we consider a holomorphic family $\pi : X\to P$ of compact complex manifolds over a compact base $P$. Let also $E$ be a vector bundle over $X$, with the projection $\pi': E\to X.$\\

Given $p\in P$ we identify the fiber $\pi^{-1}(p)$ with $X_p$ and restrict $\pi'$ to a map $\pi'_p: E_p \to X_p$. Now there is a family of $\overline{\partial}$-operators, denoted by $\overline{\partial}_p$ acting on the sections of the bundle $(E_p, \pi', X_p )$. This family defines an index bundle $\mathrm{ind}\,\overline{\partial}$ over $P$ in the sense of K-theory, where $\mathrm{ind}\,\overline{\partial}= \mathrm{Ker}\,\overline{\partial}-\mathrm{Coker}\,\overline{\partial}.$\\

If $\mathrm{ch}(\cdot)$ denotes the Chern character of a bundle and $\mathrm{td}(T_{\mathrm{vert}}X)$ denotes the Todd class of the vertical tangent bundle of $X$, then as a result of the Atiyah-Singer family index theorem \cite{AtiyahSinger4}, one has that

\begin{equation}\label{atiyahsinger familyindex}
  \mathrm{ch}(\overline{\partial})=\pi_{\int}[ \mathrm{ch}(E)\cdot\,\mathrm{td}(T_{\mathrm{vert}}X)],
\end{equation}

where $\pi_{\int}[\cdot] : H^{*}(X)\to H^{*-\mathrm{dim}\, \pi^{-1}(p)} (P)$ is the operation of integration along the fibers. If the bundles $T_{\mathrm{vert}}X$ and $E$ are also equipped with  Hermitian metrics, then these bundles will carry canonical unitary connections compatible with their metrics and therefore, using these connections, $\mathrm{ch}(E)$ and $\mathrm{td}(T_{\mathrm{vert}}X)$ will be represented by closed differential forms on $X$.\\

 Furthermore, if the base $P$ is not compact and the Atiyah-Singer family index theorem \eqref{atiyahsinger familyindex} doesn't apply, then if the index bundle $\mathrm{ind}\,\overline{\partial}$ happens to be a vector bundle (this is not always true because the dimensions of the fibers $\mathrm{Ker}\,\overline{\partial}_p$ and $\mathrm{Coker}\,\overline{\partial}_p$ may jump), one could hope for a local version of \eqref{atiyahsinger familyindex}, called a \emph{local index theorem}.   This is obtained by using a particular connection on the vector bundle $\mathrm{ind}\,\overline{\partial}$ (if one such exists), so that \eqref{atiyahsinger familyindex} holds as an equality between the corresponding differential forms of its left- and right-hand sides.\\

 As it turns out, proving a local index theorem is intimately connected to the problem of expressing the Chern character and the Todd class explicitly in terms of the maps $\pi$ and $\pi'$ defined above. But one confronts complications in producing local index theorems if $\mathrm{ind}\,\overline{\partial}$  is not a vector bundle. So to remedy this situation, D. Quillen \cite{Quillen} suggested the study of the \emph{determinant line bundle}, which in our case of study is defined by
$$\mathrm{det\, ind }\,\overline{\partial} = \wedge^{\mathrm{max}}\mathrm{Ker}\,\overline{\partial}\otimes (\wedge^{\mathrm{max}}\mathrm{Coker}\,\overline{\partial})^{-1},$$
and is a holomorphic line bundle over $P$ (this of course requires the validity of certain conditions on the bundles $(X, \pi, P )$ and $(E, \pi', X )$ ). If $\Vert\cdot\Vert$ denotes the $L^2$ norm on $\mathrm{det\, ind }\,\overline{\partial}$ induced by the metrics on $E$ and $T_{\mathrm{vert}}X$ and if $\overline{\partial_p}^{*}\overline{\partial_p}$ is the associated Laplacian acting on the sections of the bundle $(E_p, \pi', X_p )$, then $\mathrm{det}\,\overline{\partial_p}^{*}\overline{\partial_p}$ is given by the Laplacian's zeta function determinant regarded as a function on $P$, and one can define the \emph{Quillen metric} as
$$\Vert\cdot\Vert_{Q}= \frac{\Vert\cdot\Vert}{(\mathrm{det}\,\overline{\partial}^{*}\overline{\partial})^{1/2}}.$$

One of the interesting features of Quillen's metric is that the curvature form (first Chern class) of the line bundle $(\mathrm{det\, ind }\,\overline{\partial} ;\Vert\cdot\Vert_{Q})$ on $P$ is up to a multiplicative constant, the canonical K\"ahler form on $P$.

Now in order to obtain results that are interesting from geometric and global analytic point of view in Teichm\"uller theory, the choices of $P$ and $X$ and a specific bundle $E$ over $X$ are very significant. For example in the abstract setting mentioned above let us choose $P= T_g $ i.e. the Teichm\"uller space of compact genus $g$ Riemann surfaces, $X=\mathscr{T}_g$ i.e. the Bers fiber space of Teichm\"uller curves, and $E=T_{\mathrm{vert}}^{-s}\mathscr{T}_g $ i.e. the $-s\,$th power of the vertical holomorphic line bundle of $\mathscr{T}_g$ (the sections of this bundle are the holomorphic $(-s,0)$-differentials). Once again, there is a family of $\overline\partial$ operators associated to this setting, which is denoted here by $\overline{\partial}_{s}$ and act on the $(s,0)$-differentials on Riemann surfaces. In this situation, A. A. Belavin and V.G. Knizhnik \cite{Belavin} showed that the first Chern form of the corresponding $\mathrm{det\, ind }\,\overline{\partial}_{s}$ can be expressed as
\begin{equation}\label{Belavin-Knizhnik index}
  \mathrm{c}_{1}(\mathrm{det\, ind }\,\overline{\partial}_{s})=\int_{\pi^{-1}(p)}( \mathrm{ch}(T_{\mathrm{vert}}^{-s}\mathscr{T}_g)\cdot\,\mathrm{td}(T_{\mathrm{vert}}\mathscr{T}_g))_{2,2},
\end{equation}
where $(\cdot)_{2,2}$ denotes the $(2,2)$ component of a differential form on the Teichm\"uller curve, and $\pi^{-1}(p)$, for $p\in T_g $, denotes the fibers of the bundle $\pi: \mathscr{T}_g\to T_g$. Formula \eqref{Belavin-Knizhnik index} can be regarded as a local index theorem because it establishes an equality between the $(1,1)$ forms of the left and right hand sides of \eqref{Belavin-Knizhnik index}. Now if one endows $T_{\mathrm{vert}}^{-s}\mathscr{T}_g$ with a metric in such a way that for every $p\in T_g$ each $\pi^{-1}(p)$ would be equipped with the Poincar\'e metric, then as was shown by S. Wolpert \cite{Wolpert}, one has the following analogue of the Quillen's result, namely

\begin{equation}\label{Wolpert index}
  \mathrm{c}_{1}(\mathrm{det\, ind }\,\overline{\partial}_{s})=\frac{6s^2 -6s+1}{12\pi^2}\,\omega_{\mathrm{WP}},
\end{equation}
where $\omega_{\mathrm{WP}}$ is the Weil-Petersson K\"ahler form on  $T_g (\Sigma)$.
\\

If one instead takes $E=T_{\mathrm{vert}}\mathscr{T}(\disk_-)$, then the hyperbolic metric on $w^{\mu} (\mathbb{D}_{-})$ defines a Hermitian metric on $T_{\mathrm{vert}}\mathscr{T}(\disk_-)$, as we mentioned earlier this can be used to produce a canonical unitary connection on $T_{\mathrm{vert}}\mathscr{T}(\disk_-)$, compatible with the metric, such that $\mathrm{ch}(T_{\mathrm{vert}}\mathscr{T}(\disk_-))$ and $\mathrm{td}(T_{\mathrm{vert}}\mathscr{T}(\disk_-))$ will be represented by closed differential forms on $\mathscr{T}(\disk_-)$. Thus calculating the curvature form of this connection and using integration on the fibers, one hopes to connect fiber integrals of powers of the curvature form $\Theta$ (first Chern class) to some invariant quantities, such as the metric.\\

To this end, following Wolpert \cite{Wolpert}, Takhtajan and Teo \cite{Takhtajan_Teo_Memoirs}  considered the \emph{Miller-Morita-Mumford characteristic forms}. These are $(n,n)$-forms on
the Hilbert manifold $T(\disk_-)$, defined by
\begin{equation}\label{Mumford}
  \kappa_n = (-1)^{n+1} \pi_{\int} [(\mathrm{c}_{1}(T_{\mathrm{vert}}\mathscr{T}(\disk_-)))^{n+1}],
\end{equation}
where $\pi_{\int} : H^{*}(\mathscr{T}(\disk_-))\to H^{*-2}(T(\disk_-))$ is the operation of integration along the fibers. Observe that, \eqref{Mumford} yields for instance that $\kappa_1$ can be calculated by fiber-integration of the square of the curvature form $\Theta$.
Takhtajan and Teo showed the following.

The Miller-Morita-Mumford characteristic forms \(\kappa_{n}\) are right-invariant on the Hilbert manifold \(T(\mathbb{D}_-)\)
and for \(\mu_{1}, \ldots, \mu_{n}, \nu_{1}, \ldots, \nu_{n} \in H_{-1,1}\left(\mathbb{D}_{-}\right) \simeq T_{[0]} T(\mathbb{D}_-)\)
one has that
\begin{equation*}
   \kappa_{n}\left(\mu_{1}, \ldots, \mu_{n}, \bar{\nu}_{1}, \ldots, \bar{\nu}_{n}\right) ==\frac{i^{n}(n+1) !}{(2 \pi)^{n+1}} \sum_{\sigma \in S_{n}} \operatorname{sgn}(\sigma) \iint_{\mathbb{D}_{-}} G\left(\mu_{1} \bar{\nu}_{\sigma(1)}\right) \ldots G\left(\mu_{n} \bar{\nu}_{\sigma(n)}\right) dA_{\mathrm{hyp}},
\end{equation*}
where $G$ is given by \eqref{resolvent}.
Moreover, if $T(\disk_-)$ is equipped with its Hilbert manifold structure then
\begin{equation}\label{TT index 1}
 \kappa_1=\frac{1}{\pi^2} \omega_{\mathrm{WP}}.
\end{equation}

In \cite {TakhtajanZograf} L. Takhtajan and P. Zograf  proved a local index theorem for a family of \(\bar{\partial}\)-operators on a type \((g, n)\) punctured surface $\Sigma$, i.e., a compact surface of genus \(g\) with \(n\) punctures. 
This case corresponds to the general construction above with the following choices, namely take \(P=T_{g, n}\) i.e. the Teichmüller space of punctured Riemann
surfaces of type \((g, n)\), \(X=\mathscr{T}_{g, n}\) i.e. the
corresponding universal family (so that the fibers of the fibration \(\pi: \mathscr{T}_{g, n} \rightarrow T_{g, n}\)
are Riemann surfaces of type \((g, n)\)), and \({E}=T_{\mathrm{vert}}^{-k} \mathscr{T}_{g, n}\) i.e. the \(k^{\text {th }}\) power of the
vertical line bundle on \(\mathscr{T}_{g,n}\).\\

Now, in contrast to the case of compact surfaces, the Laplace operator \(\Delta_{k}=\bar{\partial}_{k}^{*} \bar{\partial}_{k}\)
associated with the Poincar\'e metric on punctured surfaces of type \((g, n)\) has a continuous spectrum.\\
In connection to this fact, let us also recall the definition of Selberg's Zeta function which for \(\operatorname{Re} s>1\) is given by the absolutely convergent product
$$
Z(s)=\prod_{\{\ell\}} \prod_{m=0}^{\infty}\left(1-e^{-(s+m)|\ell|}\right),
$$
where \(\ell\) runs over the set of all simple closed geodesics on \(\Sigma\) with respect
to the Poincare metric, and \(|\ell|\) is the length of \(\ell\). The function \(Z(s)\) admits a
meromorphic continuation to the whole complex s-plane with a simple zero at
\(s=1\).

In the case of compact Riemann surfaces
E. D'Hoker and D. Phong \cite {DHokerPhong} showed that, the determinant
of \(\Delta_{k}\) is defined via its Selberg zeta function \(Z(s)\) up to a constant multiplier depending
only on \(g\) and \(k\), and is equal to \(Z^{\prime}(1)\) for \(k=0,1\) and to \(Z(k)\) for \(k \geq 2\). The same goes also for punctured Riemann surfaces of type \((g, n)\).\\

 Using this, Takhtajan and Zograf calculated the first Chern form of the determinant
line bundle $\det \mathrm{ind}\, \bar{\partial}_{k}$ on the Teichmüller space \(T_{g, n}\) of type $(g,n)$ punctured surfaces, endowed with the Quillen metric, and showed that 
\begin{equation}
\mathrm{c}_{1}(\det \mathrm{ind}\, \bar{\partial}_{k})=\frac{6k^2-6k+1}{12\pi^2}\omega_{\mathrm{WP}}-\frac{1}{9} \omega_{\mathrm{cusp }}.
\end{equation}
Here we see that the result differs from \eqref{Wolpert index} by an additional term, which turns out to be the Kähler form of a new
Kähler metric on the moduli space of punctured Riemann surfaces of type $(g, n)$. In fact \(\omega_{\text {cusp }}\) is the symplectic form of a Kähler metric \(\langle,\rangle_{\text {cusp }}\) on \(T_{g, n} (\Sigma)\)
which is invariant with respect to the modular group on \(T_{g, n} (\Sigma)\). 

The metric \(\langle,\rangle_{\text {cusp }}\)
is defined by means of the Eisenstein-Maa{\ss} series related to \(n\) punctures. In order to give a precise definition of the cusp-metric, consider a Riemann surface \(\Sigma\) of type \((g, n)\) and assume that it is
equipped with the Poincaré metric \(\varrho\) and let \(\Gamma\) be a torsion-free Fuchsian group
 such that $\Sigma= \mathbb{H}/\Gamma$ where \(\mathbb{H}\) is the open upper half-plane. Denote by \(\Gamma_{1}, \ldots, \Gamma_{n}\) the set of non-conjugate parabolic subgroups in \(\Gamma\),
and for every \(i=1, \ldots, n\) fix an element \(\sigma_{i} \in\) PSL \((2, \mathbb{R})\) such that \(\sigma_{i}^{-1} \Gamma_{i} \sigma_{i}=\Gamma_{\infty}\)
where the group \(\Gamma_{\infty}\) is generated by the parabolic transformation \(z \mapsto z+1\). The Eisenstein-Maa{\ss} series \(E_{i}(z, s)\) corresponding to the \(i^{\text {th }}\) cusp of the group \(\Gamma\) is
defined for \(\operatorname{Re} s>1\) by the formula
$$
E_{i}(z, s)=\sum_{\gamma \in \Gamma/\Gamma_{i}} \operatorname{Im}\left(\sigma_{i}^{-1} \gamma z\right)^{s}, \quad i=1, \ldots, n
$$
It is a well-known fact that $E_{i}(z, s)$ can be meromorphically continued to the whole complex s-plane, and for
\(\operatorname{Re} s=\frac{1}{2}\) the Eisenstein-Maa{\ss} series \(E_{i}(z, s), i=1, \ldots, n\), form a complete set of
eigenfunctions of the continuous spectrum of the Laplace operator \(\Delta_{0}\). 
Now one defines the cusp-metric as follows; first set 
\begin{equation}
\langle\mu, \nu\rangle_{i}=\int_{\Sigma} \mu \bar{\nu}\, E_{i}(\cdot, 2)\, \quad \mu, \nu \in L^2_{-1,1}(\Sigma), \quad i=1, \ldots, 
\end{equation}

and it turns out that each scalar product $\langle,\rangle_{i}$ gives rise to a Kähler metric on $T_{g, n}$. The cusp-metric is the sum
\begin{equation}
\langle,\rangle_{\operatorname{cusp}}=\sum_{i=1}^{n}\langle,\rangle_{i}
\end{equation}
which is a K\"ahler metric invariant under the Teichmüller
modular group $\mathrm{Mod}_{g, n}.$
\\

In \cite{RatiuTodorov} T. Ratiu and A. Todorov constructed a determinant line bundle on $\text{Diff}^{+}(\mathbb{S}^1)/\text{M\"ob} (\mathbb{S}^1)$, together with the corresponding Quillen norm. They further proposed that, using Quillen's original construction of the determinant line bundle in \cite{Quillen}, the curvature \((1,1)\)-form of this Quillen metric is the Kähler form of the Weil-Petersson metric on $\text{Diff}^{+}(\mathbb{S}^1)/\text{M\"ob} (\mathbb{S}^1)$, which would therefore be the infinite dimensional analog of Takhtajan-Zograf's result mentioned earlier. 
In connection to this circle of ideas, we would also mention the paper by A. Fujiki and G. Schumacher \cite{FujikiSchumacher}, where the authors 
used a significant generalization of Quillen's work due to J. M. Bismut and D. Freed \cite{BismutFreed} and Bismut, H. Gillet and Ch. Soul\'e \cite{BismutGilletSoule} (concerning the first Chern form of a determinant bundle with Quillen metric), and a certain fiber integration formula, to construct a generalized Weil-Petersson form and to show that it is (up to a constant) the Chern form of a Hermitian line bundle on the moduli space.
\end{subsection}

\end{section}

\begin{section}{Weil-Petersson beyond Teichm\"uller theory} \label{se:in_the_wild} 
\begin{subsection}{Conformal field theory and string theory}   \label{se:conformal_field_theory}

  In two-dimensional conformal field theory, D. Friedan and S. Shenker pointed out the role of moduli spaces of Riemann surfaces with punctures \cite{FriedanShenker}. C. Vafa gave a formulation of conformal field theories using a moduli space of punctured surfaces endowed with local coordinates vanishing at the punctures \cite{Vafa}.  Note that the ``coordinates'' are extra data beyond that of the conformal equivalence class of the surface. 
  Segal sketched a definition of conformal field theory in \cite{Segal87} and \cite{SegalPreprint}; the latter was a well-circulated preprint which was eventually published in \cite{SegalPublished}. In this definition, Riemann surfaces with boundary parametrizations (up to conformal equivalence) are the morphisms of a category whose objects are disjoint unions of circles (some labelled incoming and some labelled outgoing). Such surfaces can be sewn using the parametrizations.  A heuristic correspondence between the moduli spaces of Vafa and Segal can be obtained if one sews on punctured disks using the boundary parametrizations, and the inclusion map becomes the inverse of the local coordinates.  
\begin{remark}
    This heuristic correspondence between the two moduli spaces is one of the motivations for the fiber model of the complex structure on Teichm\"uller space. The idea is that the boundary values of a quasiconformal deformation $F:\riem \rightarrow \riem_1$ of a bordered surface $\riem$ induce boundary parametrizations in the Segal model. In Vafa's model, this data takes the form of elements of $\mathcal{O}^{\mathrm{qc}}(\riem_1^P)$ for the corresponding compact surface $\riem^P_1$ obtained by sewing on caps.
\end{remark}
  
  In Segal's definition, a conformal field theory associates an operator between tensor products of Hilbert spaces associated to the incoming and outgoing boundary components. In categorical language, this is supposed to be a projective holomorphic functor. The association of the operator to the Riemann surface is required to be holomorphic, and thus one requires a complex structure on the moduli space and holomorphic sewing operation.  

  Radnell-Schippers \cite{RadnellSchippers_monster} showed that the moduli space can be identified with a quotient of Teichm\"uller space. To do this, the set of boundary parametrizations was enlarged to the set of quasisymmetric parametrizations; in the Vafa model, the extra coordinate data were conformal maps of the disk with quasiconformal extensions. 
  
  More precisely, and with alternate analytic choices, we make the following definitions.
\begin{definition} 
 Consider pairs $(\riem^P,f)$ where $\riem^P$ is a compact Riemann surface of genus $g$ with $n$ labelled punctures $p=(p_1,\ldots,p_n)$ and  $f=(f_1,\ldots,f_n) \in \mathcal{O}^{\mathrm{qc}}({\riem}^P)$. We say that two pairs $(\riem^P_1,f_1)\sim (\riem^P_2,f_2)$ are equivalent if there is a biholomorphism $\sigma:\riem^P_1 \rightarrow \riem^P_2$ such that $\sigma \circ f_1 = f_2$.  The puncture model of the rigged moduli space is 
 \[ \widetilde{\mathcal{M}}_P(g,n) = \{ (\riem^P,f) \}/\sim.   \]
\end{definition}  
The subscript $P$ stands for ``puncture'' model.  
The rigged moduli space in the Segal border model is as follows. 
\begin{definition} Consider pairs $(\riem,\phi)$ where $\riem$ is a bordered surface of type $(g,n)$ with ordered boundary curves $\partial_k \riem$, $k=1,\ldots,n$ and $\phi=(\phi_1,\ldots,\phi_n)$ is a collection of
quasisymmetries $\phi_k:\mathbb{S}^1 \rightarrow \partial_k \riem$. We say that two pairs $(\riem_1,\phi) \sim (\riem_2,\psi)$ are equivalent if there is a biholomorphism $\sigma:\riem_1 \rightarrow \riem_2$ such that $\psi=\sigma \circ \phi$. The border model of the rigged moduli space is 
\[ \widetilde{\mathcal{M}}_B(g,n) = \{ (\riem,\phi) \}/\sim.   \]
\end{definition}
These are in one-to-one correspondence \cite{RadnellSchippers_monster}. 

This relates to the Teichm\"uller space as follows. Recall that $\mathrm{ModI}$ is the subgroup of the modular group, whose elements fix the boundary pointwise.
\begin{theorem}[\cite{RadnellSchippers_monster}] \label{th:RMS_Teich_corresp} 
 Let $\riem$ be a bordered surface of type $(g,n)$. Then there is a bijection  
 \[  T(\riem)/\mathrm{ModI} \cong \widetilde{\mathcal{M}}_B(g,n) \cong \widetilde{\mathcal{M}}_P(g,n).   \]

 $\widetilde{\mathcal{M}}_B(g,n)$ and  $\widetilde{\mathcal{M}}_P(g,n)$ are complex Banach manifolds, and sewing is holomorphic. 
\end{theorem}

{
\begin{remark} In some citations of \cite{RadnellSchippers_monster}, it is stated that  Radnell-Schippers extend the sewing operation to quasisymmetries, but this is a routine consequence of conformal welding. The interesting fact is rather the geometric connection between the Teichm\"uller space and the rigged moduli spaces. This brings insights in both directions: for example in Teichm\"uller theory, the CFT point of view leads directly to the fiber model of Teichm\"uller space. 
\end{remark}}

Customarily the rigging is chosen to be more regular; for example the boundary parametrizations are chosen to be analytic or smooth. The question is then what is the appropriate regularity for conformal field theory? The case for Weil-Petersson class parametrizations was advocated at length in \cite{RadnellSchippersStaubach_CFT_review}. We mention here four motivations for choosing Weil-Petersson class moduli space. 

The first two of these motivations are simple to state: $\mathrm{QS}_{\mathrm{WP}}(\mathbb{S}^1)$ is the completion of $\mathrm{Diff}(\mathbb{S}^1)$ (Theorem \ref{th:WP_universal_Teich_complete}), and, $\mathrm{QS}_{\mathrm{WP}}(\mathbb{S}^1)$ is a topological group (Theorem \ref{th:WP_is_top_group}).  

The remaining two motivations relate to the determinant line bundle of $\overline{\partial} \oplus \mathrm{pr}$, which is part of the construction of certain CFTs satisfying the definition of Segal \cite{Huang}.  The definition of Segal involves the polarization of functions on the boundary $\partial \riem$ induced by positive and negative Fourier modes of the pull-back under the boundary parametrization. It is required that this polarization be an element of the {\it restricted} Siegel disk. Thus the boundary parametrization must be Weil-Petersson.  
Finally, the construction of the central extension of the group $\mathrm{Diff}(\mathbb{S}^1)$ and computation of the central charge requires that $\mathrm{Diff}(\mathbb{S}^1)$ be in the restricted group \cite[Appendix D]{Huang}. One can enlarge the group, but only if its action is still in the Shale group. 
{We are grateful to Yi-Zhi Huang for explaining this point to us.} Again, this is precisely the condition that the parametrization be a Weil-Petersson quasisymmetry. 

One can indeed substitute Weil-Petersson class objects in the definitions of the rigged moduli spaces, and obtain the following. 
\begin{theorem}[\cite{RSS_WP2}] \label{th:RMS_Teich_corresp_WP}
 Let $\riem$ be a bordered surface of type $(g,n)$ for $2g+2-n>0$. Then there is a bijection  
 \[  T_{\mathrm{WP}}(\riem)/\mathrm{ModI} \cong \widetilde{\mathcal{M}_{\mathrm{WP}}}_B(g,n)  \cong \widetilde{\mathcal{M}_{\mathrm{WP}}}_P(g,n).   \]

 Thus $\widetilde{\mathcal{M}_{\mathrm{WP}}}_B(g,n)$ and  $\widetilde{\mathcal{M}_{\mathrm{WP}}}_P(g,n)$ are complex Hilbert manifolds. 
\end{theorem}

  As we touched on in Section \ref{se:period}, the diffeomorphism group of the circle plays a role in representation theory. It also appears in various physical models. We indicate some literature as it connects to Weil-Petersson Teichm\"uller theory. 
  
  In string theory, $\mathrm{Diff}(\mathbb{S}^1)$ is the group of reparametrizations of a string, and invariance under this reparametrization can be required in different settings.  The Lie algebra of $\mathrm{Diff}(\mathbb{S}^1)$ can be identified with the smooth vector fields on the circle. The generators of this Lie algebra are $-i e^{i n \theta} \partial/\partial \theta$. The complexification is generated by elements of the form $z^{n-1}\partial/\partial z$, $n \in \mathbb{Z}$.   The vector space of polynomial expressions in these generators, with the vector field Lie bracket, is called the Witt algebra.   
  The Witt algebra is viewed in the physics literature as the conformal symmetry group \cite{BPZ,Schottenloher}\footnote{This symmetry group is twice as large as one might expect from a differential-geometric point of view of local symmetries; local biholomorphisms fixing a point are generated by $z^n \partial/\partial z$ with $n\geq 1$.}. Although the group $\mathrm{Diff}(\mathbb{S}^1)$ has a Lie algebra $\mathrm{Vect}(\mathbb{S}^1)$, there is no complexification. That is, there is no Lie group whose Lie algebra is the Witt algebra (or a completion) \cite{Schottenloher}. 
  
  Segal \cite{Segal87,SegalPreprint} argues that, the object which is closest to filling the role of the complexification of $\mathrm{Diff}(\mathbb{S}^1)$, is the rigged moduli space of the twice punctured sphere with coordinates. Equivalently, it is the rigged moduli space of the annulus with parametrizations. This is a semigroup, since sewing preserves the type $(0,2)$ and thus induces a multiplication, but there is no inversion. If one uses quasisymmetric parametrizations or equivalently elements of $\mathcal{O}^{\mathrm{qc}}(\sphere \backslash \{0,\infty \})$ for local coordinates, one obtains from Theorem \ref{th:RMS_Teich_corresp} that this semigroup can be identified with the Teichm\"uller space of an annulus $\mathbb{A}$ modulo a $\mathbb{Z}$ action by boundary twists, that is  
  \[  \widetilde{M}_B(0,2) \cong \widetilde{M}_P(0,2) \cong T(A)/\mathbb{Z}     \] 
  (see \cite{RadnellSchippers_annulus}), and sewing in this semigroup is holomorphic with respect to the complex structure obtained from $T(\mathbb{A})$. By Theorem \ref{th:RMS_Teich_corresp_WP} one could replace the moduli spaces above by their Weil-Petersson versions.
  \\ 

  It has been suggested that the universal Teichm\"uller space might play a role in a non-perturbative approach to string theory \cite{BowickRajeev,BowickRajeev_short,HongRajeev,Pekonen}.  Bowick-Rajeev \cite{BowickRajeev_short} consider the loop space in $\mathbb{R}^{d-1,1}$ space-time, and $\mathrm{Diff}(\mathbb{S}^1)/\mathbb{S}^1$ appears in association to reparametrization invariance and changes of complex structure (see Section \ref{se:polarizations_Siegel}). They state ``An approach
to string theory also based on complex geometry has
been proposed by Friedan and Shenker \cite{FriedanShenker}. 
It is of interest to know the relation between the `universal
Teichm\"uller space' with which they work and $\text{Diff}(\mathbb{S}^1)/\mathbb{S}^1$.''  Note that this remark does not actually refer to the universal Teichm\"uller space: Friedan and Shenker use the term ``universal moduli space'' and argue that this ``cannot be the universal Teichm\"uller space described in the mathematics literature'' - that is, $T(\disk_-)$.   

On the other hand, in \cite{HongRajeev}, the connection of $\text{Diff}(\mathbb{S}^1)/\mathbb{S}^1$ to $T(\disk_-)$ is observed by Hong-Rajeev. A kind of sum over paths is taken over only the compact surfaces, and the fact that the universal Teichm\"uller space contains the moduli spaces of compact Riemann surfaces is advanced as evidence for the idea of a non-perturbative approach using $\text{Diff}(\mathbb{S}^1)/\mathbb{S}^1$, and a possible action is advanced. The relevance to Weil-Petersson geometry is that the embedding into the infinite restricted Siegel disk of polarizations of Segal plays an explicit role.  Indeed Hong-Rajeev (referring to \cite{NagVerjovsky}) point out that the lack of completeness of $\mathrm{Diff}(\mathbb{S}^1)/\mathbb{S}^1$ is a technical obstacle. The results of Cui and Takhtajan{-Teo}  
(Theorem \ref{th:WP_universal_Teich_complete}) have removed this obstacle. In this context, the result of Shen and Takhtajan-Teo that an element of the Teichm\"uller space $p \in T(\disk_-)$ maps into the restricted Siegel disk if and only if $p \in T_{\mathrm{WP}}(\disk_-)$ (Theorem \ref{th:symp_res_is_WP}), is very satisfying. 

\begin{remark}
Hong-Rajeev \cite{HongRajeev} state that $\mathrm{Diff}(\mathbb{S}^1)/\mathbb{S}^1$ is a dense subset of the universal Teichm\"uller space. Setting aside different interpretations of  ``dense'', it is rather the quotient $\mathrm{Diff}(\mathbb{S}^1)/\text{M\"ob}(\mathbb{S}^1)$ which is a subset of the universal Teichm\"uller space $T(\disk_-) \cong \text{QS}(\mathbb{S}^1)/\text{M\"ob}(\mathbb{S}^1)$. However, L.P. Teo \cite{TeoVelling} showed that the universal Teichm\"uller curve $\mathscr{T}(\disk_-)$ (that is, the Bers fiber space over $T(\disk_-)$ mentioned earlier in Section \ref{se:Chern_Quillen_Zeta}), can be naturally identified with $\text{QS}(\mathbb{S}^1)/\mathbb{S}^1$. 
Takhtajan and Teo \cite{Takhtajan_Teo_Memoirs} showed that $\text{QS}_{\text{WP}}(\mathbb{S}^1)/\mathbb{S}^1$ is a Hilbert manifold and a topological group, and in fact $\mathrm{Diff}(\mathbb{S}^1)/\mathbb{S}^1$ is a dense subset of the Weil-Petersson universal Teichm\"uller curve in its associated topology. This could be seen to revive the idea of Hong and Rajeev.  
\end{remark}

In regards to the sum over compact surfaces, the universal hyperbolic lamination is also of interest and might provide another approach; see \cite{Sullivan_lamination,NagSullivan,Burgos_Verjovsky}. 
\end{subsection}
\begin{subsection}{Fluid mechanics}
A topic of interest in mathematical fluid dynamics is the study of geometric aspects of $\text{Diff}(\mathbb{S}^1)/\text{M\"ob}(\mathbb{S}^1)$.  M. Schonbek A. Todorov and J. Zubelli \cite{SchonTodZub} consider this in association with the the polarizations of Segal \cite{SegalUnitary} and Nag-Sullivan \cite{NagSullivan} and the emebedding into the Sato-Segal-Wilson Grassmannian \cite{SegalWilson}.  
A Hilbert manifold structure with tangent spaces modelled on $H^{3/2}$ is proposed there.  Motivated by Arnol'd's approach to Euler's equation, {as well as Segal-Wilson \cite{SegalWilson}},  Schonbek, Todorov and Zubelli 
utilized the remarkable fact that the geodesics of the Teichm\"uller space give information about the behaviour of solutions to the
periodic Korteweg-de Vries (KdV) equation, and thereby proved existence of periodic
solutions to the KdV equation with initial data in the Sobolev space on the unit circle. They thereafter used infinite-dimensional 
Cartan-Hadamard theory, to prove the Arnol'd exponential instability of the geodesic flow.\\

With the manifold structure of Cui and Takhtajan-Teo described in  Section \ref{se:WP_universal_Cui_Takhtajan}, 
the tangent space at the identity is the space of functions on the circle of class \(H^{\frac{3}{2}}\), in other words, the {Teichm\"uller space \(T_{\mathrm{WP}}(\mathbb{D}_-)\)} has a Hilbert manifold structure whose tangent space can be identified
with the \(H^{3 / 2}\left(\mathbb{S}^{1}\right)\)-vector fields. In \cite{Figalli} A. Figalli characterized the flows generated by the \(H^{3 / 2}\left(\mathbb{S}^{1}\right)\)-vector fields in terms of
fractional Sobolev norms. More precisely, he showed that a flow generated by an \(H^{3 / 2}\left(\mathbb{S}^{1}\right)\)-vector
field belongs to Sobolev space \(W^{1, p}\left(\mathbb{S}^{1}\right)\) for all \(p \geq 1\) and belongs to \(W^{1+r, q}\left(\mathbb{S}^{1}\right)\) for all \(0<r<1 / 2\) and
\(1 \leq q<1 / r\). However Figalli also showed that there exists an autonomous vector field \(u \in H^{3 / 2}\left(\mathbb{S}^{1}\right)\) such that its
flow map is neither Lipschitz nor \(W^{1+r, 1 / r}(\mathbb{S}^{1})\) for all \(r \in(0,1)\). In particular taking $r=1/2$ , one sees that there is a flow generated by an \(H^{3 / 2}\left(\mathbb{S}^{1}\right)\)-vector field that is neither in \(W^{3/2, 2}(\mathbb{S}^{1})=H^{3 / 2}(\mathbb{S}^{1})\), nor in the Lipschitz space.\\

By the aforementioned work of Takhtajan and Teo, the identity component of universal Teichm\"uller space takes the place of
diffeomorphisms of critical Sobolev class \(H^{\frac{3}{2}}\).

In \cite{GMR} and \cite{GR} F. Gay-Balmaz, J. Marsden and T. Ratiu studied this group from the point of view of
manifolds of maps by identifying it with a subgroup of the quasisymmetric homeomorphisms of the circle and proved that all elements of this group are of class \(H^{\frac{3}{2}-\varepsilon}\) for all
\(\varepsilon>0\). Another outcome of this was that for WP-class quasisymmetries $\phi$ on $\mathbb{S}^1$, one has that both $\phi$ and its inverse $\phi^{-1}$ are in $H^{\frac{3}{2}-\varepsilon}(\mathbb{S}^1)$ for all $\varepsilon >0$. Now since the Weil-Peterson metric induces the Hilbert space topology on each tangent space, Gay-Balmaz and Ratiu were also able to show that the metric is complete i.e. all geodesics of the Weil-Petersson metric exist for all time. They also proved that this space is Cauchy complete relative to the distance function defined by the Weil-Petersson metric, which earlier had been shown by Cui \cite{Cui}.\footnote{Because of the Hopf-Rinow theorem, given Cauchy completeness one expects geodesic completeness. However this is not obvious in the infinite-dimensional setting.} 
In connection to the Sobolev regularity mentioned above, Shen \cite{Shen_characterization} proved that there exists some quasisymmetric homeomorphism of the Weil-Petersson class
which belongs neither to the Sobolev space \(H^{\frac{3}{2}}\) nor to the Lipschitz space \(\Lambda^{1}\). Based on
this new characterization of the Weil-Petersson class, Shen introduced a metric on the Weil-Petersson Teichmüller space and used it to give a new proof of Takhtajan-Teo's result that the Weil-Petersson Teichmüller space is a topological group.

Moreover {Gay-Balmaz and Ratiu} obtained a new equation, the so-called {\it{Euler-Weil-Petersson equation}}, and showed that the solutions are \(C^{0}\) in \(H^{\frac{3}{2}}\)
and \(C^{1}\) in \(H^{\frac{1}{2}}\).  Another achievement of the Weil-Petersson techniques in the paper \cite{Gay}, based specifically on the long-time existence of Weil-Petersson geodesics, was the solution to a problem posed by E. Sharon and D. Mumford in \cite{SharonMumford}, regarding the existence of a unique geodesic between two shapes in the plane.
     \end{subsection}
     \begin{subsection}{Loewner energy}\label{Loewner subsection}
The Loewner energy of a Jordan curve is the Dirichlet energy of its Loewner driving term, and was introduced by P. Friz and A. Shekar \cite{FrizShekar} and Y. Wang \cite{Wang}. More precisely let \(\eta\) be a Jordan curve on the Riemann sphere $\sphere.$ Loewner energy of \(\eta\), denoted by $I^{L}(\eta)$, is by definition the Dirichlet energy of this driving term. It was shown by Wang in \cite{Wang} that if \(\eta\) passes through \(\infty\), then one has

\begin{equation}
I^{L}(\eta)=\frac{1}{\pi} \int_{\mathbb{H}}|\nabla \log | f^{\prime}(z)||^{2}\, d A_z+\frac{1}{\pi} \int_{\mathbb{H}^{*}}|\nabla \log | g^{\prime}(z)||^{2}\, d\mathrm{A}_z,
\end{equation}
where \(f\) and \(g\) map conformally the upper and lower half-planes $\mathbb{H}$ and $\mathbb{H}^{*}$ onto the two components of \(\sphere \backslash \eta\), while fixing \(\infty\) (i.e. $(f,g)$ is a welding pair). Furthermore it was  shown in \cite{Wang} that a Jordan curve has finite Loewner energy if and only if it is
a Weil-Petersson quasicircle, i.e. its normalized welding homeomorphism belongs to the
Weil-Petersson Teichmüller space.\\
There is an interesting connection between zeta-regularized determinants and Loewner energy, which was explored by Wang in \cite{wangzeta}. First let us briefly recall some basic facts about zeta function regularization of determinants.
To this end let \(\Delta_g\) be the Laplace-Beltrami operator with Dirichlet boundary
condition on a compact surface with conformal metric \((\mathscr{R}, g)\) with smooth (say at least \(C^{1,1}\)  boundary).
The zeta function associated to $-\Delta_g$ is defined by
\begin{equation}
    \zeta_{g}(s)=\sum_{j=1}^{\infty} \lambda_{j}^{-s}=\frac{1}{\Gamma(s)} \int_{0}^{\infty} t^{s-1} \sum_{j=1}^{\infty} e^{-t \lambda_{j}} d t
\end{equation}

where \(0<\lambda_{1} \leq \lambda_{2} \ldots\) is the discrete spectrum of \(-\Delta_g\). Now as a consequence of a result of H. Weyl \cite{Weyl}  one has that 
$\lambda_{j}\sim C\, j$ for a constant $C$, and \(\zeta_{g}\) is therefore holomorphic in \(\{s\in \mathbb{C};\,\operatorname{Re}(s)>1\}\). Thereafter \(
\zeta_{g}\) could be continued meromorphically to \(\mathbb{C}\) (with a simple pole at $s=1$), in particular \(\zeta_{g}\) is holomorphic at $s=0$  and hence \(\zeta_{g}^{\prime}(0)\) exists and is given by
$$
\zeta_{g}^{\prime}(0)=\lim _{s \rightarrow 0} \frac{\zeta_{g}(s)-\zeta_{g}(0)}{s-0}.
$$
More explicitly one has that
\begin{equation*}
    \zeta_{g}^{\prime}(s)=\sum_{j=1}^{\infty}\frac{-\log \left(\lambda_{j}\right)}{\lambda_{j}^{s}},
\end{equation*}
which for \(s=0\) yields that
$$
\zeta_{g}^{\prime}(0)=-\sum_{j=1}^{\infty} \log\left(\lambda_{j}\right)=-\log \left(\prod_{j=1}^{\infty} \lambda_{j}\right).
$$
Therefore we obtain
$$
e^{-\zeta_{g}^{\prime}(0)}=\prod_{j=1}^{\infty} \lambda_{j}=\operatorname{det}\left(-\Delta_g\right),
$$
which allows one to define the determinant of \(-\Delta_g\) by the formula
\begin{equation}
  \operatorname{det}\left(-\Delta_g\right):=e^{-\zeta_{g}^{\prime}(0)}.  
\end{equation}
When \((\mathscr{S}, g)\) is a compact surface without boundary, the Laplace-Beltrami operator \(\Delta_g\) has a discrete spectrum and a one-dimensional kernel, and its {\it regularized determinant} denoted by $\det_{\zeta}^{\prime}(-\Delta_g)$ is defined in a similar way as above
by considering only the non-zero eigenvalues.
An interesting fact in this context is that the zeta-regularized determinant of the Laplacian depends on both the Riemannian and the conformal structures of the surface.\\
Indeed if \(g=e^{2 \lambda} g_{0}\), with $\lambda\in C^{\infty}(\mathscr{R})$ (or in $C^{\infty}(\mathscr{S})$), is a metric
conformally equivalent to \(g_{0}\), then the relationship between the regularized determinants of $-\Delta_g$ and $-\Delta_{g_0}$ is given by the following \emph{Polyakov-Alvarez conformal anomaly formula} \cite{OSSarPhil}:
\begin{theorem}\label{PolyakovAlvarez}
   For a
compact surface \((\mathscr{S}, g)\) without boundary, one has that

\begin{equation}
\begin{aligned} \log \Big(\frac{\operatorname{det}_{\zeta}^{\prime}\left(-\Delta_{g}\right)}{ \operatorname{det}_{\zeta}^{\prime}\left(-\Delta_{g_0}\right)}\Big)=&-\frac{1}{6 \pi}\left[\frac{1}{2} \int_{\mathscr{S}}\left|\nabla_{g_0} \lambda\right|_{g_0}^{2}\, d\mathrm{A}_{g_0}+\int_{\mathscr{S}} K_{g_0} \lambda\, d\mathrm{A}_{g_0}\right] \\ &+\log\Big(\frac{ \operatorname{A}_{g}(\mathscr{S})}{\operatorname{A}_{g_0}(\mathscr{S})}\Big) 
\end{aligned}
\end{equation}
For a surface $(\mathscr{R},g)$ with boundary, the analogous formula is
\begin{equation}
\begin{aligned} 
\log \Big(\frac{\mathrm{det}_{\zeta}^{}\left(-\Delta_{g}\right)}{ \mathrm{det}_{\zeta}^{}\left(-\Delta_{g_0}\right)}\Big)=&-\frac{1}{6 \pi}\Big[\frac{1}{2} \int_{\mathscr{R}}\left|\nabla_{g_0} \lambda\right|_{g_0}^{2}\, d\mathrm{A}_{g_0}+\int_{\mathscr{R}} K_{g_0} \lambda\, d \mathrm{A}_{g_0}+\int_{\partial \mathscr{R}} k_{g_0} \lambda\, d \ell_{g_0}\Big]\\ & -\frac{1}{4 \pi} \int_{\partial \mathscr{R}} \partial_{n} \lambda\, d\ell_{g_0},
\end{aligned}
\end{equation}
where $d\mathrm{A}_{g_0}$ is the area measure, $d\ell_{g_0}$ is the boundary measure, $K_{g_0}$ is the Gauss curvature, and $k_{g_0}$ the geodesic curvature on the boundary, all associated to the metric $g_0$, and $\partial_{n}$ is the derivative along the unit outer normal. Moreover $\nabla_{g_0}$ is the Riemannian gradient, and $\mathrm{A}_{g_0}(\cdot)$ is the area of the surface $($both with respect to the metric $g_0$.$)$ 

\end{theorem}

In \cite{wangzeta} Wang  showed that for smooth simple loops, the Loewner energy can be
expressed in terms of the zeta-regularized determinants introduced above, of a certain Neumann jump operator. More explicitly, consider the 2-sphere  \(\sphere\) with a Riemannian metric $g$, and let \(\gamma \subset S^{2}\)
be a smooth Jordan curve dividing \(\sphere\) into two components \(D_{1}\) and \(D_{2}\). Denote
by \(\Delta_{D_{i}, g}\) the Laplacian with Dirichlet boundary condition on \(\left(D_{i}, g\right)\), $i=1,\,2$. Next define the functional \(\mathscr{H}(\cdot, g)\) on the space of smooth Jordan curves by
\begin{equation*}
\mathscr{H}(\gamma, g):=\log \operatorname{det}_{\zeta}^{\prime}\left(-\Delta_{\sphere, g}\right)-\log \operatorname{vol}_{g}\left(\sphere\right)  -\log \operatorname{det}_{\zeta}\left(-\Delta_{D_{1}, g}\right)-\log \operatorname{det}_{\zeta}\left(-\Delta_{D_{2}, g}\right).
\end{equation*}

Given the information above, Wang \cite{wangzeta} proved the following theorem   \begin{theorem}\label{wangdetsewing}
 If \(g=e^{2 \varphi} g_{0}\) is a metric conformally equivalent to the spherical
metric \(g_{0}\) on \(\sphere\), then:
\begin{enumerate}[font=\upshape]
    \item $\mathbb{S}^1$ minimizes \(\mathscr{H}(\cdot, g)\) among all smooth Jordan curves $\gamma\subset \sphere$ 
    \item Let \(\gamma\) be a smooth Jordan curve on \(\sphere\). We have the identity
\begin{equation}\label{wangloewner}
    I^{L}(\gamma)=12 \mathscr{H}(\gamma, g)-12 \mathscr{H}\left(\mathbb{S}^{1}, g\right)=
    12 \log \frac{\operatorname{det}_{\zeta}\left(-\Delta_{\mathbb{D}_{+}, g}\right) \operatorname{det}_{\zeta}\left(-\Delta_{\mathbb{D}_{-}, g}\right)}{\operatorname{det}_{\zeta}\left(-\Delta_{D_{1}, g}\right) \operatorname{det}_{\zeta}\left(-\Delta_{D_{2}, g}\right)},
\end{equation}

\end{enumerate}
 
where \(\mathbb{D}_{+}\) and \(\mathbb{D}_{-}\) are the two connected components of the complement
of \(\mathbb{S}^{1}\).
\end{theorem}

Moreover Wang showed that the Loewner
energy equals a multiple of the universal Liouville action of
Takhtajan-Teo given by \eqref{Liouville}.\\ Prior to \cite{wangzeta}, following the work by Schiffer \cite{Schifferplane,Schiffer_multiply}, Takhtajan-Teo \cite{Takhtajan_Teo_Memoirs} defined the Fredholm determinant of a quasicircle. {The Fredholm determinant of a curve has its origins in classical potential theory in association with Fredholm's solution of the Dirichlet and Neumann problems via the jump formula.  Schiffer \cite{Schifferplane,Schiffer_multiply,Schiffer_rend_Fredholm,Schiffer_Pol}, considering $C^3$ curves, showed that the classical Fredholm determinant could be expressed in terms of the Grunsky operator.  Takhtajan-Teo's work naturally extends this to quasicircles. }\\ Briefly, given a welding pair $(f,g)$ such that the corresponding operators $\mathbf{Gr}_f \mathbf{Gr}_{f}^{\ast}$ and $\mathbf{Gr}_g \mathbf{Gr}_g^{\ast}$ are trace-class, the Fredholm determinant of the corresponding quasicircle $\mathscr{C}={f}\left(\mathbb{S}^{1}\right)$ is defined by

$$\mathrm{Det}_{F}(\mathscr{C})= \det (\mathbf{I}-\mathbf{Gr}_f \mathbf{Gr}_{f}^{\ast})= \det (\mathbf{I}-\mathbf{Gr}_g \mathbf{Gr}_{g}^{\ast}).$$
Takhtajan and Teo then showed that the following sewing formula is valid for the determinants
\begin{theorem}\label{TTdetsewing}
Let \(\gamma={g}^{-1} \circ {f} \in T_{\mathrm{WP}}(\mathbb{D}_-)\) be in \(C^{3}(\mathbb{S}^1)\). Then for
\(\mathscr{C}={f}\left(\mathbb{S}^{1}\right)\) one has that
\begin{equation}\label{TTsewingformula}
    \operatorname{Det}_{F}(\mathcal{C})=\frac{\operatorname{det}_{\zeta}(-\Delta_{\tilde{\kap}}) \operatorname{det}_\zeta(- \Delta_{\tilde{\kap}^{*}})}{\operatorname{det}_{\zeta} (-\Delta_{\mathbb{D}_+}) \operatorname{det}_{\zeta}(-\Delta_{\mathbb{D}_{-}})}
.
\end{equation}

This could be interpreted as a sewing formula for the Laplace operator of a conformal metric
on the sphere \(\sphere\), which is the metric $|dw|^2$ on the interior domain \(\tilde{\kap}=\imath\left(\kap^{*}\right)\), with $\imath(z):=\frac{1}{z}$, and is the metric \(\frac{|dw|^2}{\left|\mathrm{f}^{-1}(w)\right|^{4}}\) on the exterior domain \(\tilde{\kap}^{*}=\imath(\kap)\). The Fredholm determinant \(\operatorname{Det}_{F}(\mathscr{C})\) is the inverse of the determinant of the Neumann jump operator which corresponds to cutting of \(\sphere\)
along the closed curve \(\mathscr{C}\) and considering Dirichlet boundary conditions for interior and exterior Laplace operators.
\end{theorem}

{Therefore putting Theorems \ref{wangdetsewing} and  \ref{TTdetsewing} together, one also sees a relation between the Loewner energy and the Fredholm determinant of a quaiscircle.}

\end{subsection}
\end{section}

\end{document}